\documentclass [11pt, twoside] {uwthesis}

\usepackage{amsmath}
\usepackage{amsthm}
\usepackage{amsfonts}
\usepackage{amssymb,latexsym}
\usepackage[all]{xy}
\usepackage{graphicx}
\usepackage{color}
\usepackage[mathscr]{eucal}
\usepackage{verbatim}

\theoremstyle{plain}
    \newtheorem{thm}{Theorem}[section]
    \newtheorem{prop}[thm]{Proposition}
    \newtheorem{lemma}[thm]{Lemma}
    \newtheorem{conj}[thm]{Conjecture}
    \newtheorem{cor}[thm]{Corollary}

\theoremstyle{definition}
    \newtheorem{defn}[thm]{Definition}

\theoremstyle{remark}
    \newtheorem{rem}[thm]{Remark}
    \newtheorem{example}[thm]{Example}

\numberwithin{section}{chapter}
\numberwithin{equation}{section}

\newtheorem*{namedtheorem}{\theoremname}
\newcommand{\theoremname}{testing}
\newenvironment{named}[1]{\renewcommand{\theoremname}{#1}
\begin{namedtheorem}}
{\end{namedtheorem}}

\newcommand{\A}{\mathcal A}
\newcommand{\B}{\mathcal B}
\newcommand{\C}{\mathcal C}
\newcommand{\D}{\mathcal D}
\newcommand{\E}{\mathcal E}
\newcommand{\F}{\mathcal F}

\newcommand{\T}{\mathcal T}
\newcommand{\W}{\mathcal W}
\newcommand{\V}{\mathcal V}

\newcommand{\FT}{\mathcal{F}_{tor}}
\newcommand{\stk}[1]{\stackrel{#1}{\longrightarrow}}
\newcommand{\la}{\ensuremath{\langle}}
\newcommand{\ra}{\ensuremath{\rangle}}

\newcommand{\ints}{\ensuremath{\mathbb{Z}}}

\newcommand{\rats}{\ensuremath{\mathbb{Q}}}

\newcommand{\ftwo}{\ensuremath{\mathbb{F}_2}}

\newcommand{\rar}{\ensuremath{\rightarrow}}

\newcommand{\n}{\noindent}

\newcommand{\lspace}{\hspace{2pt}}

\DeclareMathOperator{\Supp}{Supp}
\DeclareMathOperator{\Spec}{Spec}
\DeclareMathOperator{\Hom}{Hom}
\DeclareMathOperator{\RHom}{RHom}

\DeclareMathOperator{\proj}{proj}
\DeclareMathOperator{\Proj}{Proj}
\DeclareMathOperator{\stmod}{stmod}
\DeclareMathOperator{\StMod}{StMod}
\DeclareMathOperator{\Wide}{Wide}
\DeclareMathOperator{\Thick}{Thick}
\DeclareMathOperator{\Ker}{Ker}
\DeclareMathOperator{\im}{Im}
\DeclareMathOperator{\Stable}{Stable}
\DeclareMathOperator{\Cone}{Cone}
\DeclareMathOperator{\Sq}{Sq}
\DeclareMathOperator{\Ext}{Ext}

\DeclareMathOperator{\Pic}{Pic}

\DeclareMathOperator{\plainlim}{lim}
\newcommand{\clim}[1]{\ensuremath{\underset{#1}{\plainlim}\lspace}}

\DeclareMathOperator{\plainholim}{holim}
\newcommand{\holim}[1]{\ensuremath{\underset{\scriptscriptstyle #1}{\plainholim}\lspace}}

\DeclareMathOperator{\plaincolim}{colim}
\newcommand{\colim}[1]{\ensuremath{\underset{#1}{\plaincolim}\lspace}}

\newcommand{\rlim}[1]{\ensuremath{\underset{#1}{\text{Rlim}}\lspace}}

\DeclareMathOperator{\plainhocolim}{hocolim}
\newcommand{\hocolim}[1]{\ensuremath{\underset{\scriptscriptstyle #1}{\plainhocolim}\lspace}}

\newcommand{\sus}[1]{\ensuremath{\Sigma ^{#1}}}

\newcommand{\bp}[1]{\ensuremath{ \BP \langle {#1} \rangle}}
\newcommand{\BP}{\ensuremath{ B\!P}}

\begin{document}

%

\prelimpages

%
%

\Title{Refinements of chromatic towers and Krull-Schmidt decompositions in stable homotopy categories}
\Author{Sunil Kumar Chebolu}
\Year{2005}
\Program{Mathematics}
\titlepage

%
%





%
%
\setcounter{page}{-1}
\abstract{

 We study the triangulated subcategories of compact objects in stable homotopy categories such as the homotopy category of spectra, the derived categories 
of rings, and the stable module categories of Hopf algebras. In the first part of this thesis we use a $K$-theory recipe of Thomason to classify these 
subcategories. This recipe when applied to the category of finite $p$-local spectra gives a refinement of the ``chromatic tower''. This refinement has some interesting 
consequences. In particular, it gives new evidence to a conjecture of Frank Adams that the thick subcategory $\C_2$ can be generated
by iterated cofiberings of the Smith-Toda complex $V(1)$. 
Similarly by applying this $K$-theory recipe to derived categories, we obtain a complete classification of the triangulated subcategories
of perfect complexes over some noetherian rings.  Motivated by these classifications, in the second part of the thesis, we study Krull-Schmidt decompositions for thick 
subcategories. More precisely, we show that the thick subcategories of compact objects in the aforementioned stable homotopy categories decompose uniquely into 
indecomposable thick subcategories. Some 
consequences of these decompositions are also discussed. In particular, it is shown that all these decompositions respect $K$-theory. Finally in the last chapter we mimic 
some of these ideas in the category of $R$-modules. Here we consider abelian subcategories of $R$-modules that are closed under extensions and study their $K$-theory and 
decompositions.

} 
%
%
\tableofcontents

\acknowledgments{

I have had many fruitful discussions on this work with several mathematicians to whom I am 
greatly indebted; to my advisor John Palmieri more than anyone else. Among others, I thank
Steve Mitchell and Ethan Devinatz, for helping me understand several ideas in stable homotopy theory;  
Mark Hovey, for bringing to my attention  a paper of Frank Adams; and James Zhang, for being 
my encyclopedia of algebra.

Special thanks to Henning Krause and Srikanth Iyengar for their invaluable help and encouragement.
Chapter 5 was motivated from a small conversation with Henning Krause at a summer school in the University of Chicago. 
I had innumerable conversations with Srikanth Iyengar who helped me clarify my understanding
of several concepts in algebra. He provided some counterexamples when I was on the wrong track and also  
helped me improve a splitting result in section \ref{sec:KS-derivedcategory}. 

Finally I want to thank the mathematics department at the University of Washington for
providing me the McFarlan fellowship which came in invaluable during the final stages of this work.

}

\dedication{
\begin{center}
\vskip 40 mm
To my parents.
\end{center}
}

%
%

%
%

\textpages
 
 
\chapter{Introduction}
 
\section{Overview} In recent years there has been a lot of interest in the interactions between homotopy theory and its allied fields. Many fruitful interactions have been 
revealed  by the study of triangulated structures naturally  arising in  homotopy theory, homological algebra, algebraic geometry, and
modular representation theory. The far-reaching applications and generalisations of Brown representability \cite{ne5} \cite{ChKeNe} \cite{kr1},
Bousfield localisation \cite{Ri1} \cite{AlLeJe}, and nilpotence detection \cite{Ho} \cite{car}, to other fields are some striking examples of such interactions.
Several mathematicians brought to light many amazing a priori different theories 
by classifying various subcategories of triangulated categories arising in these various fields. Following the seminal work of Devinatz, 
Hopkins, and Smith \cite{dhs} in stable homotopy theory, this line of research was initiated 
by Hopkins in the 80s. In his famous 1987 paper \cite{Ho}, Hopkins classified the thick subcategories (triangulated categories that are closed under retractions) 
of the finite $p$-local spectra and those of perfect complexes over a noetherian ring. He showed that thick subcategories of the finite spectra are determined by the 
Morava $K$-theories and those of perfect complexes by the prime spectrum of the ring. These results have had tremendous impacts in their respective fields.
The thick subcategory theorem for finite spectra played a vital role in the study of nilpotence and periodicity.  For example, using this theorem
Hopkins and Smith \cite{hs} were able to settle the class-invariance conjecture of Ravenel \cite{rav} which classified the Bousfield classes of finite spectra. 
Similarly the thick subcategory theorem for the derived category  establishes a surprising connection between stable homotopy theory and algebraic 
geometry; using this theorem one is able to recover the spectrum of a ring from the homotopy structure of its derived category! These ideas were later pushed further 
into the world of derived categories of rings and schemes by Neeman \cite{Ne} and Thomason \cite{Th}, and into modular representation theory by Benson, Carlson and 
Rickard \cite{bcr}. Motivated by the work of Hopkins \cite{Ho}, Neeman \cite{Ne} classified the Bousfield classes and localising subcategories in the derived category of a 
noetherian ring. In modular representation theory, the Benson-Carlson-Rickard classification of the thick subcategories of stable modules over group algebras has led to 
some deep structural information on the representation theory of finite groups. Finally the birth of axiomatic stable homotopy theory \cite{mps} in the mid 90s encompassed 
all these various theories and ideas and studied them all in a more general framework. With all these developments over the last 30 years, the importance of triangulated 
categories in modern mathematics is by now abundantly clear.

Motivated by these results, in this thesis we focus on the study of triangulated subcategories in stable homotopy categories. (Loosely speaking, a stable 
homotopy category is a sufficiently well-behaved triangulated category that is formally similar to the homotopy category of spectra, e.g., the derived category of 
a ring; for more details, see section \ref{sec:tc}) We focus on classifications of triangulated subcategories and decompositions of thick subcategories of compact 
objects in stable homotopy categories which include the homotopy category of spectra, derived categories of rings, and the stable module categories of some finite 
dimensional co-commutative Hopf algebras. Our main tools will be the thick subcategory theorems and $K$-theory.

\section{Main Results} 

One of our goals is to classify the triangulated subcategories of finite spectra and those of perfect complexes over a noetherian ring. To this end,
we adopt a $K$-theoretic approach of Thomason. To set the stage, let $\T$ denote a triangulated category (we have the compact objects of a stable homotopy category
in mind) and let $K_0(\T)$ denote the \emph{Grothendieck group} of $\T$. Thomason
established in \cite{Th} a bijection between the dense (see definition \ref{defn:dense}) triangulated subcategories of $\T$ and the subgroups of $K_0(\T)$. 
This theorem gives a key to the problem of classifying triangulated subcategories of $\T$:
\emph{First classify all thick subcategories of $\T$ and then compute the Grothendieck groups of these thick subcategories.
For every subgroup of each of these Grothendieck groups we get a triangulated subcategory of $\T$, and every triangulated subcategory
arises this way.} Thus we can try to refine all the aforementioned thick subcategory theorems (\cite{hs}, \cite{Ho,Ne}, \cite{bcr}) to obtain a classification of 
triangulated subcategories. By  ``chromatic tower for a stable homotopy category'', we mean the lattice of all the thick subcategories of compact objects. Since a 
thick subcategory is also a triangulated subcategory, this $K$-theoretic approach gives a refinement of the chromatic towers, which explains the first part of the 
title of this thesis.

In chapter 3 we apply this $K$-theoretic technique to the category $\F_p$ of finite $p$-local spectra. 
This gives rise to the following triangulated subcategories of $\F_p$. Let $\C_n$ denote the thick subcategory of $\F_p$ consisting of $K(n-1)$-acyclics, and 
for non-negative integers $k$ and $n \ge 1$, define a full subcategory of $\F_p$ as
\[ \C_n^k = \{ X \in \C_n \; | \; \chi_n(X) := \sum_i (-1)^i \log_p |\bp{n-1}_i X| \equiv 0 \ \text{mod} \, l_n k \},\]
where $l_n$ is the smallest positive value of $\chi_n(X)$ on $\C_n$.
Define $\C_0^k$ to be the full subcategory of spectra in $\F_p$ whose (rational) Euler characteristic is divisible by $k$.
Our understanding of the lattice of triangulated subcategories of $\F_p$ can be summarised in the following theorem.

\begin{named}{\textbf{Theorem A}} \emph{For each $n \ge 0$ and $k \ne 1$, the subcategory $\C_n^k$ is a dense triangulated subcategory of $\C_n$. 
These triangulated subcategories satisfy the inclusions
\[ \C_{n+1} \subsetneq \C_n^k \subsetneq \C_n .\] 
Moreover, when $n=0$ or $1$, every dense triangulated subcategory of $\C_n$ is equal to $\C_n^k$ for some $k$. }
\end{named}

One of the interesting parts of this theorem is to establish the inclusion $\C_{n+1} \subseteq \C_n^k$. This is shown using a spectral sequence calculation.
For a fixed $X \in \C_{n+1}$, we construct a strongly convergent Bockstein spectral sequence 
\[ E_r^{*,*} \Rightarrow E(n)_*(X)\]
whose $E_1$ term is build out of $\bp{n-1}_* X$, and which abuts to 
\[ \holim{v_n} \bp{n}_* X = E(n)_*(X) \quad (= 0 \; \mbox{by assumption}).\]
Now working backward with this spectral sequence we are able to show that $X \in \C_n^k$ for all $k$.
This theorem gives some new evidence to the following conjecture of Adams.

\begin{conj} \cite[Page 529]{MayTho} The Smith-Toda complex $V(1)$ generates the thick subcategory $\C_2$ by iterated cofiberings.
\end{conj}
It is not hard to see that this conjecture implies that $\C_2 \subseteq \C_1^k$, which is part of the above theorem.

Among other things, the above theorem is saying that the thick subcategory $\C_1$ can be generated by iterated cofiberings of $M(p)$, the mod-$p$ Moore
spectrum; if we also allow retractions, then this is a well-known consequence of the thick subcategory theorem of Hopkins and Smith \cite{hs}.
This gives the following interesting corollary. For a $p$-torsion spectrum $X$, define its Euler characteristic by
\[\chi_1(X): = \sum_i (-1)^i \log_p |H\mathbb{Z}_i (X)|. \]

\begin{cor} Let $X$ be a type-$1$ spectrum and $Y$ an arbitrary finite $p$-torsion spectrum. Then $Y$ can be generated by $X$ via cofibrations
if and only if $\chi_1(X)$ divides $\chi_1(Y).$
\end{cor}

This corollary can be compared with the following result of Hopkins \cite{Ho}:  Let $X$ and $Y$ be finite $p$-local spectra. Then $Y$ can be generated by 
$X$ via cofibrations and retractions if and only $\Supp(Y) \subseteq \Supp(X)$. ($\Supp(X)$ is the chromatic support of $X$;  the set of all
non-negative integers $n$ such that $K(n)_* X \ne 0$.)


In chapter 4, we do similar $K$-theoretic analysis in the derived category and obtain complete classifications of perfect complexes over some noetherian rings 
(PIDs and Artin rings). See section \ref{se:classifications} for these results.\vskip 3mm \noindent

We now move on to the second part of this thesis -- Krull-Schmidt decompositions for thick subcategories; see definition \ref{defn:KS}. Our motivation for studying such 
decompositions comes partly from the aforementioned classifications of triangulated subcategories. We have noticed  several situations where the Grothendieck group 
under study is a direct sum of infinite cyclic  groups. A careful investigation has led us to the observation that this behaviour is in fact induced by a decomposition 
of the underlying category. More precisely, if $\A = \coprod \A_i$ is a Krull-Schmidt decomposition of a thick subcategory $\A$, then under certain conditions (which are 
always satisfied in the examples we consider, 
e.g., $\Hom(\A_i, \, \A_j) = 0$ for all $i \ne j$),  $K_0(\A) \cong \bigoplus_i K_0(\A_i)$. In short, Krull-Schmidt decompositions respect $K$-theory. 

We now summarise our Krull-Schmidt results in the next theorem.

\begin{named} {\textbf{Theorem B}} \emph{Let $R$ denote a noetherian ring and let $B$ denote a finite dimensional graded co-commutative Hopf algebra satisfying the 
tensor product property (e.g., finite dimensional sub-Hopf algebras of the mod-$2$ Steenrod algebra). Then we have the following.
\begin{enumerate}
\item Every thick subcategory of perfect complexes over $R$ admits a Krull-Schmidt decomposition.
\item Every thick  ideal of small objects in the chain homotopy category of projective $B$-modules is indecomposable.
\item Every thick ideal of finite dimensional stable $B$-modules admits a Krull-Schmidt decomposition.
\end{enumerate}
Further all these Krull-Schmidt decompositions are unique and they respect $K$-theory. }
\end{named}

The essential idea behind the proofs of the decompositions in the above theorem is quite simple. We use thick subcategory theorems which establish a bijection 
between thick subcategories of small objects and some geometric spaces. We obtain the above decompositions by decomposing the corresponding geometric spaces in an
appropriate way. These decompositions give some interesting splitting results. For example, 

\begin{cor}   Let $X$ be a perfect complex over a noetherian ring $R$. Then $X$ admits a unique decomposition into perfect complexes,
\[ X \cong \underset{i \in I}{\bigoplus}\; X_i,\]
such that  the supports of the $X_i$ are pairwise disjoint and indecomposable.
\end{cor}

We now move on to the last chapter. Since computing the Grothendieck groups of thick subcategories is the key to the problem of classifying triangulated subcategories, 
it is natural to investigate different ways to compute these groups. 
So the motivating question now is the following. Can one compute the Grothendieck groups of thick subcategories of perfect complexes in the 
$R$-module category? More precisely, if $\C$ is a thick subcategory of perfect complexes over $R$, is there a subcategory of
$R$-modules whose Grothendieck group is naturally isomorphic to that of $\C$? 

Towards this, it is natural to investigate the right analogue of a thick subcategory of complexes
in the world of modules. This has been investigated by Mark Hovey \cite{wide}. It turns out that the right 
notion of ``thickness for modules'' is ``wideness'':  A subcategory $\C$ of an abelian category
is said to be \emph{wide} if it is an abelian subcategory that is closed under extensions. 

Hovey \cite{wide} showed that when $R$ is regular and coherent, the wide subcategories of $\Wide (R)$ (the wide subcategory generated by $R$) 
are in bijection with the thick subcategories of $\Thick (R)$ (the thick subcategory generated by $R$).
More precisely, he showed that the map 
            \[f: L_{\Wide} (R) \rightarrow L_{\Thick}(R),\]
from the lattice of wide subcategories to that of thick subcategories defined by 
\[f(\C) = \{ X \in \Thick (R) : H_n(X) \in \C \;\; \text{for all}\; \; n \},\]
is an isomorphism when $R$ is regular and coherent. We use this classification result of Hovey to prove the following theorem.

\begin{named} {\textbf{Theorem C}} \emph{Let $R$ be a regular coherent ring. If $\C$ is any wide subcategory of $\Wide(R)$, then
\[ K_0(\C) \cong K_0(f(\C)). \]
In addition if $R$ is noetherian, then every wide subcategory $\mathcal{W}$ of finitely generated modules
admits a  Krull-Schmidt decomposition: $\mathcal{W} = \coprod_{i \in I} \mathcal{W}_i$. Moreover,
$K_0(\mathcal{W}) \cong \underset{i \in I}{\bigoplus} K_0(\mathcal{W}_i).$ 
}
\end{named}

As a corollary to this theorem we recover a known splitting result for finitely generated modules over a noetherian ring; see corollary \ref{cor:splittingformodules}.

\section{Organisation}

This thesis is organised as follows. In chapter 2 we set up the categorical stage for our work. Here we recall some basic 
facts about triangulated categories and Grothendieck groups, and  state the main $K$-theory recipe of Thomason. In chapters 3 and 4, we apply these tools to the study of 
classifications of triangulated subcategories of finite spectra and perfect complexes respectively. Chapter 5 deals with Krull-Schmidt decompositions in these various 
stable homotopy categories. Finally in the last chapter we study wide subcategories where we mimic some of these ideas (from the previous chapters) in the 
$R$-module category. 

Chapters 3 and 4 can be read independent of each other. Although most of Chapter 5 can be read separately, it would be well motivated if read after chapters 3 and 4. 
Chapter 6 takes some input from chapter 4 and very little from chapter 5. The following picture shows the logical interdependencies of the chapters. 
(A dotted arrow indicates a mild dependency.)

\begin{figure}[!h]
\[
\xymatrix{
                            &  *++[o][F-]{1} \ar[d]                      &                     \\
                            &  *++[o][F-]{2} \ar[dl] \ar[dr] \ar[dd]     &                     \\
 *++[o][F-]{3} \ar@{..>}[dr]      &                                      &               *++[o][F-]{4} \ar@{..>}[dl] \ar[dd] \\
                            &  *++[o][F-] {5} \ar@{..>}[dr]              &                     \\
                            &                                            &   *++[o][F-]{6}
}
\]
\label{logic}
\caption{Interdependencies of chapters}
\end{figure}
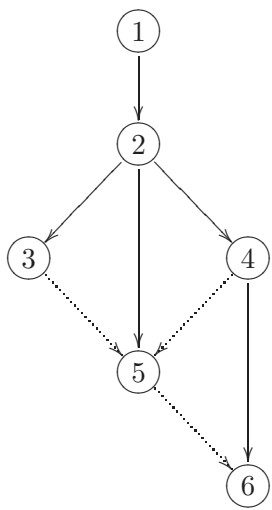

Every chapter begins with an agenda and ends (except chapter 2) with a discussion of questions and further directions.


 
\chapter{Triangulated categories and $K$-theory}

In this chapter we will set up the categorical stage for our work. In section \ref{sec:tc} we review some background material related to triangulated
categories and introduce the  triangulated categories that we will be working in the later chapters. 
Sections \ref{sec:groups} and \ref{sec:recipe} are the most important sections of this chapter. In these sections we introduce a
$K$-theory technique of Thomason which will be used for classifying triangulated subcategories. In the last section we consider triangulated categories which
are equipped with a nice product and observe that the Grothendieck groups of such categories naturally admit 
the structure of a commutative ring and explore some consequences of this observation.

\section{Triangulated categories} \label{sec:tc}
Roughly speaking, triangulated categories are those 
categories where the categorical notions like kernel, cokernel, limit etc. exist only in some weak sense. Typically they arise when one passes to homotopy by
inverting some class of maps. Derived categories of rings (obtained by inverting the quasi-isomorphisms) are good examples of triangulated categories.
These categories were introduced by Puppe \cite{puppe} about 40 years ago and their importance is becoming more and more visible all over mathematics,
more so in homotopy theory and algebra. Formally, a triangulated category is an additive category that is equipped with a suspension endofunctor $\Sigma$
which is an equivalence, and a collection of sequences $(A \rar B \rar C \rar \Sigma A)$ called exact triangles (a.k.a. cofibre sequences) which satisfy some axioms;
see \cite[Appendix A]{mps}.

\begin{example} The following are the main examples of triangulated categories that we are interested in.
These will be discussed in more detail in the subsequent chapters.
\begin{enumerate}
\item $\mathcal{S}$ -- The stable homotopy category of $p$-local spectra.
\item $D(R)$ -- The unbounded derived category of $R$-modules.
\item $K(\Proj\,B)$ -- The chain homotopy category of unbounded chain complexes of graded projective $B$-modules over a finite dimensional co-commutative Hopf algebra
$B$.
\item $\StMod(B)$ -- The stable module category over a finite dimensional co-commutative Hopf algebra $B$.
\end{enumerate}
\end{example}

Stable homotopy theory, following Hovey-Palmieri-Strickland \cite{mps}, is the study of \emph{``sufficiently well-behaved''} 
triangulated categories. 
All the categories mentioned in the above example are such categories.
Loosely speaking,  a triangulated category is  sufficiently well-behaved if it is formally similar to the homotopy category of spectra.
We will make this a little more precise below but we want to emphasise that this will be our view point. These sufficiently well-behaved
triangulated categories  will be called stable homotopy categories and our goal is to understand more about the structure of such 
categories by analysing their subcategories.

A sufficiently well-behaved triangulated category (i.e., stable homotopy category) has the following extra structure and properties.

\begin{defn} \cite{mps} A triangulated category $\T$ is a \emph{tensor triangulated category} if $\T$ is equipped with a triangulated 
functor (smash product) 
\[ \wedge:\, \T \times \T \rar \T \]
that is covariant in both variables and which is unital ($\exists \,S$  such that $S \wedge X \cong X$ for all 
objects $X$), symmetric ($A \wedge B \cong B \wedge A$), and associative ($A \wedge (B \wedge C) \cong (A \wedge B) \wedge C$). 
In addition, if $\T$ has arbitrary coproducts, a  closed symmetric  monoidal structure (\cite[Appendix A]{mps}) that is compatible with the 
triangulation, and satisfies some technical conditions like the existence of strongly dualisable objects and Brown representability for cohomology
functors, then we  say that $\T$ is a \emph{stable homotopy category}. 
(These later conditions are not so important for us right now.)
\end{defn}

The following table shows the smash products and  unit objects for the categories listed in the above example.

\begin{table}
\begin{center}
\caption{Unit objects and smash products in stable homotopy categories}
\begin{tabular}{c|l|l}
Category &  \hspace{20 mm} Product & \hspace{12 mm} Unit object  \\
\hline
$\mathcal{S}$  &  \;\;$\wedge$        \hspace{4.6 mm} smash product of spectra          &  \;\; $S$ (sphere spectrum)\\
$D(R)$         &  $^L\otimes_R$  \hspace{1.5 mm} derived tensor product                 &  \;\; $R$   (concentrated in degree 0) \\
$K(\Proj\,B)$   &  \;$\otimes_k$    \hspace{3 mm} tensor product of chain complexes        &  \;\; projective resolution of $k$ \\
               &                 \hspace{10 mm} with diagonal $B$-action.                &          \\
$\StMod(B)$    &  \;$\otimes_k$    \hspace{3 mm} tensor product of vector spaces          &  \;\; $k$    (trivial module)  \\
               &                 \hspace{10 mm} with  diagonal $B$-action             & 
\end{tabular}
\end{center}
\end{table}

\begin{defn} A  full subcategory $\A$ of a triangulated category $\T$  is a \emph{triangulated subcategory} if it is an additive subcategory that is
closed under isomorphisms, suspensions $(\Sigma^n, n \in \mathbb{Z})$, and cofibrations $(A \stk{f} B \rar \Cone{f})$.  The triangulated subcategory
generated by a set $X$ of objects in $\T$ is the smallest triangulated subcategory  of $\T$ that contains $X$. An additive functor between
triangulated categories is said to be an \emph{exact functor} if it commutes with $\Sigma$ and preserves exact triangles.
\end{defn}

Unless stated otherwise, all our subcategories will be full and closed under isomorphisms.

\section{Grothendieck groups} \label{sec:groups}

We begin by recalling some definitions and results from
\cite{Th}. Let $\T$ denote a  triangulated category that is
essentially small (i.e., it has only a set of isomorphism classes of
objects). Then the \emph{Grothendieck group} $ K_0(\T) $ is defined to
be the free abelian group on the isomorphism classes of $\T$ modulo
the Euler relations $[B]=[A]+[C]$, whenever $A\rightarrow B
\rightarrow C \rightarrow \Sigma A $ is an exact triangle in $\T$
(here $[X]$ denotes the element in the Grothendieck group  that is
represented by the isomorphism class of the object $X$). This is
clearly an abelian group with $[0]$ as the identity element and
$[\Sigma X]$ as the inverse of $[X]$. The identity
$[A]+[B]=[A \amalg B]$ holds in the Grothendieck group. Also note that any
element of $K_0(\T)$ is of the form $[X]$ for some $X\in \T$. All
these facts follow easily from the axioms for a triangulated category. 

Grothendieck groups have the following universal property: Any map $\sigma$ from the set of  isomorphism
classes of $\T$  to an abelian group $G$ such that the Euler relations hold in $G$
factors through a unique  homomorphism $ f: K_0(\T)
\rightarrow G $; see the diagram below.
\[
\xymatrix{ {\{\mbox{isomorphism classes of objects in}\, \T\}} \ar[d]^\pi \ar[dr]^
\sigma \\  {K_0(\T)} \ar@{..>}[r]^f & G }
\]
Therefore the map $\pi$ will be called the \emph{universal Euler characteristic function}.
$K_0(-)$ is clearly a covariant functor from the category of small triangulated categories to the 
category of abelian groups.

Note that these Grothendieck groups are defined only for essentially
small triangulated categories. For any ring $R$, the bounded derived category of $R$-modules is
essentially small, and so is the category of finite spectra. The unbounded derived category of $R$-modules, however, is not essentially
small.

The following important lemma due to Landsburg \cite{la}  will be used several
times. This is a very nice criterion for the equality of two
classes in the Grothendieck group.

\begin{lemma} \cite{la} \label{le:landsburg}  Let $\T$ be an essentially small triangulated category. If
$X$ and $Y$ are objects in $\T$, then $[X]=[Y]$ in $K_0(\T)$ if
and only if there are objects $A, B$, $C$ in  $\T$ and maps such
that there are exact triangles
\[A \stk{f}  B \amalg X \stk{g} C \stk{h} \Sigma A, \]
\[A \stk{f'} B \amalg Y \stk{g'} C \stk{h'} \Sigma A .\]
\end{lemma}

\section{Thomason's $K$-theory recipe} \label{sec:recipe}

In this section, we introduce a $K$-theory recipe due to Thomason that will be  central to chapters 3 and 4. We begin with some important 
categorical definitions.

\begin{defn} \label{defn:thick}
A triangulated subcategory $\C$ of $\T$ is said to be a \emph{thick subcategory} if it is closed under
retractions, i.e., given a commuting diagram
\[
\xymatrix{
 B \ar[r]^i \ar@/_20pt/[rr]_= & A \ar[r]^r & B
}
\]
in $\T$ such that $A$ is an object of $\C$, then so is $B$. Since retractions split in a triangulated category, 
it is easily seen that this property of $\C$ is equivalent to saying that $\C$ is closed under 
direct summands: $A \amalg B \in \C  \Rightarrow  A \in \C$ and 
$B \in \C$. 
\end{defn}

\begin{defn} \label{defn:small}
An object $X$ in a triangulated category $\T$ is \emph{small} if the natural map 
\[ \underset{\alpha \in \Lambda}{\bigoplus}\, \Hom(X, A_{\alpha}) \rightarrow  \Hom(X, \,\underset{\alpha \in \Lambda}{\coprod} A_{\alpha})\]
is an isomorphism for all set indexed collections of objects $A_{\alpha}$ in $\T$. 
(Some authors call such objects as finite or compact objects.)
\end{defn}

\begin{example} The full subcategory of small objects in any triangulated category is thick. 
\end{example}

\begin{defn} \label{defn:dense}
We say that a triangulated subcategory $\C$ is \emph{dense} in $\T$ if every object in $\T$ is a direct 
summand of some object in $\C$. 
\end{defn}

So in some sense, to be made precise below, a dense triangulated subcategory is ``very big'' relative to
the ambient category. To make this more precise, we first need the concept of a quotient triangulated category.

\begin{thm} \cite[Theorem 2.1.8]{ne1} Let $\D$ be a triangulated category and $\C \subseteq \D$ a triangulated subcategory. 
Then there exists a triangulated category $\D/\C$ and a triangulated
functor $F : \D \rightarrow \D/\C$ with $\C \subseteq \Ker(F)$  such
that the following universal  property is satisfied. If $G: \D
\rightarrow \T$ is any other triangulated functor with $\C \subset
\Ker(G)$, then it factors uniquely as
\[  \D \stk{F} \D/\C \rightarrow \T .\]

Further the kernel of the universal functor $F : \D \rightarrow \D/\C$
is precisely the thick closure of $\C$ in $\D$ (i.e., intersection of all thick subcategories  in $\D$ that contain $\C$).
\end{thm}

\begin{prop}  Let $\D$ be a triangulated category and $\C \subseteq \D$ be a 
triangulated subcategory. Then $\C$ is dense in $\D$ if and
only if $\D/\C = 0$.
\end{prop}

\begin{proof} It follows directly from the definitions that $\C$ is dense in
$\D$ if and only if $\D$ is the thick closure of $\C$; see \cite[remark 1.5]{Th}. Now from the above theorem, we note that this happens if and only if the
kernel of $\D \rightarrow \D/\C$ is equal to  $\D$. The latter is clearly
equivalent to saying that $\D/\C = 0 $. 
\end{proof} 

The following theorem due to Thomason \cite{Th} is the foundational theorem
that motivated this thesis.

\begin{thm} \cite[Theorem 2.1]{Th} \label{main} Let $\T$ be an essentially small triangulated category. 
Then there is a natural order preserving bijection between the posets,

\begin{center}
\{dense triangulated subcategories $\A$ of $\T$ \}
$\overset{f} {\underset{g}{\rightleftarrows}}$  \{subgroups $H$ of
$K_0(\T)$ \}.
\end{center}

\noindent The map $f$ sends $\A$ to image of the map $K_0(\A) \rar K_0(\T)$, 
and the map $g$ sends $H$ to the full subcategory of all objects $X$
in $\T$ such that $[X] \in H$.
\end{thm}

\subsection{Thomason's $K$-theory recipe} \label{ss:recipe}

The importance of Thomason's theorem \ref{main} can be seen from the following simple observation.
\emph{Every triangulated subcategory $\mathcal{A}$ of $\T$ is dense in a unique thick
subcategory of $\T$} --  the one obtained by taking the intersection of all the thick subcategories 
of  $\T$ that contain $\mathcal{A}$. This observation in conjunction with theorem \ref{main} gives the following brilliant 
recipe of Thomason to the problem of classifying the triangulated subcategories of $\T$:

\begin{enumerate}
\item{Classify the thick subcategories of $\T$.}
\item{Compute the Grothendieck groups of all thick subcategories.}
\item{Apply Thomason's theorem \ref{main} to each thick subcategory of $\T$.}
\end{enumerate}

We will apply this recipe to the categories of small objects in some stable homotopy categories like  
 the stable homotopy category of spectra (Chapter 3) and the derived categories of rings (Chapter 4). 
The stable homotopy structure on these categories will guide us while applying this recipe and 
consequently we will derive some structural information on these categories.

\section{Grothendieck rings} \label{s:gring}

\noindent Through out this section $\T$ will denote a tensor triangulated 
category that is essentially small. Now making use of the available smash
product, we want to define a ring structure on the Grothendieck group. This can be done in a very natural and obvious way.

\begin{defn} If $[A]$ and $[B]$
are any two elements of $K_0(\T)$, then define $[A] [B] := [A \wedge B] .$
\end{defn}

\begin{lemma} This is a well-defined
operation and makes $K_0(\T)$ into a commutative ring.
\end{lemma}

\begin{proof} We show that the above product is well-defined and
leave the remaining straightforward details to the
reader. Suppose $[A]=[A']$ and $[B]=[B']$. We have to show that
$[A \wedge A']=[B \wedge B'] $. Since $[A]=[A']$, Landsburg's criterion
gives us objects $X,Y,Z$, and maps such that the following 
sequences are exact triangles.
\noindent
\[X \rightarrow Y \amalg A \rightarrow Z \rightarrow \Sigma X \]
\[X \rightarrow Y \amalg A' \rightarrow Z \rightarrow \Sigma X .\]

\noindent
Similarly using $[B]=[B']$, we get exact triangles (for some objects $U,V,W$ and some maps),

\noindent
\[U \rightarrow V \amalg B \rightarrow W \rightarrow \Sigma U \]
\[U \rightarrow V \amalg B' \rightarrow W \rightarrow \Sigma U. \]

Now smashing the first two triangles with $B$ tells us that $[A \wedge B]=[A' \wedge B]$, and smashing the last two triangles with $A'$  tells us that
$[A' \wedge B]=[A' \wedge B']$. By combining these two equations we get $[A \wedge B] = [A'\wedge B']$.  
Commutativity of the ring follows from the equivalence given by the twist map $(A \wedge B \cong B \wedge A)$ and the Grothendieck class
$[S]$ of the unit object $S$ will be the identity element for multiplication.
\end{proof}

The following definitions are motivated by their analogues in commutative ring theory.

\begin{defn} A full triangulated
subcategory $\C$ of $\T$ is said to be \emph{$\otimes$-closed} or a
\emph{triangulated ideal} if for all $A \in \T$ and for all $B \in \C $, $B
\wedge A \in \C$. A full triangulated ideal is said to be respectively
thick or dense if it is such as a triangulated subcategory.
\end{defn}

\begin{defn} A full triangulated ideal $\B$ of $\T$ is said to be \emph{prime} if
for all $X$, $Y \in \T$ such that  $X \wedge Y \in \B$,  either $X
\in \B$ or $ Y \in \B$. Similarly $\B$ is said to be  a \emph{maximal} in $\T$ if there
is no triangulated ideal $\A$ such that $\B \subsetneq \A \subsetneq \T$.
\end{defn}

Before we can state the main result of this section, we need to recall
another lemma from \cite{Th}.

\begin{lemma} \cite[Lemma 2.2]{Th} Let $\A$ be a dense triangulated subcategory of an essentially small triangulated category
$\T$. Then for any object $X$ in $\T$, one has that $X \in \A$ if and
only if $[X] = 0 \in K_0(\T)/ \im (K_0(\A) \rar K_0(\T)).$
\end{lemma}

The following theorem is now expected.

\begin{thm} \label{imain}
Let $\T$ be a tensor triangulated category that is essentially
small. Then, under Thomason's bijection (theorem \ref{main})
\begin{center}
\{dense triangulated subcategories of $\T$ \} $\longleftrightarrow$ \{subgroups of $K_0(\T)$\}, 
\end{center}
we have the following correspondence.
\begin{enumerate} 
\item The dense triangulated ideals correspond precisely to the ideals of
the ring $K_0(\T)$.
\item The dense prime triangulated ideals correspond precisely to the
prime ideals of $K_0(\T)$.
\item The dense maximal triangulated ideals correspond precisely to
the maximal ideals of $K_0(\T)$.
\end{enumerate}
\end{thm}
\begin{proof} Except possibly the second statement, everything else is straightforward.

\noindent 1. Let $I$ be any ideal in $K_0(\T)$. Then the
corresponding dense triangulated category is  $\T_{I} = \{X \in \T :[X] \in I \}$.  
Now for $A \in \T$ and $B \in A_I$, note that $[B \wedge A]=[B][A] \in I$ ($I$ being an ideal) and
therefore $B \wedge A \in \T_{I}$. This shows that $\T_{I}$ is a triangulated ideal.
The other direction is equally easy.

\noindent 2. Suppose $H$ is a prime ideal in $K_0(\T)$.  The
corresponding dense triangulated ideal is given by $\B = \{ X : [X] \in H \}$. 
Now if $A \wedge B \in  \B$, then by definition of $\B$, we have
$[A\wedge B] \in H $, or equivalently $[A][B] \in H $. Now primality of
$H$ implies that either $[A] \in H$ or $[B] \in H$, which means 
either $A \in \B $ or $ B \in \B$. For the other direction, suppose
$\B$ is dense prime triangulated ideal of $\T$. The
corresponding subgroup $H$ is the image of the map $K_0(\B) \rar K_0(\T) $. 
Now suppose the product of two elements $[A]$ and $[B]$ belongs to $H$. Then we have $[A \wedge B] \in H$ or
equivalently  $[A \wedge B] = 0 $ in $K_0(\T)/ H$.  By the above lemma, we then have
 $A \wedge B \in \B$.  Since $\B$ is prime,
this implies that either $A \in \B$ or $B \in \B$, or equivalently $[A] \in H$ or $[B] \in H$.

\noindent 3. This follows directly from the fact that the bijection is
a map of posets. 
\end{proof}


 
\chapter{Triangulated subcategories of finite spectra}

In this chapter we work in the stable homotopy category $\mathcal{S}$ of spectra. There are many good descriptions of this category; see the books by
Adams \cite{bluebook} and Margolis \cite{mar}. The objects in this category are called spectra. We work mostly $p$-locally, so all our spectra are localised at  
a fixed prime $p$. The category $\mathcal{S}$ has a symmetric and associative product known as the smash product $(\wedge)$ which is compatible with the triangulation. 
In this product structure the sphere spectrum $S$ plays the role of the unit object. The small objects in this category are precisely the $p$-localisations 
of (possibly desuspended) suspension spectra of finite CW-complexes. The full subcategory of these finite objects will be denoted by $\F_p$.

The motivating problem for this chapter is to classify the triangulated subcategories of $\F_p$ using Thomason's recipe. We begin by recalling the
celebrated thick subcategory theorem of Hopkins-Smith in section \ref{se:hsthick}. In section \ref{se:K_0groups} we analyse the Grothendieck groups of 
the thick subcategories of $\F_p$ using some $\BP$-based homology theories and discuss a conjecture of Adams.
In section \ref{se:lattice} we study the lattice of triangulated subcategories of $\F_p$ and give some evidence for Adams's  conjecture. 
In the last section we classify the dense triangulated subcategories of finite torsion spectra.

\section{Hopkins-Smith thick subcategory theorem} \label{se:hsthick}

Recall that for each non-negative integer $n$, there is a field 
spectrum called the $n$th Morava $K$-theory $K(n)$, whose coefficient ring is
$\mathbb{F}_p[v_n,v_n^{-1}]$ with $|v_n| = 2(p^n-1)$. The celebrated thick subcategory theorem of 
Hopkins and Smith then states:

\begin{thm} \cite{hs} For each non-negative integer $n$,
let $\C_n$ denote the full subcategory of all finite $p$-local
spectra  that are $K(n-1)$-acyclic. Then a non-trivial subcategory $\C$ of $\F_p$ is 
thick if and only if $\C = \C_n$ for some $n$. Further these thick subcategories
give a nested decreasing filtration of $\F_p$ \cite{rav, mit}: 
\[ \cdots \C_{n+1} \subsetneq \C_n \subsetneq \C_{n-1} \cdots \subsetneq \C_1 \subsetneq \C_0 (= \F_p).\]
\end{thm}

A property $P$ of finite spectra is said to be \emph{generic} if the collection of all spectra in $\F_p$ which satisfy the given property $P$ is $\C_n$ for some $n$.
A spectrum $X$ is said to be of type-$n$ if $X$ belongs to $\C_n - \C_{n+1}$. For example, the sphere spectrum $S$ is of type-$0$, the mod-$p$ Moore spectrum
$M(p)$ is of type-$1$, etc.

\noindent
\section{Grothendieck groups for thick subcategories} \label{se:K_0groups}

The problem of computing the Grothendieck groups of thick subcategories of the
finite $p$-local spectra was first considered, to our knowledge, by Frank Adams.
This appeared in an unpublished manuscript \cite[Page 528-529]{MayTho} of Adams on the work of Hopkins.
We begin with a recapitulation of Adams's work and then bring it to this new context of classifying
triangulated subcategories. We should also point out that our Euler characteristic functions are simpler
than the ones considered by Adams.

\subsection{$\C_0$ - Finite $p$-local spectra} 
We begin with the fundamental notion of Euler characteristic of a finite spectrum. Recall that if $X$ is any finite $p$-local spectrum, 
the Euler characteristic of $X$ with rational coefficients  is given by  
\begin{equation} \label{E:chi_0}
\chi_0 (X) = \sum_{i= -\infty}^{\infty} (-1)^i \dim_{\rats} \;H\rats_i(X).
\end{equation}

Since $X$ is a finite spectra it has homology concentrated only in some finite range and therefore this is a
well-defined function.

\begin{example} For every non-negative integer $m$, define a full subcategory $\C_0^m$ of $\C_0 \, (= \F_p)$ as:
 \[\C_0^m =\{ X \in \C_0 : \chi_0 (X)\equiv 0 \, \text{mod} \,m \} .\] 
It is an easy exercise to verify that these are all dense triangulated subcategories of $\C_0$.
\end{example}

\begin{prop} 
A triangulated subcategory $\C$ of $\C_0 (= \F_p)$ is dense in  $\C_0$ if and only if $\C = \C_0^m$ for some $m$.
\end{prop}

\begin{proof} Let $S$ denote the $p$-local sphere spectrum and recall that $\C_0$ is the union of the
following sets. 

$F^1 =\{ \Sigma^k S: \; k \in \mathbb{Z}\}$ 

$F^2  =\{ C: A \rightarrow B \rightarrow C \rightarrow \Sigma A $
\;\mbox{is an exact triangle such that}\; $ A,B \in F^1 \} $ 

$F^3 =\{ C: A \rightarrow B \rightarrow C \rightarrow \Sigma A $
\;\mbox{is an exact triangle such that}\; $ A,B \in F^2 \} $ 
 
\hskip 12mm \vdots 

$F^n =\{ C: A \rightarrow B \rightarrow C \rightarrow \Sigma A $
\;\mbox{is an exact triangle such that}\; $ A,B \in F^{n-1} \} $ 

\hskip 12mm \vdots

Now using the Euler relations, it is clear that the 
Grothendieck group is a cyclic group generated by $[S]$. 
To see that this group is infinite, consider the map                                    
\[\chi_0 :\{\mbox{isomorphism classes of }  \C_0 \}  \rightarrow \mathbb{Z}\] 
which sends an isomorphism class to its
Euler characteristic. By the universal property of the
Grothendieck group we get an induced homomorphism $\phi_0: K_0(\C_0) \rar \mathbb{Z}$.
Since $\chi_0(S)=1$, $\phi_0$ maps the generator $[S]$ of $K_0(\C_0)$ to $1$, and therefore $\phi_0$ is an isomorphism.
The proposition now follows by invoking Thomason's theorem \ref{main}.  
\end{proof}

\begin{prop} $K_0(\C_0) \cong \mathbb{Z}$ as commutative rings.
\end{prop}
\begin{proof} First note that the smash product on $\C_0$ induces a ring structure on $K_0(\C_0)$. 
 Now the isomorphism of abelian groups $\phi_0: K_0(\C_0) \rar \mathbb{Z}$ (from the proof of the previous proposition)
maps the Grothendieck class of the wedge of $n$-copies of the sphere spectrum  to $n$.  Since the smash product distributes over the wedge, we 
conclude that $\phi_0$ is a ring isomorphism.  
\end{proof}

We now give some nice consequences of this proposition.

\begin{cor} Every dense triangulated subcategory of $\F_p$ is a triangulated ideal.
\end{cor}

\begin{proof} We have seen that $K_0(\T) \cong \mathbb{Z}$ as rings. So we apply theorem \ref{imain}
and observe the simple fact that every subgroup of $\mathbb{Z}$ is also an ideal of $\mathbb{Z}$. 
\end{proof}

\begin{cor} If $\B$ is a dense triangulated subcategory of $\F_p$, then $\B$ is prime if and only if $\B =\C_0^0 $ or $\C_0^p$ for some
prime number $p$.  In particular, a dense triangulated subcategory of $\F_p$ is maximal if and only if it is $\C_0^p$ for some prime $p$.
\end{cor}

\begin{proof}  By theorem \ref{imain}, we have a correspondence between the prime (maximal) dense triangulated subcategories
of $\F_p$ and the prime (maximal) ideals of $K_0(\F_p) \,(\cong \ints)$. The corollary now follows by comparing 
the prime ideals and maximal ideals of $\mathbb{Z}$.
\end{proof}

\begin{rem} The last two corollaries are structural results on the subcategories of $\F_p$. It is not clear how one would
establish such results without using this $K$-theory approach of Thomason.  
\end{rem}

\subsection{$\C_1$ - Finite $p$-torsion spectra.}
Having classified the dense triangulated subcategories of $\C_0$, we now
look at the thick subcategory $\C_1$ consisting of the finite $p$-torsion spectra.
A potential candidate which might generate $\C_2$ as a triangulated category is the Moore spectrum $M(p)$.
So it is natural to ask if $M(p^2)$ can be generated from $M(p)$ by iterated cofiberings. The fact that this is possible can be seen
easily using the octahedral axiom which gives the following commutative diagram of cofibre sequences. 
\[
\xymatrix{
S \ar[r]^p  \ar@{=}[d]  & S \ar[r] \ar[d]^p & M(p) \ar@{.>}[d] \\
S \ar[r]^{p^2} \ar[d]^0 & S \ar[r] \ar[d] & M(p^2) \ar@{.>}[d] \\
0 \ar[r]& M(p)  \ar[r] &  M(p)
}
\]
The exact triangle in the far right is the desired cofibre sequence.
Inductively, it is easy to see that $M(p^i)$ can be generated from $M(p)$ using cofibre sequences.
This example motivates the next proposition.

\begin{prop} \label{prop:steve} Every spectrum in $\C_1$ can be generated by iterated cofiberings of $M(p)$.
\end{prop}

\begin{proof} First observe that the integral homology of any spectrum in $\C_1$ consists of 
finite abelian $p$-groups. Further these spectra have  homology concentrated only in a finite range, so we induct on $|H\mathbb{Z}_*(-)|$. 
If $|H\mathbb{Z}_*(X)| = 1 $, that means $X$ is a trivial spectrum and it is obviously generated by $M(p)$.  So assume that 
$|H\mathbb{Z}_*(X)| > 1$ and let $k$ be the smallest integer such that $H\mathbb{Z}_k(X)$ is non-zero. 
Then by the Hurewicz theorem, we know that $\pi_k(X) \cong H\mathbb{Z}_k(X)$. 
So pick an element of order $p$ in $H\ints_k (X)$ and represent it by a map 
$\alpha: S^k \rightarrow X$. Since $p\alpha = 0$, the composite $S^k \stk{p} S^k \stk{\alpha} X$ 
is zero and hence the map $\alpha$ factors through $\Sigma^k M(p)$.  This gives the following diagram 
where the vertical sequence is a cofibre sequence.
\[
\xymatrix{
\Sigma^k S \ar[drr]_0 \ar[r]^p & \Sigma^k S \ar[dr]^{\alpha} \ar[r]  & \Sigma^k M(p) \ar@{.>}[d]^{\alpha'}\\
& & X \ar[d] \\
& & Y 
}
\]
It is easily seen that $H\mathbb{Z}_i(\alpha')$ is non-zero (and hence injective) if $i=k$, and is zero otherwise.
Therefore, from the long exact sequence in integral homology induced by the vertical cofibre 
sequence, it follows that $|H\mathbb{Z}_*(Y)| = |H\mathbb{Z}_*(X)| - p$. The induction
hypothesis tells us  that $Y$ can be generated by cofibre sequences using $M(p)$ and then the above vertical 
cofibre sequence tells us that $X$ can also be generated by cofibre sequences using $M(p)$.
\end{proof}

Note that the regular Euler characteristic is not good for spectra in $\C_1$ because these
spectra are all rationally acyclic and therefore their Euler characteristic is always
zero (no matter which field coefficients we use). The integral homology of these spectra 
consists of finite abelian $p$-groups and that motivates the following definition.

\begin{defn} For any spectrum $X$ in $\C_1$, define 
\begin{equation} \label{E:chi_1}
\chi_1(X) := \sum_{i = -\infty}^{\infty} (-1)^i \log_p |H\mathbb{Z}_i X|,
\end{equation}
and for every non-negative integer $m$, define a full subcategory 
\[\C_1^m := \{ X \in \C_1 : \chi_1(X) \equiv 0\ \text{mod} \, m \}.\]
\end{defn}

\begin{thm} \label{thm:steve}A triangulated subcategory $\C$ of $\C_1$ is dense in $\C_1$ if and only if $\C = \C_1^m$ for some non-negative integer $m$.
\end{thm}

\begin{proof} The function $\log_p |-|$ is clearly an additive function on the abelian category of
finite abelian $p$-groups, i.e., whenever $0 \rightarrow A \rightarrow B \rightarrow C \rightarrow 0$
is a short exact sequence of finite abelian $p$-groups, we have $\log_p |B| = \log_p |A| + \log_p |C|$.
Using this it is elementary to show that if we have a bounded exact sequence of finite abelian $p$-groups,
then the alternating sum of $\log_p|-|$'s is zero. This applies, in particular, to the long exact sequence
in integral homology for a cofibre sequence in $\C_1$. In other words, $\chi_1(-)$ respects the
Euler relations in $\C_1$. So we get an induced a map $\phi_1 : K_0(\C_1) \rar \ints$. From  
proposition \ref{prop:steve}, we know that $K_0(\C_1)$ is a cyclic group generated by $[M(p)]$. 
Since $\phi_1([M(p)]) = 1$, it follows that $\phi_1$ is an isomorphism. The given classification
is now clear from theorem \ref{main}

\end{proof}

The following are  some easy corollaries of the above theorem.

\begin{cor} \label{cor:testforp-torsionspectra} Let $\Delta(X_1,X_2,\cdots,X_k)$ denote the triangulated subcategory generated by 
spectra $X_1, X_2, \cdots,$ and $X_k$. Then we have the following.
\begin{enumerate}
\item If $A$ is a  type-$1$ spectrum, then $\Delta(A)$ consists of all spectra $X$ in $\C_1$ for which $\chi_1(X)$ is divisible by $\chi_1(A)$.
\item If $A$ is a type-$1$ spectrum and $B$ is a spectrum in $\C_1$, then $B$ can be generated from $A$ using cofibre sequences if and only if $\chi_1(B)$  is divisible by 
$\chi_1(A)$.
\item If $A$ is a type-$1$ spectra and $B$ is a spectrum in $\C_1$, then $\Delta(A,B)$ consists of all spectra $X$ in $\C_1$ for which 
$\chi_1(X)$ is divisible by the 
highest common factor of $\chi_1(A)$ and $\chi_1(B)$.
\end{enumerate}
\end{cor}

\begin{rem} Similar results can be proved for  type-$0$ spectra using the Euler characteristic function $\chi_0$.
\end{rem}

We now examine the Grothendieck group of the Verdier Quotient category $\C_0/\C_1$. We begin with some generalities.
If $\A$ is any thick subcategory of a triangulated category $\C$, then we have natural functors fitting
into an exact sequence
\[ \A \rightarrow \C \rightarrow \C/\A \rightarrow 0,\]
where the first functor is the inclusion functor and the second one is the quotient functor into the
Verdier quotient. Applying $K_0(-)$ to the above sequence induces an exact sequence \cite[Page 355, Proposition 3.1]{sga5}, 
\begin{equation} \label{eq:K_0exactseq}
 K_0(\A) \rightarrow K_0(C) \rightarrow K_0(\C/\A) \rightarrow 0.
\end{equation}
The first map in this sequence is in general not injective. Here is an example where such a map
fails to be injective. If $\A = \C_1$ and $\C = \C_0$, then the inclusion functor 
$\C_1 \hookrightarrow \C_0$  induces a map 
\[K_0(\C_1)\cong \mathbb{Z} \rightarrow K_0(\C_0) \cong \mathbb{Z}.\] 
This map is evidently the zero map because the Euler characteristic $\chi_0$ applied to any rational acyclic gives zero. 
Exactness of the sequence $K_0(\C_1) \rar K_0(\C_0) \rar K_0(\C_0/\C_1) \rar 0$ gives the following corollary.

\begin{cor} $K_0(\C_0/\C_1) \cong \mathbb{Z}$ and is generated by the Grothendieck class of the image of
the $p$-local sphere spectrum under the quotient functor $\C_0 \rightarrow \C_0/\C_1$.
\end{cor}

Since the thick subcategories of $\F_p$ are all nested ($\C_{n+1} \subseteq \C_n$), we get an exact 
sequence of triangulated functors,
\[ \C_{n+1} \rightarrow \C_n \rightarrow \C_n/\C_{n+1}.\]
Applying the functor $K_0(-)$ gives an exact sequence of groups,
\[ K_0(\C_{n+1})\rightarrow K_0(\C_n) \rightarrow K_0(\C_n/\C_{n+1}) \rightarrow 0.\]
We do not know much about these groups beyond the fact that they are countably generated abelian groups.
The hard and interesting thing here is to determine these groups and understand the map  
$K_0(\C_{n+1})\rightarrow K_0(\C_n)$ in the above exact sequence.

\subsection{$\C_n$ - Higher thick subcategories.} 
Now we want to study the triangulated subcategories of the thick subcategories $\C_n$, for $n > 1$. 
To this end, we make use of the spectra related to $\BP$ and their homology theories to construct some Euler 
characteristic functions.  Recall that there is a ring spectrum called the Brown-Peterson spectrum (denoted by $\BP$) whose coefficient ring is given by
$\mathbb{Z}_{(p)}[v_1,v_2,\cdots,v_n,\cdots]$, with $|v_i| = 2p^i-2$. Associated to $\BP$ we have, for each $n \ge 1$, the Johnson-Wilson spectrum
$\bp{n}$ whose coefficient ring is  $\mathbb{Z}_{(p)}[v_1,\cdots,v_n]$.
We begin with some lemmas. 

\begin{lemma} If $X$ is any finite $p$-torsion spectrum, then $\bp{n}_i(X)$ is
a finite abelian $p$-group for all $i$.
\end{lemma}
\begin{proof} Consider the Atiyah-Hirzebruch spectral sequence,
\[ E^2_{s,t} = H_s(X, \bp{n}_t) \implies \bp{n}_{s+t}(X).\]
Since $\bp{n} _* = \mathbb{Z}_{(p)}[v_1,\cdots,v_n] $, we observe that 
$\bp{n}_t$ is a free $\mathbb{Z}_{(p)}$-module of finite rank  and hence is torsion free. 
Since $X$ is a finite $p$-torsion spectrum, it follows that the $p$-local  homology  $H_s(X, \bp{n}_t)$ consists of finite abelian $p$-groups.
So we conclude that the $E^2$ term, and hence the $E^n$ term for all $n \ge 2$, consists of finite abelian $p$-groups.  
Since the spectral sequence collapses after a finite stage, clearly the same is true for $E^{\infty}$. The groups (abelian $p$-groups) on the lines with
slope $-1$ are the associated graded groups of $\bp{n}_i(X)$. From this it is easy to see that $\bp{n}_i(X)$ is
a finite abelian $p$-group for all $i$.
\end{proof}

For each $n\ge1$, recall that there is a generalised Moore spectrum of type-$n$. For $n=1$ this is 
just the Moore spectrum $M(p)$. For $n=2$ this is the cofibre of a self map of the Moore spectrum 
\[ \Sigma^{|v_1^{i_1}|} M(p) \stk{v_1^{i_1}} M(p) \]  
that induces multiplication by $v_1^i$ (for some $i$) in $\BP_*$ homology  and is denoted by
$M(p,v_1^{i_1})$. For higher values we define this inductively: A type-$n$ generalised Moore spectrum is obtained by taking the
cofibre of a self map
\[\Sigma^{|v_n^{i_n}|} M(p,v_1^{i_1},\cdots,v_{n-1}^{i_{n-1}}) \stk{v_n^{i_n}}
  M(p,v_1^{i_1},\cdots,v_{n-1}^{i_{n-1}}) \]
that induces multiplication by $v_n^{i_n}$ (for some $i_n$) in $\BP_*$ homology and is inductively 
denoted by $ M(p,v_1^{i_1},\cdots,v_{n-1}^{i_{n-1}}, v_n^{i_n})$. For sufficiently large powers of the $v_i$'s these spectra
are known to exist \cite{hs}. However the existential problem of such spectra with specified exponents seems to be a very hard question.

\begin{lemma} Fix $n \ge 1$ and let $X$ be any spectrum in $\C_n$. Then $\bp{n-1}_i X$ is zero for 
all but finitely many $i$. 
\end{lemma}

\begin{proof} The strategy here is a thick subcategory argument. Say that a finite $p$-local 
spectrum $X$ has the property P if $\bp{n-1}_i X$ is zero for 
all but finitely many $i$. It is straightforward to verify that if 
$A \rightarrow B \rightarrow C \rightarrow \Sigma A$
is a cofibre sequence of finite $p$-local spectra and if any two spectra in this cofibre sequence have the property
P, then so does the third. It is also clear that if $A \amalg B$ has the 
property P, then so do $A$ and $B$. Thus the full subcategory of spectra in $\F_p$ that satisfy the property P is a thick 
subcategory. Now by the Hopkins-Smith thick subcategory theorem, we will be done if we can exhibit one
type $n$-spectrum for which P holds. Consider a generalised type-$n$ Moore spectrum
$M(p, v_1^{i_1}, \cdots,v_{n-1}^{i_{n-1}})$ and observe that 
\[ \bp{n-1}_* M(p, v_1^{i_1}, \cdots,v_{n-1}^{i_{n-1}} ) = 
\frac{\mathbb{Z}_{(p)}[v_1,\cdots,v_{n-1}]}{(p,v_1^{i_1},\cdots,v_{n-1}^{i_{n-1}})}.  \]
This homology is a truncated polynomial algebra and therefore is concentrated only in a finite range. So we are done.
In fact, it is not hard to see that the the full subcategory of all finite $p$-local spectra that satisfy the property P is precisely
the thick subcategory $\C_n$.
\end{proof}

With these lemmas at hand, we can define a function $\chi_n : \C_n \rightarrow \mathbb{Z}$  as
\begin{equation} \label{E:chi_n}
 \chi_n X = \sum_{i=-\infty}^{\infty} (-1)^i \log_p |\bp{n-1}_i (X) |.
\end{equation}
The previous two lemmas confirm that this is a well-defined function. It is straight forward to 
verify that this function is an Euler characteristic function. Moreover the $\bp{n-1}_*$ homology of the generalised type-$n$ Moore
spectrum $M(p, v_1^{i_1}, \cdots,v_{n-1}^{i_{n-1}})$
is non-trivial and is concentrated in a finite range of even degrees and hence
$\chi_n M(p,v_1^{i_1}, \cdots, v_{n-1}^{i_{n-1}})$ is non-zero.  So by the universal property
of the Euler characteristic function $\chi_n$, we have the following split short exact sequence.
\[
\xymatrix{
0 \ar[r] & \Ker(\phi_n) \ar[r] & K_0(\C_n) \ar[r] \ar[dr]_{\phi_n} & \im (\phi_n) \ar[r] \ar@/_16pt/[l] \ar@{^{(}->}[d] & 0 \\  
         &                     &                  & \ints                             & 
}
\]
This discussion can be summarised in the following proposition.

\begin{prop} For each $n\ge 1$, $K_0(\C_n)$ has a direct summand isomorphic to $\mathbb{Z}$.
 This gives, for each 
$k \ge 0$, a dense triangulated subcategory of $\C_n$ defined by 
\[ \C_n^k := \{ X \in \C_n : \chi_n (X) \equiv 0\!\!\! \mod{l_n k} \} \]
where $l_n$ is a generator for the cyclic group $\im (\phi_n)$.
\end{prop}

Note that these triangulated subcategories correspond to the subgroups $n\mathbb{Z} \oplus \Ker(\phi_n)$ and hence are dense in $\C_n$
by theorem \ref{main}.

\begin{rem} The above proposition  recovers the dense triangulated subcategories $\C_1^k$  of theorem \ref{thm:steve}. In fact, the 
Euler characteristic function (\ref{E:chi_n}) agrees with (\ref{E:chi_1}) when $n=1$. This is because
$\BP \la 0 \ra \cong H\mathbb{Z}_{(p)}$ (the Eilenberg-Mac Lane spectrum for the integers localised at $p$),
and for $p$-torsion spectra, ${H\mathbb{Z}_{(p)}}_* (X) \cong  H\mathbb{Z}_* (X)$. 
Thus, for $X$ in $\C_1$, we have $\bp{0}_* X \cong H\mathbb{Z}_* (X)$ and consequently the two Euler characteristic functions agree.
\end{rem}
 
It is easily seen that if the Smith-Toda complex $V(n) ( := M(p,v_1,v_2,\cdots v_{n-1}))$ exists, then the map $\phi_n: K_0(\C_n) \rar \ints$ is surjective. 
The natural thing to do now is to determine the kernel of $\phi_n$. This turns out to be a very hard problem. It is known that  
$V(1)$ exists at odd primes. In view of this Frank Adams made the following conjecture.

\begin{conj} \label{conj:Adams}
\cite[Page 529]{MayTho} The generalised Moore spectrum $M(p,v_1)$ ($p$ odd) generates the thick subcategory $\C_2$ as a triangulated category.
\end{conj}

This conjecture is equivalent to saying that $\Ker(\phi_2)=0$, or equivalently that $K_0(\C_2) \cong \mathbb{Z}$. 
Adams \cite[Page 528]{MayTho} also asked the following weaker question.

\vskip 3mm 
\noindent
\textbf{Question:} What is a good set of generators for $\C_n$ ?  (Recall that a set $A$ generates $\C_n$ if the smallest 
triangulated category that contains $A$ is $\C_n$.) 
\vskip 3mm

This is a very important question and one answer was given by Kai Xu \cite{xu}. Before we can state his result we
need to set up some terminology.

Recall that a spectrum $X$ is \emph{atomic} if it does not admit any non-trivial idempotents, i.e., if $f \in [X,X]$
is such that $f^2 = f$, then $f = 0$ or $1$. Since idempotents split in the stable homotopy category, this is
also equivalent to saying that $X$ does not have any non-trivial summands. Now if $X$ is a finite $p$-torsion atomic
spectrum, then the finite non-commutative ring $[X, X]$ of degree zero self maps of $X$ 
modulo its Jacobson radical (intersection of all the left maximal ideals)
is isomorphic to a finite field \cite{adamskuhn}. So for every $p$-torsion atomic spectrum $X$, we can define a 
non-negative integer $e(X)$ by the isomorphism 
\[ [X,X]/rad \cong \mathbb{F}_{p^k}.\]

\begin{example} It can be easily verified that for all $i \ge 1$, $M(p^i)$ is an atomic spectrum with $e$-value one.
\end{example}

The main result of \cite{xu} which uses the nilpotence results \cite{dhs,hs} as the main tools then states:

\begin{thm}\cite{xu} For any every pair of natural numbers $(n, k)$, there is an atomic spectrum $X$
of type-$n$ such that $e(X)=k$. Further, if $C(n,k)$ denotes the triangulated subcategory of $\C_n$
generated by the type-$n$ atomic spectra with $e(-) \le k$, then
$C(n,k) = \C_n$.
\end{thm}
 
So this gives an answer to Adams's question. One set of generators for $\C_n$ can be taken to be the collection
of all type $n$ atomic spectra $X$ with $e(X)=1$.  The next natural question is how big is this set? Is this a finite set? If so, 
then we can infer that the Grothendieck groups $K_0(\C_n)$ are  all finitely generated abelian groups. But we do not know the answer to this question.

Using this theorem of Xu, we can now revise Adams's conjecture as follows.

\begin{conj} At odd primes, $V(1)$ generates (by iterated cofiberings) all  type-$2$ atomic spectra $X$ such that $e(X)=1$.
\end{conj}


\section{The lattice of triangulated subcategories} \label{se:lattice}

We will now study the lattice of triangulated subcategories of 
finite $p$-local spectra. First recall that the thick subcategories of finite $p$-local
spectra are nested \cite{mit}, i.e., 
\[ \cdots \subsetneq \C_{n+1} \subsetneq \C_n \subsetneq \C_{n-1} \subsetneq \cdots \cdots \subsetneq \C_1 \subsetneq \C_0. \]
Our goal now is to understand how the triangulated subcategories $\C_n^k$ fit in this chromatic chain. We begin with a simple proposition.

\begin{prop} $\C_1$ is contained in every dense triangulated subcategory of $\C_0$, i.e., $\C_1 \subsetneq \C_0^k$ for all $k \ge 0$.
\end{prop}

\begin{proof} If $X$ is in $\C_1$, then its rational homology is trivial and hence
its rational Euler characteristic $\chi_0(X)$ is zero. Therefore
$X$ belongs to every dense triangulated subcategory of $\C_0$. It is easy to see that the spectrum $S \vee \Sigma S$ belongs  to $\C_0^k - \C_1$ 
for all $k \ge 0$, therefore the containment $\C_1 \subseteq \C_0^k$ is proper.
\end{proof}

Motivated by this proposition, we wondered if it is true that
for all non-negative integers $n$ and $k$,
\[ \C_{n+1} \subsetneq \C_n^k  \subseteq \C_n. \]
We now proceed to show that this is indeed the case. 

\subsection{An Atiyah-Hirzebruch spectral sequence.}

In this subsection, we prove that $\C_2 \subsetneq \C_1^k$ for all $k \ge 0$. 
Our main tool will be an Atiyah-Hirzebruch spectral sequence. We begin with two lemmas -- the first one is an elementary algebraic fact and the second one is a standard topological fact.

\begin{lemma}\label{parallel} If $A:= \cdots \rar  0 \rar A_1 \rar A_2 \rar \cdots \rar A_k \rar 0 \rar \cdots$ is a bounded chain complex of finite $p$-groups,
then $\sum_i (-1)^i \log_p |A_i| $= $\sum_i (-1)^i \log_p |H_i(A)|$.
\end{lemma}

\begin{proof} This is left as an easy exercise to the reader.
\end{proof}

\begin{lemma} \label{le:adams}  The thick subcategory $\C_2$ consists  of all finite $p$-torsion spectra whose complex 
$K$-theory is trivial.
\end{lemma}

\begin{proof} We use a result of Adams \cite{adams} which states that the complex $K$-theory localised at $p$ 
splits as a wedge of suspensions of $E(1)$. 
More precisely, $K_p = \bigvee _{i=0}^{i=p-2} \Sigma^{2i} E(1)$. In particular, $\la K_p \ra = \la E(1)\ra$. With
this at hand, we get the following equalities of Bousfield classes.
\begin{eqnarray*}
\la K \ra &=& \bigvee_p \la K_p \ra \\ 
          &=& \bigvee_p \la E(1)\ra \\
          &=& \bigvee_ p ( \la K(0)\ra \bigvee \la K(1) \ra ) . 
\end{eqnarray*}
The last equality follows from \cite[Theorem 2.1(d)]{rav}.
Now it is clear from these equations that for $X$ finite and $p$ torsion, $K_* X = 0$ if and only if $K(1)_* X = 0$.
\end{proof}

\begin{prop} \label{prop:newevidence}  $\C_2$ is properly contained in every dense triangulated subcategory of $\C_1$, i.e., 
$\C_2 \subsetneq \C_1^k$ for all $k \ge 0$. 
\end{prop}

\begin{proof} Recall that $\C_1^k$ is the collection
of $p$-torsion spectra $X$ for which $\chi_1(X)$ is divisible by $k$. So it is clear that $X \in \C_1^k$ for all $k$ 
if and only if $\chi_1(X)=0$. 
Therefore by lemma \ref{le:adams} we have to show that if $X$ is a finite $p$-torsion spectrum for which
$K_*(X) = 0$, then $ \chi_1(X) := \sum (-1)^i \log_p|H\mathbb{Z}_i X | = 0$.

\begin{figure}
\begin{center} 
\scalebox{.75}{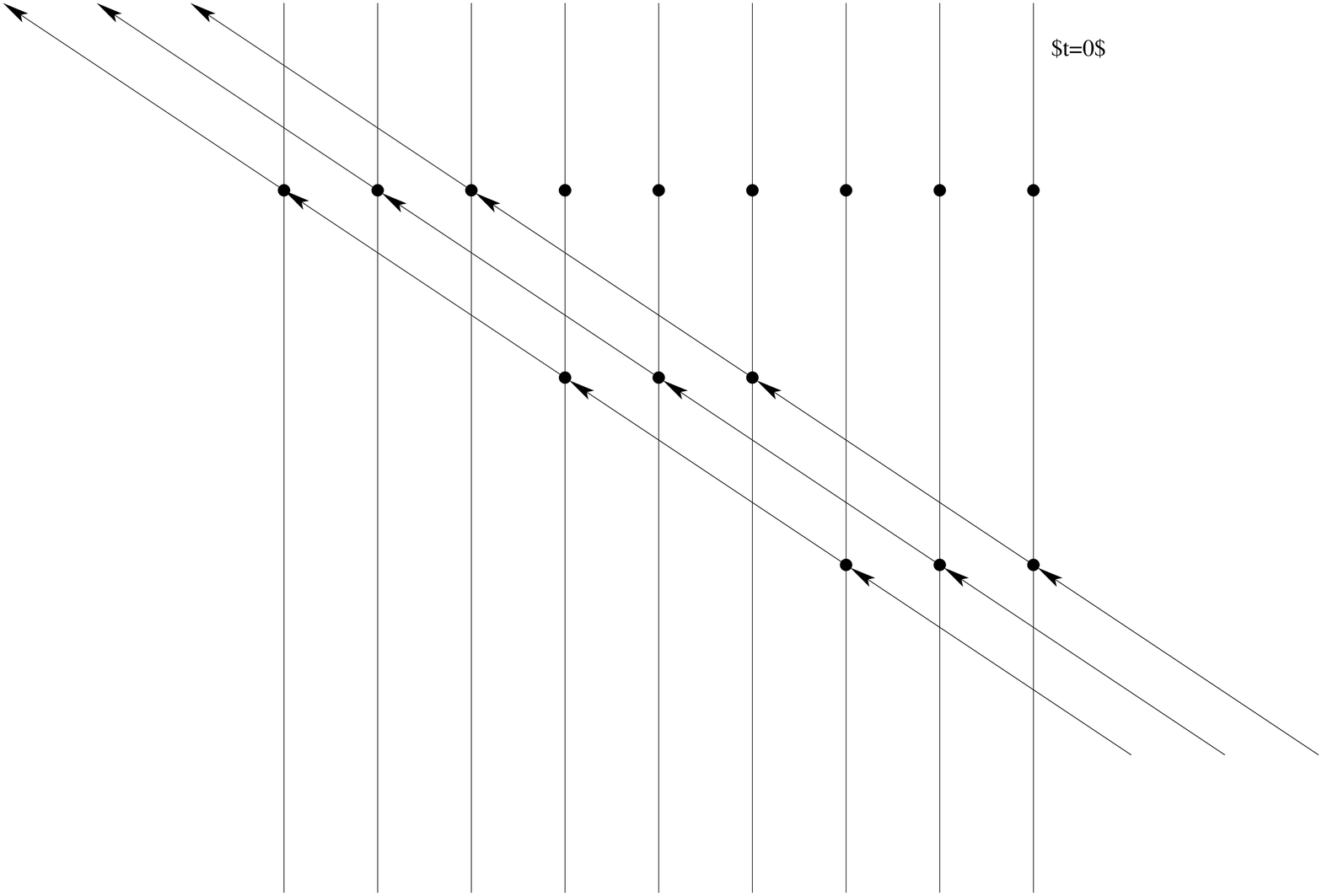}
\caption{Atiyah-Hirzebruch spectral sequence: $E^3$ - term}
\label{E3}
\end{center}
\end{figure}

Recall that the spectrum $K$ (also denoted by $BU$) is a ring spectrum whose 
coefficient group is given by the complex Bott periodicity theorem; $K_* = \mathbb{Z}[u,u^{-1}]$ where $|u|=2$.
We make use of the Atiyah-Hirzebruch spectral sequence 
\[ E^2_{s,t} = H_s(X; K_t) \implies K_{s+t} X. \]
converging strongly to the $K$-theory of $X$. The differentials ($d_r$) in this spectral sequence have bidegrees $|d_r| = (-r, r-1)$.
Note that $X$ being a finite spectrum, the $E^2$ page is concentrated in a vertical strip of finite width (see Figure \ref{E3}) and 
therefore the spectral sequence collapses after a finite stage. Since it converges to $K_*(X)$, which is zero by hypothesis,
we conclude that for all sufficiently large $n$, $E^n = E^{\infty} = 0$.

Next we claim that the function
\[ n \longmapsto \underset{i}{\sum} (-1)^i \log_p |E^n_{i,0}|\]
is a constant function.

Assuming this claim, we will finish the proof of the proposition. When $n=2$, this function 
takes the value  $\sum_i (-1)^i \log_p|H\mathbb{Z}_i\, X | = \chi_1(X)$  and for large enough $n$,
the function takes the value $0$ (since $E^n = E^{\infty} = 0$). Since the function is constant (by the above claim), we get 
$\chi_1(X) = 0$.

Now we prove our claim. First note that this spectral sequence is a module over the coefficient ring $K_* \cong \mathbb{Z}[u,u^{-1}]$. Therefore
the differentials commute with this ring action. Also since $u$ is a unit, it acts isomorphically on the spectral sequence; see diagram below.

\[
\xymatrix{
H_{*+3}(X, K_{t-2}) \ar@/^2pc/[rr]^{1}  \ar[r]^u \ar[d]^{d_3} & H_{*+3}(X, K_t) \ar[r]^{u^{-1}} \ar[d]^{d_3} & H_{*+3} (X, K_{t-2}) \ar[d]^{d_3} \\
H_*(X, K_t)  \ar[d]^{d_3} \ar[r]^{u}  & H_*(X, K_{t+2}) \ar[d]^{d_3} \ar[r]^{u^{-1}}  & H_*(X, K_t) \ar[d]^{d_3}   \\
H_{*-3}(X, K_{t+2}) \ar@/_2pc/[rr]_{1} \ar[r]^u & H_{*-3}(X, K_{t+4}) \ar[r]^{u^{-1}} & H_{*-3}(X, K_{t+2})  
}
\]
It is clear from this diagram that $u$ induces periodicity on $E^4$: $E^4_{i,*} \cong E^4_{i-2,*}$. 
Also note that just for degree reasons, all the even differentials
are zero. So at $E_3$, where the first potential non-zero differentials occurs, the alternating sum $\sum_i (-1)^i \text{log}_p |E^3_{i,0}|$
can be broken into three parts as shown in the equation below.s
\[ \sum_i (-1)^i \text{log}_p |E^3_{i,0}| = 
   \sum_{i\equiv 0 (3)} (-1)^i \text{log}_p |E^3_{i,0}| +
   \sum_{i\equiv 1 (3)} (-1)^i \text{log}_p |E^3_{i,0}| +
   \sum_{i\equiv 2 (3)} (-1)^i \text{log}_p |E^3_{i,0}| \]
Now because of the periodicity of the $E^3$ page, we can assemble these terms in the above equation along three parallel lines as shown in the $E_3$ page, where the
term with congruence class $l$ modulo three corresponds to the line with $t$-intercept $l$; see figure \ref{E3}.
Now we can apply lemma \ref{parallel} along each of these lines which are bounded 
complexes of finite abelian $p$ groups to pass to the homology groups. Invoking the periodicity of the 
differentials again, we conclude that the new alternating sum thus obtained is
equal to  $\sum_i (-1)^i \text{log}_p |E^4_{i,0}|$. Now an easy induction will complete the proof of the claim:
At ($E^r, d_r$),  for $r$ odd, we decompose the alternating sum $\sum_i (-1)^i \log_p|E_{i,0}^r|$ into $r$ parts, one for each congruence class modulo
$r$, and assemble these terms along $r$ parallel lines on the $E^r$ page such that the term whose congruence class is 
$l$ modulo $r$ corresponds to the line with $t$-intercept $l$. Periodicity of the differentials can be used (as before) 
to complete the induction step.

So we have shown that $\C_2 \subseteq \C_1^k$ for all $k \ge 0$. To see that this inclusion is strict, observe that the 
$p$ torsion spectrum $M(p) \vee \Sigma M(p)$ belongs to $\C_1^k - \C_2$ for all $k$.
\end{proof}

\begin{cor} $K_0(\C_1/C_2) \cong \ints$ generated by the Grothendieck class of the image of the Moore spectrum under the quotient functor
$\C_1 \rar \C_1/\C_2.$
\end{cor}

\begin{proof} On applying $K_0(-)$ to the sequence $\C_2 \rar \C_1 \rar \C_1/\C_2 \rar 0$, we get an exact sequence of abelian groups:
$K_0(\C_2) \rar K_0(\C_1) \rar K_0(\C_1/\C_2) \rar 0$; see (\ref{eq:K_0exactseq}).  The first map in this sequence is the zero map by proposition 
\ref{prop:newevidence}, and $K_0(\C_1)$ was shown to be infinite cyclic on $[M(p)]$. So the corollary follows by combining these two facts.
\end{proof}

The inclusion $\C_2 \subsetneq \C_1^k \; \forall \, k$, which we have just established, gives some new evidence to conjecture \ref{conj:Adams} of Adams. 
To see this, first note that since $M(p,v_1)$ is a cofibre of an even degree self map 
of $M(p)$,  $\chi_1(M(p,v_1)) = 0$. Now if  $M(p,v_1)$ generates $\C_2$  by cofibre sequences (Adam's conjecture), 
then  $\chi_1(X)=0$ for all $X \in \C_2$. This now clearly implies that $\C_2 \subsetneq \C_1^k \; \forall \, k$.

One can try to test this conjecture by asking whether we can generate $M(p,v_1^i)$ (they all exist at the odd primes) from $M(p,v_1)$.
The octahedral axiom gives the following commutative diagram of exact triangles. (We are ignoring all suspensions in the pictures
below for clarity.) The exact triangle on the far right, together with a straight forward induction, tells us that the
generation in question is possible.
\[
\xymatrix{
M(p) \ar[r]^{v_1}  \ar@{=}[d]  & M(p) \ar[r] \ar[d]^{v_1} & M(p,v_1) \ar@{..>}[d] \\
M(p) \ar[r]^{v_1^2} \ar[d]^0 & M(p) \ar[r] \ar[d] & M(p, v_1^2) \ar@{..>}[d] \\
0 \ar[r]& M(p,v_1)  \ar[r] &  M(p,v_1)
}
\]
We can go one step further and ask if it is possible to  generate $M(p^i,v_1^j)$  from $M(p,v_1)$ using cofibre sequences?
Results from nilpotence and periodicity \cite{hs} give the following  partial answer to this question.

\begin{prop} For every fixed positive integer $k > 0$, there exists infinitely many positive integers $j$ for which
$M(p^k,v_1^j)$ can be generated from $M(p,v_1)$ using cofibre sequences.
\end{prop}

\begin{proof} It suffices to show that for each fixed $k$, there exists an integer $j_0$ such that  $M(p^k, v_1^{j_0})$ can be generated from 
$M(p,v_1)$. For, it then follows from the octahedral axiom that we can obtain $M(p^k,v_1^{j})$ from $M(p,v_1)$ whenever $j$ is a multiple of $j_0$.

First note that there is a cofibre sequence 
\[ M(p) \stk{p} M(p^2) \stk{} M(p).\]
Since both $M(p)$ and $M(p^2)$ are type $1$ spectra, we know that they admit $v_1$ self maps \cite[Theorem 9]{hs}. By taking iterated
compositions of these self maps ($v_1$ and $z$ for $M(p)$ and $M(p^2)$ respectively), we can assume that there is commutative diagram 
\cite[Corollary 3.8]{hs} as shown below.
\[ 
\xymatrix{
M(p) \ar[r]^{v_1^l} \ar[d]  & M(p) \ar[d]\\
M(p^2) \ar[r]^{z^m} & M(p^2)
}
\]
This diagram can be extended to a diagram of exact triangles using the octahedral axiom.
\[
\xymatrix{
M(p) \ar[r]^{v_1^l}  \ar[d]  & M(p) \ar[r] \ar[d] & M(p,v_1^l) \ar@{..>}[d] \\
M(p^2) \ar[r]^{z^m} \ar[d]  & M(p^2) \ar[r] \ar[d] & M(p^2, v_1^l) \ar@{..>}[d] \\
M(p) \ar[r]^{t} & M(p)  \ar[r] &  \Cone(t)
}
\]

By \cite[Corollary 3.9]{hs}, we know that the fill in map $t$ in the above diagram is a $v_1$ map
and therefore its cone is a generalised type-$2$ Moore spectrum. Again by the uniqueness of
$v_n$ maps \cite[Corollary  3.7]{hs}, we know that there exist integers $a$ and $b$ for which $t^a=v_1^b$;
see the diagram below.
\[ 
\xymatrix{
M(p) \ar[r]^{t^a} \ar@{=}[d]  & M(p) \ar@{=}[d]\\
M(p) \ar[r]^{v_1^b} & M(p)
}
\]
Now we can replace the bottom left hand map of the above $3 \times 3$ diagram with $v_1^b$ and iterate all
the maps further to get commutativity of both the left hand squares. This gives us the desired cofibre
sequence (the far right triangle in the diagram below).
\[
\xymatrix{
M(p) \ar[r]^{v_1^l}  \ar[d]  & M(p) \ar[r] \ar[d] & M(p,v_1^l) \ar@{..>}[d] \\
M(p^2) \ar[r]^{z^m} \ar[d]  & M(p^2) \ar[r] \ar[d] & M(p^2, v_1^l) \ar@{..>}[d] \\
M(p) \ar[r]^{v_1^l} & M(p)  \ar[r] &  M(p,v_1^l)
}
\]
Finally note that we can generate $M(p,v_1^l)$ from $M(p,v_1)$. So this completes the proof of the 
proposition for $k=2$. The general case is a straight forward generalisation.
\end{proof}

All these results give only some evidence for Adams's conjecture. The conjecture, however, still remains open.

\subsection{A Bockstein spectral sequence}

Our goal now is to prove: $\C_{n+1} \subsetneq \C_n^k$ for all $n$ and $k$. We mimic our strategy for the case
$n=1$ by replacing the Atiyah-Hirzebruch spectral sequence with a Bockstein spectral sequence.  The above inclusion 
is an easy corollary of the following theorem.

\begin{thm} If $X$ is spectrum in $\C_{n+1}$, then  $\chi_n (X) := \sum (-1)^i \log_p |\bp{n-1}_i X| = 0$.
\end{thm}

\begin{cor} $\C_{n+1} \subsetneq \C_n^k$ for all $k \ge 0$ and all $n \ge 1$.
\end{cor}

\begin{proof} By the above theorem, for a spectrum $X$ in $\C_{n+1}$, $\chi_n(X) = 0$. Therefore $X$ belongs to $\C_n^k$ for all $k$.
The spectrum $F \wedge \Sigma F$, where $F$ is a  of type-$n$, belongs to $\C_n^k - \C_{n+1}$. So $\C_{n+1} \subsetneq \C_n^k$ for all $k \ge 0$.
\end{proof}

\n
We now outline the strategy for proving the above theorem. This is very similar to the proof of proposition \ref{prop:newevidence}.
First recall that $\C_{n+1}$ can also be characterised as 
the collection of finite $p$-local spectra that are acyclic with respect to $E(n)$. So we seek a strongly convergent spectral sequence 
                                 \[ E_r^{**} \implies E(n)_* X,\] 
whose $E_1$ term is build out of  $\bp{n-1}_* X$. We show that a certain Bockstein spectral sequence has this property
and that it collapses after a finite stage. Finally we work backward to conclude that
$\chi_n X = 0$. 

\begin{figure}
\begin{center}
\scalebox{.75}{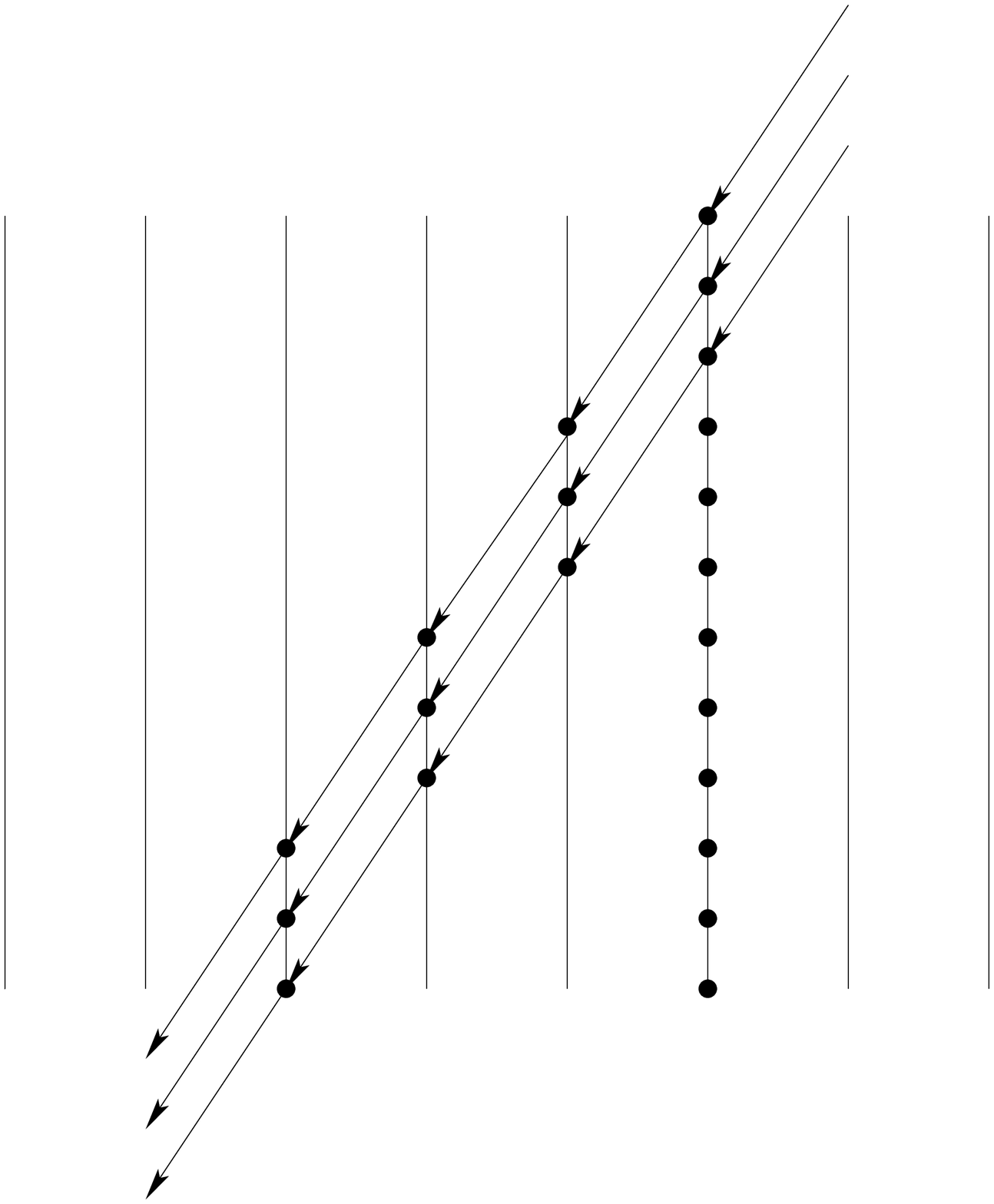}
\caption{Bockstein spectral sequence: $E_1$ - term}
\label{E1}
\end{center}
\end{figure}

Now we proceed to construct such a spectral sequence. Fix an integer $n \ge 1$ and recall that 
\[E(n) = v_n^{-1} \bp{n} = \hocolim{v_n} \bp{n}.\]
This gives a sequence of cofibre sequences that fit into a diagram extending to infinity in both directions as shown below.
\[
\xymatrix{
\cdots \ar[r]^{v_n\ \ \ \ \ \ } & \sus{|v_n|} \bp{n} \ar[r]^{\ \ v_n} \ar[d] & \bp{n} \ar[r]^{v_n \ \ \ \  } \ar[d] & \sus{-|v_n|} \bp{n} \ar[r]^{\ \ \ \ v_n} \ar[d] 
& \cdots \cdots \ar[r] & E(n) \\
\cdots              &  \sus{|v_n|} \bp{n-1} \ar[ul]|{\circ}           &        \bp{n-1} \ar[ul]|{\circ}  &  \sus{-|v_n|} \bp{n-1} \ar[ul]|{\circ}  
& \cdots  &                 
}
\]
Since both the functors  $(-)\wedge X$ and $\pi_*(-)$ commute with the functor $\hocolim{v_n}(-)$, smashing the above
diagram with $X$ and taking $\pi_*$ gives an exact couple of graded abelian groups:
\[
\xymatrix{
\cdots \ar[r]^{v_n \ \ \ \ \ \ \ \ } & \sus{|v_n|} \bp{n}_*X \ar[r]^{\ \ \ \ v_n} \ar[d] & \bp{n}_*X \ar[r]^{v_n\ \ \ \ \ \ \  } \ar[d] & \sus{-|v_n|} \bp{n}_*X  \ar[d] 
 \cdots\cdots  \ar[r] & E(n)_*X \\
\cdots              &  \sus{|v_n|} \bp{n-1}_*X \ar[ul]|{\circ}          &        \bp{n-1}_*X \ar[ul]|{\circ}  &  \sus{-|v_n|} 
\bp{n-1}_*X \ar[ul]|{\circ}  & \cdots                   
}
\]
This exact couple gives rise to a (Bockstein) spectral sequence $E_r^{*,*}$ in the usual way. We choose a convenient grading so that 
the $E_1$ term is concentrated in a horizontal strip of finite width, i.e., $E_1^{*,q} = 0$ for $|q| >> 0$ (see figure \ref{E1}). This can be ensured by setting
\[D_1^{p,q} = \Sigma^{-p|v_n|} \bp{n}_{-p|v_n|+q} X, \]
\[E_1^{p,q} = \Sigma^{-p|v_n|} \bp{n-1}_{-p|v_n|+q} X.\]
With this grading one can easily verify that the differentials have bidegrees given by $|d_r|=(-r, -r|v_n|-1)$. 
It is clear from figure \ref{E1} that after a finite stage all the differentials
exit the horizontal strip and therefore the spectral sequence collapses. Another important fact that we need about this
Bockstein spectral sequence is the periodicity of all the differentials. More precisely, for all integers $p$ and $r \ge 0$,
$E_r^{p,*} = E_r^{p+1,*}$, and further the following diagram commutes.
\[
\xymatrix{
E_r^{p,*}\;\;\;  \ar[d]^{=}  \ar[r]^{d_r\;\;\;\;\;\;\;} & \;\;E_r^{p-r,*-1-r|v_n|} \ar[d]^{=}\\
E_r^{p+1,*}  \ar[r]^{d_r\;\;\;\;\;\;\;\;} & \;  E_r^{p+1-r,*-1-r|v_n|}
}
\]
Now we show that this spectral sequence converges strongly to $E(n)_*(X)$. For this, we make use of a theorem
of Boardman \cite{boardman}. Before we can state his result we have to recall some of his terminology. Consider an exact couple 
of graded abelian groups:

\[
\xymatrix{
\cdots \ar[r]^{i} & A^{s+1} \ar[r]^{i} \ar[d]^j & A^s \ar[r]^{i}\ar[d]^j & A^{s-1} \ar[r]^{i} \ar[d]^j 
& A^{s-2} \ar[r]^{i} \ar[d]^j & \cdots  \ar[r] & A^{-\infty} \\
\cdots              &  E^{s+1} \ar[ul]^k|{\circ}           &   E^s \ar[ul]^k|{\circ}  &  E^{s-1} \ar[ul]^k|{\circ}  
& E^{s-2} \ar[ul]^k|{\circ}  &\cdots \cdots  &                 
}
\]

This gives the following filtrations of the groups $E^s$ by cycles and boundaries of the differentials
in the spectral sequence that arises from this exact couple,
\[ 0 = B_1^s \subset B_2^s \subset B_3^s \subset \cdots \
   \cdots \subset Z_3^s \subset Z_2^s \subset Z_1^s=E^s, \]
where $Z_r^s := k^{-1}(\text{Im}\, [i^{(r-1)}: A^{s+r} \rar A^{s+1}])$, and $B_r^s:= j \ker[i^{(r-1)}: A^s \rar A^{s-r+1}]$.

Associated to this exact couple and the resulting spectral sequence, \cite{boardman} defines the following groups.
These definitions also hold when $n=\infty$.
\begin{itemize}
\item $A^{-\infty}:= \colim{s} A^s$, \hspace{3 mm}  $A^{\infty}:= \clim{s} A^s$, \hspace{3 mm} $RA^{\infty}:= \rlim{s} A^s$
\item $K_n A^s:= \ker[i^{(n)}: A^s \rar A^{s-n}] $, \hspace{3mm} $\text{Im}^r A^s:= \text{Im}[i^{(r)}: A^{s+r} \rar A^s ]$
\item $K_n \text{Im}^r A^s := K_n A^s \cap \text{Im}^r A^s$
\item $W=\colim{s} \rlim{r} K_{\infty} \text{Im}^r A^s $
\item $RE_{\infty}^s = \rlim{r} Z_r^s$.
\end{itemize}

The main theorem of \cite{boardman} then states:

\begin{thm}\cite[Theorem 8.10]{boardman}
The spectral sequence arising from the above exact couple converges strongly to $A^{-\infty}$ if the          
obstruction groups $A^{\infty}$, $RA^{\infty}$, $W$, and $RE_{\infty}$ are zero.
\end{thm}

We recall a few elementary facts about inverse limits before we can apply this theorem to our Bockstein 
spectral sequence. The proofs follow quite easily from the universal properties of these
limit functors. Parts of this lemma can also be derived from the more general Mittag-Leffler condition.

\begin{lemma} \label{facts-limits} Consider a sequence of groups (graded) and homomorphisms:
\[ \cdots \rar A^{s+1} \rar A^s \rar A^{s-1} \rar \cdots.\]
Then the following statements hold.
\begin{itemize}
\item If, for some integer $M$, the composite of $M$ consecutive maps in this sequence is always null, 
then  $\clim{s} A^s = \rlim{s} A^s = 0$.
\item If there is an integer $L$ such that for all $s \ge L$ the map $A^{s+1} \rar A^s$ is the identity
map, then $\clim{s} A^s = A^{L}$ and $\rlim{s} A^s = 0$
\end{itemize}
\end{lemma}

We now show that our Bockstein spectral sequence converges strongly to $E(n)_* X$ by verifying the hypothesis of
Boardman's theorem. So we apply his theorem to the sequence,
\begin{equation}\label{vnseq}
\cdots \stk{v_n}  \sus{|v_n|} \bp{n} \stk{v_n} \bp{n} \stk{v_n} \sus{-|v_n|} \bp{n} \stk{v_n} \cdots  
\end{equation}
For brevity, we shall denote $\sus{-s|v_n|} \bp{n}_* X$ by $A^s$.

\n
(a) $A^{\infty}=0$: By lemma \ref{facts-limits}, all we need to show is that there is some integer $M$
such that the composite of any $M$ consecutive maps in the sequence \eqref{vnseq}
is zero. $X$ being a spectrum in $\C_{n+1}$, $\bp{n}_*X$ is concentrated only
in a finite range. Now since $v_n$ is a graded map of degree $2(p^n-1)$, a sufficiently large iterate of $v_n$
vanishes, whence $A^{\infty} = \clim{s} \sus{-p|v_n|} \bp{n}_* X = 0$.

\n
(b) $RA^{\infty}=0$: This is also immediate from part (1) of lemma \ref{facts-limits} because we have already seen that a sufficiently large
iterate of $v_n$ is zero in part (a).

\n
(c) $W=0$:  The colimit of the sequence \eqref{vnseq} is $E(n)_*X$ and this latter group is zero by hypothesis.
It now follows that $K_{\infty} A^s = A^s$. We make use of this fact in the third equality below. 
\begin{eqnarray*}
W &=& \colim{s} \rlim{r} \left( K_{\infty} \text{Im}^r A^s \right) \\
  &=& \colim{s} \rlim{r} \left( K_{\infty} A^s \cap \text{Im}^r A^s \right) \\
  &=& \colim{s} \rlim{r} \left( A^s \cap \text{Im}^r A^s \right) \\
  &=& \colim{s} \rlim{r} \; \text{Im}^r A^s \\
  &=& \colim{s} \rlim{r} ( \cdots 0 \rar 0 \rar \cdots \subseteq \text{Im}^2 A^s \subseteq \text{Im}^1 A^s) \\
  &=& \colim{s} 0  = 0 
\end{eqnarray*}
\n
(d) $RE_{\infty}=0$: This follows from the fact that the spectral sequence collapses after a finite stage. For each fixed $s$, 
all inclusions in the sequence 
\[ \cdots \stk{=} Z_{r+1}^s \stk{=} Z_{r}^s \subseteq \cdots \subseteq Z_3^s \subseteq Z_2^s \subseteq Z_1^s = E_1^s ,  \]
become equalities eventually. So we now invoke part (2) of  lemma \ref{facts-limits} to conclude that 
$RE_{\infty}^s = \rlim{r} Z_r^s = 0$. \\

So we have shown that all the obstruction groups vanish and therefore Boardman's theorem tells us that our
spectral sequence converges strongly to $E(n)_*X$. \\

We now move to the final part of the theorem: $\chi_n(X)=0$. We mimic the argument given for
the fact ``$\C_2 \subseteq \C_1^k$''. The crux of the proof lies in the key observation that the $q$ component
of the differential $d_r$ is always an odd number ($-r|v_n|-1$). We claim that the function
\[k \overset{\phi}{\longmapsto} \sum(-1)^i \log_p |E_k^{0,i}|\]
is constant. Toward this, we decompose the
alternating sum  $\sum(-1)^i \log_p |E_1^{0,i}|$ into $t:=|v_n|+1$ parts as follows:
\[ \sum_{i} (-1)^i \log_p |E_1^{i,0}| = \sum_{i\equiv 0 (t)} + \sum_{i \equiv 1 (t)} + \cdots + \sum_{i \equiv t-1 (t)}\] 
Now the periodicity of all the differentials coupled with the fact that $|v_n|+1$ is an odd number will enable us to  assemble all these terms 
on the right hand side of this equation
along $t$ parallel complexes of differentials on the $E_1$ term. (The term corresponding to the
congruence class $l$ modulo $t$ will correspond to the parallel complex with $q$ intercept $l$; see figure \ref{E1}.) 
We can now apply lemma \ref{parallel} to each of these complexes and pass on to the homology groups
without changing the underlying alternating sum. 
Again using the periodicity of the differentials, we can reassemble all the terms after taking homology to obtain $\phi(2)$. This shows that
$\phi(1)=\phi(2)$. Now a straight forward induction will complete the proof of the claim.

Finally, to see that $\chi_n(X)=0$, observe that when $k=1$, $\phi$ takes the value $\chi_n(X)$, and for all sufficiently
large values of $k$, $\phi$ is zero because our spectral sequence collapses at a finite stage and  converges 
strongly to $E(n)_*X$, which is zero by hypothesis. Since we know that $\phi$ is constant, this completes the
proof of the theorem.
\\

\section{Dense triangulated subcategories of torsion spectra} \label{se:torsionspectra}
 
We consider the category $\FT$ of finite torsion spectra and study 
the dense triangulated subcategories in it. (Note that here we are not working $p$-locally.) 

\begin{prop} \label{prop:K_0FT} $\FT$ can be generated by the set $\{M(p)\}$ using cofibre sequences.
\end{prop}

\begin{proof} This is very similar to the $\C_1$ case. Note that the integral homology of any spectrum 
in $\FT$ is a finite abelian group. So exactly as before we induct on the cardinality $|H\mathbb{Z}_* (-)|$.
Fix some $X \in \FT$ and let $k$ be the smallest dimension in which the integral homology of $X$ is non-zero.
Then by the Hurewicz theorem we know that $\pi_k(X)$ is isomorphic to $H\mathbb{Z}_k (X)$. Since $H\ints_k(X)$
is finite abelian group, for some prime $p$, it has an element of order $p$. Represent this element by a 
map $S^k \rightarrow X$ and proceed exactly as in proposition \ref{prop:steve}  to kill this bottom dimensional homology class.
Now a straightforward induction will complete the proof. 
\end{proof}

\begin{prop} $K_0(\FT) \cong \bigoplus_p \mathbb{Z}$. Further the free generators of the Grothendieck group
are the classes of $M(p)$, for various primes $p$.
\end{prop}

\begin{proof} For each prime $p$, define an Euler characteristic function $\chi_p$ as follows.
\[\chi_p (X) = \sum_{-\infty}^{\infty} (-1)^i \; \log_p|H\mathbb{Z}_i(X) \otimes \mathbb{Z}_{(p)}|.\]
Since localisation is exact, this defines an Euler characteristic function on $\FT$. It is also 
easy to see that for any given $X \in \FT$, there are only finitely many primes $p$ for which
$\chi_p (X)$ is non-zero. Thus we get a map $ \chi (:= \bigoplus_p \chi_p) : \FT \rar \bigoplus_p \mathbb{Z}$ which is clearly 
an Euler characteristic function. For any two primes $p$ and $q$, $\chi_p(M(q)) = \delta_p^q$ (Kronecker function), therefore the map $\chi$ is surjective.
So by the universal property of the Grothendieck groups, we get a surjective group 
homomorphism from $K_0 (\FT)$ to $\bigoplus_p \mathbb{Z}$. The previous proposition tells us that 
$K_0(\FT)$ is generated by the classes $[M(p)]$. The fact that the Grothendieck group admits
a surjective map onto $\bigoplus_p \mathbb{Z}$ implies that there are no non-trivial relations among the classes $[M(p)]$. 
By combining these two facts, we conclude that  $K_0(\FT) \cong \bigoplus_p \mathbb{Z}.$
\end{proof}

It is now clear that every dense triangulated subcategory of $\FT$ is of the form
\[ \FT^k:= \{ X \in \FT \,: \chi(X) \, \equiv 0 \ \text{mod} \, k \} \;\; \mbox{\;for some\;} k.\]

Say that a finite torsion spectrum $X$ is \emph{integrally type-$1$} if $X_{(p)}$ is  of type-$1$ ($K_*X_{(p)}\ne 0$) whenever $X_{(p)} \ne 0$.
The \emph{prime support} $\Supp(X)$ of a spectrum $X$ is the collection of all primes $p$ such that $X_{(p)}\ne 0$. 
Then we have the following corollary of the above classification of dense triangulated subcategories of $\FT$.

\begin{cor} Let $X$ and $Y$ be finite torsion spectra  with $X$ integrally type-$1$. Then $Y$ can generated by iterated cofiberings of $X$ if and only if
$\Supp(Y) \subseteq \Supp(X)$ and $\chi_p(Y)$ is a multiple of $\chi_p(X)$ for all $p$.
\end{cor}

\begin{proof} Since we are working with torsion spectra, we can write $X = \coprod_p X_{(p)}$ and $Y = \coprod_p Y_{(p)}$. Now it is clear that 
$Y$ can be generated by iterated cofiberings of $X$ if and only if $Y_{(p)}$ can be generated by iterated cofiberings of $X_{(p)}$ for all $p$.
By corollary \ref{cor:testforp-torsionspectra}, the latter is possible if and only $\chi_p(X)$ divides $\chi_p(Y)$.
\end{proof}

\section{Questions}

\subsection{Grothendieck groups and Adams's Conjecture}
After all the work we have done in this chapter, we are still far from our goal; our classification of these subcategories is incomplete.
The biggest problem that needs to settled is the classification of all triangulated subcategories of $\F_p$. We believe that such a 
classification  will reveal some hidden ideas behind the chromatic tower which might shed some new light on the stable homotopy category.  
We have seen that $K_0(\C_0)$ is infinite cyclic with the sphere as the generator, and $K_0(\C_1)$ is infinite cyclic with the mod-$p$ 
Moore spectrum as the generator. Computing the Grothendieck groups of higher thick subcategories ($\C_n$, $n \ge 2$), as we have seen, 
is much more complicated. In particular, the  conjecture of Adams that the Smith-Toda complex $V(1)$ 
(at odd primes) is a generator for the thick subcategory $\C_2$ seems intractable. 
We gave only some evidence for this conjecture. The conjecture, however, still remains open.  

Another interesting and related question at this point is the following. In $\C_n$, what is the triangulated
subcategory generated by atomic spectra $X$ of type-$n$ for which $e(X)=k$ ($k$ some fixed positive integer)? By \cite{xu}, this is the 
whole of $\C_n$ if $k=1$; otherwise, this will be some dense triangulated subcategory of $\C_n$.  The immediate question that springs out 
now is whether these triangulated subcategories, when $n=2$, are precisely the subcategories $\C_2^k$? 
A non-affirmative answer to this question will settle Adams's conjecture in the negative. Similarly when $n = 1$, it will be interesting 
to match these dense triangulated subcategories with the subcategories $\C_1^k$.

One can also try to investigate some properties of these Grothendieck groups. For instance, are they
finitely generated? are they torsion-free?  These are some questions which we think merit further study in this direction.

\subsection{Euler characteristics and the lattice of triangulated subcategories}

It would be interesting to find an Euler characteristic defined on $\C_n$ that is not a multiple of $\chi_n$ 
(see equation (\ref{E:chi_n})). Such an Euler characteristic function 
might tell something new about $K_0(\C_n)$.

Recall that $l_n$ was defined to the generator of the image of $\phi_n: K_0(\C_n) \rar \ints$. It is clear that $l_n = 1$ if $V(n)$ exists. In general, $l_n$
can be a very large integer and not much is known. 
The following is a conjecture of Ravenel (personal communication): If $(p, f) \ne (2, 1)$, then $p^f$ divides 
$l_{(p-1)f+1}$.

We have seen that the triangulated subcategories $\C_n^k$ are sandwiched between $\C_{n+1}$ and $\C_n$. Is the same true for all triangulated subcategories 
of $\F_p$?


 
\chapter{Triangulated subcategories of perfect complexes}

In this chapter we work in the derived categories of rings.  
There are many beautiful constructions of this category; see \cite{wei} for the classical approach, or \cite{hoveymodel} for a
model category theoretic approach.  We briefly review some preliminaries on the derived category $D(R)$ of a commutative ring $R$. 
It is obtained from the category of unbounded chain complexes of $R$-modules and chain maps 
by  inverting the \emph{quasi-isomorphisms} (maps that induce an isomorphism in 
homology). $D(R)$ is a tensor triangulated category with the derived
tensor product as the smash product and the ring $R$ (in degree $0$) as the unit object. 
It is a standard fact that the small objects of  $D(R)$ are precisely those complexes that are quasi-isomorphic to 
\emph{perfect complexes} (bounded chain complexes of finitely 
generated projective $R$-modules); see \cite[Prop. 9.6]{ch} for a nice proof of this fact. 
It follows from \cite[Corollary 10.4.7]{wei} that the full subcategory of small objects in $D(R)$ is 
equivalent (as a triangulated category) to the chain homotopy category of perfect complexes.
The latter will be denoted by $D^b(\proj\,R)$ and it provides a nice framework for studying small objects.
The full subcategory of small objects in $D(R)$ can also be characterised as the thick subcategory generated by $R$; see \cite[Prop. 9.6]{ch}.

In this chapter our goal is to classify the triangulated subcategories of perfect complexes over some rings (PIDs, self-injective 
noetherian rings and more generally Artinian rings) using Thomason's recipe. We begin  with some classical algebraic $K$-theory of rings
and show how this is related to the classification of dense triangulated subcategories of $D^b(\proj\,R)$. We then introduce a
theorem of Hopkins and Neeman which classifies the thick subcategories of perfect complexes over noetherian rings.
Finally in section \ref{se:classifications} we carry out our computations of the Grothendieck groups of thick subcategories
to obtain classifications of triangulated subcategories of perfect complexes over the aforementioned rings.

\section{Algebraic $K$-theory of rings} \label{se:algktheory}

Now we recall some classical algebraic $K$-theory of rings and connect it to the problem of classifying (dense) triangulated subcategories
of perfect complexes. 
If $R$ is any commutative ring, the \emph{$K$-group of the ring $R$}, which is denoted by $K_0(R)$, is defined to be the free abelian
group on isomorphism classes of finitely generated projective modules modulo the  
subgroup generated by the relations $[P]-[P \oplus Q] - [Q] = 0$, where $P$ and $Q$
are finitely generated projective $R$-modules.
A folklore result says that these $K$-groups are isomorphic to the Grothendieck groups of $D^b(\proj\,R)$.

\begin{prop} (well-known) \label{prop:well-known} If $R$ is any ring (not necessarily commutative), then there is a natural isomorphism of abelian groups
   \[ K_0(R)\cong K_0(D^b(\proj\,R)). \]
\end{prop}
\noindent
\begin{proof}(sketch) Let $A$ denote the free abelian group on the isomorphism classes of perfect complexes 
in the derived category of $R$ and let $B$ denote the free abelian group on isomorphism classes of
finitely generated projective $R$-modules. Since the Grothendieck groups under consideration are quotients of these
free groups, we define maps on $A$ and $B$ that descend to give the desired bijections.
Define $f: A \rightarrow B$ by $f(\la X \ra) = \sum_{i \in \mathbb{Z}} (-1)^i \la X_i \ra$ and
 $g: B \rightarrow A$ by $g(\la M \ra) = \la M[0]\ra$. Now one can easily verify that
these maps descend to the Grothendieck groups and that the descended maps are inverses of each other.
\end{proof}

The importance of this folklore result, for our purpose, can be seen from the following observation.

\begin{rem} This folklore result, along with the theorem \ref{main} of Thomason, connects the 
the problem of classifying dense triangulated subcategories in $D^b(\proj\,R)$ with the
$K$-theory of $R$. More precisely, there is a 1-1 correspondence between the subgroups of $K_0(R)$ and 
the dense triangulated subcategories of $D^b(\proj\,R)$. So this leads us naturally to algebraic $K$-theory -- a subject that has
been extensively studied.
\end{rem}

For the remainder of this section, we review some well-known computations of $K$-groups of rings that will be relevant to us.

\begin{example} \label{ex:K-groups}
If $R$ is any commutative local ring, then it is a fact that every finitely generated projective $R$-module is free. Thus
the monoid $\proj\,R$ (the category of finitely generated projective $R$-modules) is equivalent to the  monoid consisting of whole numbers.
The Grothendieck group of the later is clearly isomorphic to $\mathbb{Z}$. Thus for local rings, 
\[ K_0(D^b(\proj \, R) \cong K_0(R) \cong \mathbb{Z}. \]
Similarly, it is easy it to see that if $R$ is any principal ideal domain, then $K_0(R) \cong \ints.$
\end{example}

Computing these $K$-groups in general can be very hard and it is one of the central problems in algebraic $K$-theory.

Now we state some results on the $K$-groups of some low  dimensional commutative rings.
These results are of interest because they connect $K$-groups and classical Picard groups of rings.
Before we state the result, we need to recall the definition of the Picard groups of rings.

\begin{defn} The \emph{Picard group}, $\Pic(R)$ of a ring $R$, is defined to be the group of isomorphism classes 
of $R$-modules 
that are invertible under the tensor product. Similarly the \emph{Picard group} of $D(R)$, denoted by $\Pic(D(R))$, 
is the group of isomorphism classes of objects in $D(R)$ that are invertible under the derived tensor product. In both 
these cases, note that the ring $R$ acts as the identity element.
\end{defn}

\begin{thm}\cite{wei1} \label{thm:K_0picard} Let $[\Spec(R), \mathbb{Z}]$ denote the additive group of continuous functions from $\Spec(R)$ to the ring of integers
with the discrete topology. Then the following holds:

\begin{enumerate}
\item For every $0$-dimensional commutative ring $R$, $K_0(D^b(\proj\,R)) \cong [\Spec(R), \mathbb{Z}].$

\item For every $1$-dimensional commutative noetherian ring, $K_0(D^b(\proj\,R)) \cong \Pic(D(R)).$
\end{enumerate}
\end{thm}

\begin{proof} The first statement is Pierce's theorem; see \cite[Theorem 2.2.2]{wei1}. The second statement can be seen as a corollary of a theorem due to  
Fausk \cite{fau}: There is a natural split short exact sequence (for any commutative ring),
\[ 0 \rightarrow \Pic(R) \rightarrow \Pic(D(R)) \rightarrow [\Spec(R), \mathbb{Z}] \rightarrow 0.\]
Therefore $\Pic(D(R)) \cong \Pic(R) \oplus [\Spec(R),\mathbb{Z}].$ It is shown in \cite[Corollary 2.6.3]{wei1}
that 
\[K_0(D^b(\proj\,R) \cong  \Pic(R) \oplus [\Spec(R),\mathbb{Z}].\] 
So the second statement follows by combining these two results.
\end{proof}

With these and related results from algebraic $K$-theory, one can compute the $K$-groups of various families of rings and that will help
understand the dense triangulated subcategories of perfect complexes over all those rings. However, in order to classify all 
triangulated subcategories of perfect complexes, one has to compute the Grothendieck groups of the thick subcategories of these 
complexes.

\noindent
\section{Hopkins-Neeman thick subcategory theorem.} \label{se:hopkinsneeman}

In this section we state a theorem due to Hopkins and Neeman \cite{Ne}
which classifies the thick subcategories of perfect complexes over a noetherian ring.
We begin with some fundamental definitions.

\begin{defn} \cite{Ne} Given a complex $X$ in $D(R)$, define the \emph{support} of $X$, denoted by 
$\Supp(X)$,  to be the set $\{ p \in \Spec(R) : X \otimes^L K(p) \ne 0 \}$, where $K(p)$ is the fraction field
of the domain $R/p$.
\end{defn}

The support of a perfect complex has several other useful characterisations which we collect in the
next lemma. 

\begin{lemma}\cite{mps} Let $X$ be a perfect complex over a noetherian ring. Then we have the following equivalent characterisations of its support.

(a) $ \Supp(X) = \{ p \in \Spec(R): X \otimes_R R_p \ne 0 \}$.

(b) $ \Supp(X) = \{ p \in \Spec(R): X \otimes_R E(R/p) \ne 0 \}$.

(c) $ \Supp(X) = \bigcup_i \Supp\, H_i(X)$.
\end{lemma}

\begin{proof} Part (a) follows from \cite[Propositions 6.1.7, 9.3.2]{mps}. Part (b) is also 
immediate from \cite[Proposition 6.1.7]{mps} because in this proposition they show that $\la E(R/p)\ra = \la K(p) \ra$.
Finally for part (c), note that $p \in \Supp(X)$ if and only if $X \otimes R_p \ne 0$ by part (a). Then we have
the following implications.
\[X \otimes R_p \ne 0  \Leftrightarrow H_* (X \otimes R_p) \ne 0 
                     \Leftrightarrow H_*(X) \otimes R_p \ne 0 
                     \Leftrightarrow p \in \Supp H_i(X) \text{\;for some $i$}. \]
So this proves part (c).
\end{proof}

The following corollary is immediate from part (c) of the above lemma.

\begin{cor} If  $X$ is a perfect complex over a noetherian ring, then $\Supp(X)$ is a Zariski closed set.
\end{cor}

\begin{defn}
A subset of $\Spec(R)$ is said to be \emph{closed under specialisation} if it
is a union of closed sets under the Zariski topology. Equivalently, and more explicitly, 
a subset $S$ of $\Spec(R)$ is specialisation closed if whenever a prime ideal $p$ is in $S$, then
so is every prime ideal $q$ that contain $p$.
\end{defn}

Now we are ready to state the celebrated thick subcategory theorem of Hopkins and Neeman.

\begin{thm}\cite{Ne} If $R$ is any noetherian ring, then there is a natural order preserving bijection between
the sets,
\begin{center}
\{thick subcategories $\A$ of $D^b(\proj \, R)$\}
\begin{center}
$f\downarrow \ \ \uparrow g$
\end{center}
\{subsets $S$ of $\Spec(R)$ that are  closed under specialisation\}.
\end{center}
The map $f$ sends a thick subcategory $\A$ to  $\bigcup_{X \in \A} \Supp(X)$,
and the map $g$ sends a specialisation-closed subset $S$ to the thick subcategory  $ \T_S:= \{ X \in D^b(\proj \, R): \Supp(X) \in S \}$. 
\end{thm}

\noindent
\emph{Historical Note:} This theorem was first announced by Hopkins in \cite{Ho} where he claimed that this result is true for
all commutative rings. Later Neeman pointed out a mistake in Hopkins's proof and gave a correct proof that works for noetherian rings. 
He also gave a counter example \cite[Example 4.1]{Ne} which shows that the above thick-subcategory theorem breaks if the ring is not noetherian.

The following corollary is immediate from the above theorem.
\begin{cor}  Every thick subcategory of $D^b(\proj\,R)$ is a thick ideal.
\end{cor}

\section{Classification of triangulated subcategories} \label{se:classifications} 

We will now apply Thomason's recipe to classify the triangulated subcategories of perfect complexes over some commutative
noetherian rings: Principal ideal domains, self-injective noetherian rings and Artin rings.

\noindent
\subsection{Principal ideal domains}

We first set up a few notations and recall some definitions and basic facts about PIDs.

For any element $x$ in $R$, we define the mod-$x$ Moore complex $M(x)$, in analogy with the stable homotopy
category of spectra, to be the cofibre of the self map (of degree 0)
\[  R \stk{x} R \] 
in $D(R)$.

Now we recall the notion of the length of a module. For an $R$-module $M$, a chain
\[M = M_0 \supsetneq M_1 \supsetneq M_2 \supsetneq \cdots \supsetneq M_r = 0\]
is called a composition series if each $M_i/M_{i+1}$ is a simple module (one that does not have any non-trivial submodules).
The length of a module, denoted by $l(M)$, is defined to be the length of any composition series
of $M$. The fact that this is well-defined is part of the Jordan-Holder theorem.

The function $l(-)$ is an additive function on the subcategory of $R$-modules which have finite length, i.e., if 
\[ 0 \rightarrow M_1 \rightarrow \cdots \rightarrow M_s \rightarrow 0\]
is an exact sequence of $R$-modules of finite length, then
\[ \sum_i (-1)^i l(M_i) = 0.\]
Also note that when $R$ is a PID, every finitely generated torsion module has
finite length. (This can be seen as an immediate consequence of the structure theorem for finitely generated 
modules over a PID.) 

Finally, we need the following easy exercise. If $p$ and $q$ are two distinct (nonzero) prime elements in a PID $R$, then
for any $i \ge 1$,
\[  
R/(p^i) \otimes_R R_{(q)} =
 \left\{
   \begin{array}{ll}
      R/(p^i) & \mbox{when $p=q$},\\
      0       & \mbox{when $p \neq q$}.
   \end{array}
 \right.
\]

Now we are ready to compute the Grothendieck groups for thick subcategories  of perfect complexes over a PID.
Given a subset $S$ of $\Spec(R)$ that is closed under specialisation, the thick subcategory that 
corresponds to the subset $S$ (under the Hopkins-Neeman bijection) will be denoted by $\T_S$.

\begin{thm} \label{thm:K_0PID}
Let $R$ be a PID and let $S$ be a specialisation closed subset of $\Spec(R)$. Then we have the following.
\begin{enumerate}
\item If $S = \Spec(R)$, then $K_0(\T_S)$ is an infinite cyclic group generated by $R$. 
\item  If $S \ne \Spec(R)$, then $K_0(\T_S)$ is a free abelian group on the Grothendieck classes of the Moore complexes $M(p)$, for 
$p \in S$.
\end{enumerate}
\end{thm}

\begin{proof} 
The first part follows from the fact that every finitely generated projective module over a PID is free;
consequently its Grothendieck group is an infinite cyclic group (see proposition \ref{prop:well-known} and example \ref{ex:K-groups}).
For the second part, first note that a specialisation closed subset $S \ne \Spec(R)$ is a subset of maximal ideals in $R$ 
(because non-zero prime ideals in a PID are also maximal). For each prime element $p$ in $S$, define an Euler characteristic function
$\lambda_p : \T_S  \rightarrow \mathbb{Z}$ by 
\[\lambda_p(X) := \sum_i (-1)^i \; l[H_i(X) \otimes_R R_{(p)}].\]
(Since $S$ does not contain $(0)$, it follows that $H_*(X)$ is a torsion $R$-module. Therefore $\lambda_p(-)$ is a well-defined
Euler characteristic function.)

Also, since $\lambda_p(M(q)) = \delta_p^q$,  the Euler characteristic map 
\[\bigoplus_{p \in S} \lambda_p :  K_0(\T_S) \rightarrow \bigoplus_{p \in S} \mathbb{Z}\]
is clearly surjective.

To see that this map is injective, it suffices to show that every complex
in $\T_S$ can be generated by the set $\{M(p) : p \in S \}$ using cofibre sequences (see proposition \ref{prop:K_0FT}). We do this by induction 
on $\sum_i l[H_i(-)]$. 
If $X \in \T_S$ is such that $\sum_i l[H_i(X)] = 1$, then there exists an integer $j$ such that 
$H_i(X) = 0$ for all $i \neq j$, and $H_j(X) = R/(p)$ for some prime  $p$ in $R$ (because every 
simple module over a PID is of the form $R/(p)$ for some prime  $p$). 
Such an $X$ is clearly quasi-isomorphic to $M(p)$. Now consider an $X \in \T_S$ for which  $\sum_i l[H_i(X] > 1$.
Without loss of generality we can assume that up to suspension $X$ is of the form 
\[ \cdots 0 \rightarrow P_k \rightarrow \cdots \rar  P_1 \rightarrow P_0 \rightarrow 0 \cdots\]
with $H_0(X) \neq 0$ (otherwise, we can replace $X$ with a quasi-isomorphic complex in $\T_S$ which has 
this property). Now pick a non-zero element in $H_0 (X)$
and represent it with a cycle $t$. Since $H_0(X)$ is a torsion module, there exists a prime $p$ and a positive integer $k$ such that $p^k t=0 $ in homology. Replacing 
$t$ with $p^{k-1}t$, we can assume that $pt=0$ in homology, which  means $pt$ is a boundary. So there is an element $y \in P_1$ which maps under the differential
to $pt$. Consider diagram \ref{fig:killingahomologyclass}
\begin{figure}
\[
\xymatrix{ \vdots \ar[d] \ar[r] & \vdots \ar[d] \\
           0 \ar[d] \ar[r] & P_2  \ar[d]       \\
           R \ar[d]_p \ar[r]^a & P_1  \ar[d]^c \\
           R \ar[d] \ar[r]^b & P_0  \ar[d]     \\
           0  \ar[r] & 0        
           }
\]
\caption{Killing a homology class in the Hurewicz dimension}
\label{fig:killingahomologyclass}
\end{figure}
where $a(1)= y$, $c(y) = pt$, and $b(1)=t$. This diagram shows a chain map 
between the two complexes in $\T_S$ such that the induced map in homology in dimension $0$ 
sends $1$ to $t$ (by construction). So the class $t$ is killed. Now if we extend this morphism to a triangle
 ($M(p) \rightarrow X \stk{d} Y \rightarrow \Sigma M(p)$)
and look at the long exact sequence in homology, it is clear that 
$X$ and $Y$ have the same homology in all dimensions except dimension 0.
In dimension $0$, part of the long exact sequence gives a short exact 
sequence $ 0 \rightarrow R/(p) \rightarrow H_0(X) \rightarrow H_0(Y) \rightarrow 0$.
Therefore $l[H_0(Y)] = l[H_0(X)] -1$. By induction hypothesis we know that 
$Y$ can be generated by the set $\{M(p): p \in S \}$ using cofibre sequences. The above exact triangle then
tells us that $X$ can also be generated using cofibre sequences in this way.
So we are done.
\end{proof}

\begin{rem} In the next chapter, we will recover the above calculations using Krull-Schmidt decompositions; see theorem \ref{th:main} and theorem \ref{th:mainalg}
for details.
\end{rem}

\begin{cor} If $R$ is any PID, then $\Pic(D(R)) \cong \ints.$
\end{cor}
\begin{proof} Since $K_0(D^b(\proj\,R)) \cong \ints$, the corollary follows by invoking  theorem \ref{thm:K_0picard}.
\end{proof}

\noindent
\emph{Classification of the triangulated subcategories of $D^b(\proj\,R)$ when $R$ is a PID.} There are two families of triangulated subcategories:
triangulated subcategories that correspond to $S = \Spec(R)$ and  the ones that correspond to $S \ne \Spec(R)$ (subsets of maximal ideals). \\
\noindent 
1.  $S = \Spec(R)$:   Consider the Euler characteristic function 
\[\chi(X) = \sum_{-\infty}^{\infty} (-1)^i \dim_F \{H_i(X) \otimes_R F\},\]
where $F$ is the field of fractions of our domain $R$. For every integer $k$, we define 
\[ D_k = \{ X : \chi(X) \equiv 0 \! \!\mod{k} \}. \]
These are all the triangulated subcategories that are dense in $D^b(\proj\,R)$.\\
\noindent 
2. $S \ne \Spec(R)$:  Given such a subset $S$
   and a subgroup $H$ of $\bigoplus_{p \in S} \mathbb{Z}$, we define 
\[\T(S,H) = \{X \in \T_S: \left(\oplus_{p \in S} \lambda_p \right)(X) \in H \}. \] 
These are all the triangulated subcategories that are dense in $\T_S$.

It is clear from theorem \ref{thm:K_0PID} and theorem \ref{main} that every triangulated subcategory of $D^b(\proj\,R)$ is one of these
two types. 

Here is an interesting consequence of the above theorem. 

\begin{cor} Let $X$ and $Y$ be perfect complexes over a PID. Then $Y$ can be generated from
$X$ using cofibrations if and only if
\begin{itemize}
\item $\Supp(Y) \subseteq \Supp(X)$,  and
\item  If $(0) \in \Supp(X)$, then $\lambda_0(X)$ divides $\lambda_0(Y)$; otherwise, 
  $\lambda_p(X)$ divides $\lambda_p(Y)$ for all $p \in \Supp(X)$.
\end{itemize}
\end{cor}

The above classification should be compared with the classification theorem for triangulated subcategories of finite spectra.
The comparison can be made more formal once we  make the following definition.
Given a triangulated category $\A$, define the \emph{triangular lattice} $\T(\A)$ of $\A$ to be the lattice of all triangulated subcategories 
of $\A$. The meet of a set of triangulated subcategories is given by intersection,
and the join is obtained by taking the smallest triangulated subcategory containing the union of all the triangulated subcategories.
Then we have the following proposition.

\begin{prop} Let $\F$ denotes the category of finite spectra. Then there is a sub-lattice of $\T(\F)$ that is isomorphic to 
$\T(D^b(\proj\,\mathbb{Z}))$.
\end{prop}
\begin{proof} 
In view of theorem \ref{main} which gives an order preserving bijection between the subgroups of the Grothendieck
group and the dense triangulated subcategories, we will be done if we can show:

(a) There is a sublattice of thick subcategories of $\F$  that is isomorphic to the lattice of thick subcategories of 
   $D^b(\proj\,\ints)$; 
(b) The thick subcategories that correspond to each other in the above isomorphism have isomorphic Grothendieck groups.

For (a), map the thick subcategory $\T_S$ ($S$ is a subset of primes) to the full subcategory of $\F$ consisting of 
$S$-primary torsion spectra (finite spectra $X$ such that $\pi_*(X) \otimes \ints_{(q)} = 0$ for all $q$ not in $S$).
Part (b) is immediate from theorem \ref{thm:K_0PID} and proposition \ref{prop:K_0FT}.

\end{proof}

\noindent
\subsection{Self-injective noetherian rings}

We will now classify the triangulated subcategories of $D^b(\proj\,R)$ when $R$ is a commutative self-injective and noetherian.
Recall that a  commutative ring $R$ is \emph{self-injective} if $R$ is injective as a module over itself. In
other words, the injective dimension of $R$ is $0$. The importance of 
these rings comes from the fact that the notion of injectivity and projectivity are the same 
in the category of finitely generated $R$-modules.

\begin{example} $K[x]/(x^2)$ when $K$ is a field, the ring $\mathbb{Z}/m$, finite dimensional group algebras over a field, or more
generally finite dimensional Hopf algebras with antipode over a field \cite[Corollary 3.1.11]{be} are some examples
of self injective noetherian rings.  If $R$ is any PID and $a$ is a non-zero non-unit in $R$, then
$R/(a)$ is also self-injective.
\end{example}

It is a fact that a self-injective noetherian ring is Artinian. In particular, it is zero dimensional (every prime ideal is maximal) 
and has only finitely many prime ideals. This implies that every subset of $\Spec(R)$ is closed under specialisation. So by the
Hopkins-Neeman bijection, the non-trivial thick subcategories of  $D^b(\proj\,R)$ are given by $\T_S = \{X \in  D^b(\proj\,R) | \Supp(X)\subseteq S \}$,
where $S$ is an arbitrary non-empty subset of $\Spec(R)$.

Now we state some lemmas that  will be used in computing the Grothendieck groups of these thick subcategories.

\begin{lemma}\cite{cm} Every injective module $M$ over a noetherian ring
$R$ is a direct sum of indecomposable injective modules. In fact,
\[ M \cong \bigoplus_{p \in Ass(M)} E(R/p)^{\mu_p} \]
where $E(R/p)$ denotes the injective hull of $R/p$ and $\mu_p := \dim_{K(p)} \Hom_{R_p} (K(p)$, $ M_p)$.
\end{lemma}

\begin{lemma} \label{le:coolbutsad} Let $R$ be a self-injective noetherian ring. If $p$ and $q$ are any two distinct prime ideals in $R$,
then $\Hom_R (E(R/p), E(R/q))=0$.
\end{lemma}

\begin{proof} Since $\Supp(E(R/p)) = \{ p\}$, it follows that for all primes $s$ in $R$,
\[\Hom_{R_s}(E(R/p)_s, E(R/q)_s) = 0.\]
Therefore $\Hom_R (E(R/p), E(R/q)) = 0$. Since $E(R/p)$ (viewed as a complex concentrated in dimension zero) is a perfect complex.
It follows that $\Hom_{D^b(\proj\,R)} (E(R/p), E(R/q))$ is also zero.
\end{proof}

\begin{prop} Let $R$ be a self-injective noetherian ring and let $\T_S$ be the thick subcategory corresponding to a
subset $S$ of $\Spec(R)$. Then the Grothendieck group
 $K_0(\T_S)$ is a free abelian group on the Grothendieck classes $\{ E(R/p)\}_{p \in S}$. 
\end{prop}

\begin{proof} Since $\Supp E(R/p)=\{p\}$, it is clear that $\{E(R/p): p \in S \} \subseteq \T_S$. 
We will now show that $\T_S$ can be generated as a triangulated category by these injective hulls $E(R/p)$.
Consider any perfect complex $X$ in $\T_S$ and recall that over self-injective noetherian rings
a finitely generated projective module is also injective. So the complex $X$ is of the form 
\[\cdots 0 \rightarrow X_1 \stk{\delta_1} \cdots \stk{\delta_{k-1}} X_k \rar \cdots \rightarrow 0 \cdots\]
such that for each $i$, $X_i \cong \bigoplus_p E(R/p)^{\mu_{p_i}}$ (the direct sum runs over a subset of primes).
From lemma \ref{le:coolbutsad}, we can infer that the differentials in this complex map each $E(R/p)$ summand into $E(R/p)$. 
So the complex $X$ breaks into subcomplexes $X^p$:
\[ X = \bigoplus_{p\in \Spec(R)} X^p = \left( \bigoplus_{p\in S} X^p \right) \oplus \left(\bigoplus_{p \in S^c} X^p \right)\] 
where $X^p$ is a subcomplex of $X$ obtained by collecting all the $E(R/p)$-summands in each term. ($X^p$ can also be
thought of as being obtained from $X$ by tensoring with $R_{p}$.) 
Since  $\Supp X \subseteq S$ and $\Supp X^p \subseteq \{p\}$, it follows (from the above splitting of $X$) 
that $\bigoplus_{p \in S^c} X^p$ must be acyclic.  Therefore $X$, which is quasi-isomorphic to $\bigoplus_{p \in S} X^p$,  can be 
generated by iterated cofiberings of the injective hulls $E(R/p)$, for $p$ in  $S$. This shows that $K_0(\T_S)$ is an abelian 
group generated by $[E(R/p)]$, for $p$ in $S$. To see that it is in fact a free abelian group, note that there is a universal
Euler characteristic function on $\T_S$: 
\[ \bigoplus_{p \in S} \Lambda_p : \; K_0(\T_S) \rightarrow \bigoplus_{p \in S} \mathbb{Z},\]
where $\Lambda_p (X) = \sum_i  (-1)^i \dim_{K(p)} H_i(X \otimes K(p))$.  This Euler characteristic function is 
surjective ($\Lambda_p (E(R/q)) = \delta_p^q)$ and it is also injective (since $E(R/p)$ for $p \in S$ generates $\T_S$). 
So that completes the proof of this proposition.

\end{proof}

\noindent
\emph{Classification of the triangulated subcategories of $D^b(\proj\,R)$.} 
 For every subset $S$ of $\Spec(R)$ and every subgroup $H$ of $\bigoplus_{p \in S} \mathbb{Z}$, we define
\[\T(S,H) = \{X \in \T_S: (\oplus_{p \in S} \Lambda_p )(X) \in H \}. \] 
By theorem \ref{main}, it is now clear that these are all the triangulated subcategories of $D^b(\proj \, R)$.

\begin{rem} Now consider two subsets $A$ and $B$ of $\Spec(R)$ such that $A \subseteq B$. Then we have the 
inclusion $\T_A \hookrightarrow \T_B $ of thick subcategories. It can be easily verified 
that the induced map on the Grothendieck groups $K_0(\T_A) \rightarrow K_0(\T_B)$ is injective. In fact, 
for each prime $p$ in $R$, the Grothendieck class of $E(R/p)$ in $K_0(\T_A)$ is mapped to the Grothendieck
class of $E(R/p)$ in $K_0(\T_B)$. So that gives a short exact sequence of abelian groups
\[0 \rightarrow K_0(\T_A) \rightarrow K_0(\T_B) \rightarrow K_0(\T_B/\T_A) \rightarrow 0,\]
which implies that $K_0(\T_B/\T_A)$ is a free abelian group of rank $|B-A|$. 
\end{rem}

Having computed the Grothendieck groups, let us examine the ring structures on them.
(Since the thick subcategories of perfect complexes are also thick ideals, it makes sense to talk about the
ring structure on their Grothendieck groups.) Fix a subset $S \subseteq \Spec(R)$ and look at the associated 
thick subcategory $\T_S$. Note that the ring structure on $K_0(\T_S)$ will be completely determined by the
product on the Grothendieck classes of the free generators $E(R/p)$ for $p \in S$. 
Now we recall some facts about the injective hulls $E(R/p)$ over self-injective noetherian rings.
\begin{center}
 $E(R/p) \otimes E(R/p) \cong E(R/p)$, \\
 $E(R/p) \otimes E(R/q) \cong 0$ \hspace{3 mm} if $p \neq q$.
\end{center}
We have seen that every complex in $\T_S$ is quasi-isomorphic to one that is built out of
the injective hulls $E(R/p)$, for $p \in S$, using cofibre sequences. So this (together with the
above facts) will tell us that $\bigoplus_{p \in S} E(R/p)$ will be a unit for the smash product 
in $\T_S$. Therefore the Grothendieck ring $K_0(\T_S)$ is isomorphic to  
the direct product $\prod_{p \in S} \mathbb{Z}$. Now recall that
the dense triangulated ideals in a thick subcategory correspond precisely to the ideals in the Grothendieck
ring of the thick subcategory. We conclude that every dense triangulated subcategory of $\T_S$ is a triangulated
ideal if and only if $S$ is a one-point space.

This discussion can be summarised in the following theorem.

\begin{thm} Let $R$ be a self-injective noetherian ring. Then the thick subcategories $\T_S$ 
 of $D^b(\proj \, R)$ are closed-symmetric monoidal categories with $\bigoplus_{p \in S} E(R/p)$ as the unit for
the derived tensor product. The Grothendieck ring of $\T_S$ is isomorphic to $\prod_{p \in S} \mathbb{Z}$,
and if  $A \subseteq B$ are subsets of $\Spec(R)$, then the induced map on the Grothendieck groups is
an injective ring homomorphism and $K_0(\T_B/\T_A)$ is a free abelian group of rank $|B-A|$. Further 
every dense triangulated subcategory of $\T_S$ is tensor triangulated category if and only if $S$ is a 
one-point space.
\end{thm}

\subsection{Product of rings: Artin rings} We now  address the following question.\\
\noindent
\textbf{Question:} Suppose a commutative ring $R$ is a direct product of rings. Let us say
\[R \cong R_1 \times R_2 \times \cdots \times R_k,\]
and suppose that we have a classification of all the triangulated subcategories of $D^b(\proj\,R_i)$,
for all $i$. Using this information, how can we get a classification of all triangulated subcategories of $D^b(\proj\,R)$?  
More generally, one can ask how the spectral theory of $R$ and the spectral theories of the rings $R_i$ are related.

Before we go further, we remark that by the obvious induction, it suffices to consider only two components 
($R \cong R_1 \times  R_2)$. We collect some standard facts abouts products of triangulated categories
that will be needed.

\begin{lemma} \label{products} Let $\T_1$ and $\T_2$ be triangulated categories. Then we have the following.
\begin{itemize}
\item The product category $\T_1 \times \T_2$ admits a triangulated structure that has the following 
universal property: Given any triangulated category $\D$ and some triangulated functors 
$\T_1 \leftarrow \D \rightarrow \T_2$, there exist a unique triangulated functor $\F: \D \rightarrow \T_1 \times \T_2$ 
making the following diagram of triangulated functors commutative. 
\[
\xymatrix{
& \D \ar[dl] \ar@{.>}[d]^F \ar[dr] &\\
\T_1 & \T_1 \times \T_2 \ar[l]^{\pi_1} \ar[r]_{\pi_2} & \T_2
}
\]

 \item Every thick (localising) subcategory of $\T_1 \times \T_2$ is of the form $\B_1 \times \B_2$, where
$\B_i$ is a thick (localising) subcategory in $\T_i$.

\item If both $\T_1$ and $\T_2$ are essentially small triangulated 
(tensor triangulated) categories, then $K_0(\T_1 \times \T_2) \cong K_0(\T_1) \times K_0(\T_2)$ as groups (rings).
\end{itemize}
\end{lemma}

\begin{rem} In contrast with the thick subcategories, not every triangulated subcategory of
$\T_1 \times \T_2$ is of the form $\A_1 \times \A_2$, where $\A_i$ is a triangulated subcategory 
of $\T_i$. This is clear from the third part of the above lemma because not every subgroup of
the product group $K_0(\T_1) \times K_0(\T_2)$ is a product of subgroups. 
\end{rem}

Now we relate the category of perfect complexes over $R$ and those over the rings
$R_i$.

\begin{prop} There is a natural equivalence of triangulated categories,
\[ D^b(\proj\,R) \simeq D^b(\proj\,R_1) \times D^b(\proj\,R_2).\]
In particular $K_0[D^b(\proj\,R)]  \cong K_0 [D^b(\proj\,R_1)] \times K_0[D^b(\proj\,R_2)].$
\end{prop}

\begin{proof} Note that every module $M$ over $R_1 \times R_2$ is a direct sum (in the category of
$R_1 \times R_2$ - modules) 
\[ M \cong P_1 \oplus P_2,\]
where $P_i$ is an $R_i$-module: Take $P_1 = \la (1,0) \ra M$ and $P_2 = \la (0,1) \ra M$, where $P_i$ is also regarded as 
an $R_1 \times R_2$ module via the projection maps $R_1 \times R_2 \rightarrow R_i$. Now one can verify easily that 
this decomposition is functorial and sends finitely generated (projective) modules to finitely generated (projective) 
modules. Therefore every complex in $D^b(\proj\,R_1 \times R_2)$ splits as a direct sum of two complexes, one each in 
$D^b(\proj\,R_i)$. Conversely, given a pair of complexes $(X_1, X_2)$ with $X_i$ in $D^b(\proj\,R_i)$, their direct sum 
$X_1 \oplus X_2$ is clearly a complex in $D^b(\proj\, R_1 \times R_2)$. Now it can be verified  that these two functors 
establish the desired equivalence of triangulated categories.  The second statement follows immediately from the 
lemma \ref{products}.
\end{proof}

So in view of this proposition, the problem of classifying triangulated subcategories in $D^b(\proj\,R_1 \times R_2)$
boils down to classifying triangulated subcategories of $D^b(\proj\,R_1) \times D^b(\proj\,R_2)$. So, following
Thomason's recipe, we need to classify the thick subcategories of $D^b(\proj\,R_1) \times D^b(\proj\,R_2)$ and compute their Grothendieck groups. 
This is given by Lemma \ref{products}: Every thick subcategory $\T$ of $D^b(\proj\,R_1) \times D^b(\proj\,R_2)$
is of the form $\T_1 \times \T_2$ where $\T_i$ is thick in $D^b(\proj\,R_i)$, and 
$K_0(\T) \cong K_0(\T_1 \times \T_2) \cong K_0(\T_1) \times K_0(\T_2)$.

\subsubsection{Artin rings}
\noindent
We will now apply these ideas to Artin rings. 
Recall that a ring is said to be Artinian if every descending chain of ideals terminates. As mentioned
earlier, every self-injective noetherian ring is Artin, but the converse is not true: the ring
\[ k[x, y]/(x^2, xy, y^2)\]
is an Artin ring whose injective dimension is infinite.
Recall that Artin rings are precisely the zero dimensional noetherian rings and they have the following structure theorem. 
\begin{thm}\cite{am}
Every Artin ring $R$ is isomorphic to a finite direct product of 
      Artin local rings. Moreover, the number of local rings that appear in this
      isomorphism is equal to the cardinality of $\Spec(R)$.
\end{thm}

Thus  $R \cong \prod_{i=1}^{n} R_i$, where each $R_i$ is an Artin local ring
($n$ is the cardinality of $\Spec(R)$). We have seen above that there is an equivalence 
of triangulated categories,
\[ D^b(\proj\,R) \cong \prod_{i=1}^{n} D^b(\proj\,R_i).\]
So by the above discussion, we just have to compute the Grothendieck groups of the thick subcategories of $D^b(\proj\,R_i)$.
But since each $R_i$ is an Artinian local ring, $\Spec(R_i)$ is a one-point space. This implies (by the Hopkins-Neeman theorem) that
the only non-zero thick subcategory is $D^b(\proj\,R_i)$ itself, whose Grothendieck group is well-known to be infinite cyclic.  
Therefore we have,
 \[ K_0(\T_S) \cong \bigoplus_{p_i \in S} K_0(D^b(\proj\,R_i)) \cong \bigoplus_{p_i \in S} \mathbb{Z}. \]
The universal Euler characteristic function that gives this isomorphism is $\oplus_{p_i \in S} \Lambda_{p_i}$ where 
$\Lambda_{p_i} (X) = \sum_{t=-\infty}^{\infty} (-1)^t \dim_{R_i/p_i} H_t (X \otimes R_i/p_i) $.

So we have finally proved the promised classification result that generalises the analogous result on self-injective 
noetherian rings. 

\begin{thm} Let $R$ be any Artin ring and let $R = \prod_i R_i$ be its unique decomposition into Artin local rings. 
For every subset $S$ of $\Spec(R)$ and every subgroup  $H$ of $\bigoplus_{p_i \in S} \mathbb{Z}$, define
\[\T(S,H) := \{X \in \T_S: (\oplus_{p_i \in S} \Lambda_{p_i} )(X) \in H \}. \] 
This is a complete list of  triangulated subcategories of $D^b(\proj \, R)$.  The thick subcategories $\T_S$ of $D^b(\proj\,R)$
are isomorphic to $D^b(\proj\, \prod_{p_i \in S} R_i)$, and therefore, 
are closed symmetric monoidal categories.  The Grothendieck ring of $\T_S$ is 
isomorphic to $\bigoplus_{p_i \in S} \mathbb{Z}$ and $\prod_{p_i \in S} R_i$ is the unit for smash product.
Further, every dense triangulated subcategory of $\T_S$ is a triangulated ideal if and only if $S$ is a
one point space.
\end{thm}

\begin{rem} It is clear from the proof that this theorem also holds whenever $R \cong \prod_{i=1}^{n} R_i$, where each of the rings 
$R_i$ has exactly one prime ideal.  
\end{rem}

We now derive some easy consequences of the above theorem.

\begin{cor}
Let $X$ and $Y$ be perfect complexes over an Artin ring. Then $Y$ can be generated from
$X$ using cofibrations if and only if
\begin{itemize}
\item $\Supp(Y) \subseteq \Supp(X)$, and 
\item $\Lambda_{p_i}(X)$ divides $\Lambda_{p_i}(Y)$ for all $p_i \in \Supp(X)$,
\end{itemize}
\end{cor}

\begin{cor} An Artin ring $R$ is local if and only if every dense triangulated subcategory of $D^b(\proj\,R)$ is  a triangulated ideal.
\end{cor}

\begin{proof} By theorem \ref{imain}, it is clear that every dense triangulated subcategory is a triangulated ideal if and only
if every subgroup of $K_0(D^b(\proj\,R)) \cong \prod_{p \in \Spec(R)} \ints$ is also an ideal. Clearly the later happens if and only if 
$|\Spec(R)| = 1$, or equivalently if $R$ is local.
\end{proof}

\noindent
\subsection{Non-noetherian rings}

In order to study the problem of classifying  triangulated subcategories in the non-noetherian case,
we need a thick subcategory theorem for $D^b(\proj\,R)$, when $R$ is a non-noetherian ring.
This is given by a result of Thomason, which is a far-reaching generalisation of 
Hopkins-Neeman theorem to schemes. Before we can state this deep result, we have to recall some terminology 
from \cite{Th}.

Let $X$ be a quasi-compact and quasi-separated scheme. The prime examples of such schemes are
the affine schemes $\Spec(R)$ and noetherian schemes.
We shall denote the derived category of the abelian category of sheaves of $\mathcal{O}_X$-modules
by $D(X)$. A \emph{strict perfect complex} on $X$ is a bounded complex of locally free $\mathcal{O}_X$-modules
of finite type. A  \emph{perfect complex} on $X$ is a complex $E$ of sheaves of $\mathcal{O}_X$-modules where
$X$ can be covered by open sets $U_{\alpha}$ such that $E|_{U_{\alpha}}$ is quasi-isomorphic to a strict perfect complex on
$U_{\alpha}$. The full subcategory of perfect complexes on $X$ will be denoted by $D(X)_{\text{parf}}$.

Let $E$ be a complex of sheaves of $\mathcal{O}_X$-modules. The \emph{cohomological support} 
of $E$ is the subspace $\Supp(E) \subseteq X$ consisting of those points $x \in X$ for which the associated stalk
complex $E_x$ of $\mathcal{O}_{X,x}$-modules is not acyclic. This can also be expressed as
\[\Supp(E) = \bigcup_{n \in \mathbb{Z}} \Supp \, H_n(E), \]
the union of supports in the classical sense of the cohomology sheaves of $E$.

Now we are ready to state the grand result of Thomason which classifies the thick subcategories of $D(X)_{\text{parf}}$.

\begin{thm} \cite{Th} \label{thm:thomason-thicksubcategory} Let $X$ be a quasi-compact and quasi-separated scheme. Then there is a bijection
between the thick ideals $\A$ of $D(X)_{\text{parf}}$
and the subsets $Y$ of $X$ which are a union of closed sets $Y_{\alpha}$ such that $X-Y_{\alpha}$ is quasi-compact.
The bijection maps a thick subcategory $\A$ to the subspace $\bigcup_{x \in \A} \Supp\,x$ and the inverse map
sends $Y$ to the thick subcategory of all perfect complexes that are acyclic off $Y$, i.e., complexes whose
support is contained in $Y$.
\end{thm}

A special case of this result is the next corollary which generalises the Hopkins-Neeman thick subcategory theorem to commutative rings 
(not necessarily noetherian). Although this corollary is known, we could not find
a proof in the literature, so we provide one.

\begin{cor} \label{cor:thicksubcatrings} Let $R$ be any commutative ring. Then there is a natural order preserving bijection between
the sets
\begin{center}
\{thick subcategories $\A$ of $D^b(\proj \, R)$\}
\begin{center}
$f\downarrow  \;\;\; \uparrow g$
\end{center}
\{subsets $S$ of $\Spec(R)$ such that $ S = \bigcup_{\alpha} V(I_{\alpha})$, where $I_{\alpha}$ is finitely generated ideal\}.
\end{center}
The map $f$ sends a thick subcategory $\A$ to the set $ \bigcup_{X \in \A} \Supp(X) $
and the map $g$ sends $S$ to the thick subcategory $\{ X \in D^b(\proj \, R): \Supp(X) \in S \}$. 
\end{cor}

\begin{proof} The affine scheme $\Spec(R)$ is both quasi-compact and quasi-separated, therefore 
the above theorem of Thomason applies to it. So throughout the proof we let $X=\Spec(R)$. 

First we want to show that $D(X)_{\text{parf}}$ is equivalent to $D^b(\proj\,R)$. Recall the well-known equivalence
\cite[Corollary 5.5]{Ha} between the abelian categories of quasi-coherent sheaves of $\mathcal{O}_X$-modules
and the category of $R$-modules. This equivalence clearly extends to the bounded homotopy categories of 
chain complexes over these abelian categories. Now over affine schemes, every perfect complex is globally quasi-isomorphic
to a strict perfect complex and therefore we need to ask the following question.  What is the image of
the full subcategory of the strict perfect complexes under this equivalence of chain homotopy categories? 
Let $E$ be any strict perfect complex. Then for each $i$, $E_i$ is a locally free $\mathcal{O}_X$-module of finite type. 
Let the corresponding $R$-module be denoted by $M_i (= \Gamma(X,E_i) )$ (every locally free sheaf is quasi-coherent). 
Since $E_i$ is locally free of finite type, for all $p \in X$, the stalk module ${E_i}_p$ is a free $\mathcal{O}_{X,p}$-module
of finite rank.
We know that the latter is isomorphic to the localised module ${M_i}_p$ \cite[Prop. 5.1(b)]{Ha}. So we conclude that ${M_i}_p$ is a free 
$R_p$-module of finite rank for all $i$. That means each $M_i$ is a finitely generated projective $R$-module.
Conversely if $M$ is a perfect complex of $R$-modules, then it corresponds to a bounded complex of quasi-coherent
sheaves of $\mathcal{O}_X$-modules. Since each $M_i$ is finitely generated and projective, all its localisations, ${M_i}_p$, are free of 
finite rank over $R_p$. As before, since ${M_i}_p$ is isomorphic to the stalk module ${E_i}_p$ over $\mathcal{O}_{X,p}$, we conclude that
$E_i$ is locally free of finite rank. So that completes the proof of the equivalence.

Next we argue that the equivalence shown in the previous paragraph respects support. Under the equivalence of the abelian 
categories of quasi-coherent sheaves of $\mathcal{O}_X$-modules and $R$-modules, for any strict perfect complex $E$,
 $H_i(E)$ will correspond to $H_i(M)$. Further $H_i(E)$ is quasi-coherent, and therefore, 
\[ H_i(E)_p \cong  H_i(M)_p\] as $\mathcal{O}_{X,p}$ ($\cong R_p$)-modules. So the above equivalence respects support.

Finally we note two things: Every thick subcategory of $D^b(\proj\,R)$ is a thick ideal \cite[Corollary 3.11.1 (a)]{Th}, and a standard exercise 
\cite[Page 12, Ex.17 (v)]{am} states that the complement of a closed subset 
$S$ of $\Spec(R)$ is quasi-compact if and only if $S = V(I)$ for some finitely generated
ideal in $R$. So that gives us the bijection as stated in the theorem.

\end{proof}

\begin{rem} The subsets of $\Spec(R)$ in this corollary which determine the thick subcategories of perfect complexes will be
called \emph{Thick supports}.  If $R$ is noetherian, then every thick support is a specialisation-closed subset (since ideals in a noetherian ring
are finitely generated), therefore the above corollary recovers the Hopkins-Neeman thick subcategory theorem.
\end{rem}

We consider the simplest possible non-noetherian rings -- non-noetherian rings which have a unique
prime ideal, e.g., $R = \ftwo[X_2, X_3, \cdots ]/(X_2^2, X_3^3, \cdots )$.

\begin{prop} Let $R$ be any commutative ring with a unique prime ideal $p$. Then every triangulated subcategory of
$D^b(\proj\,R)$ is of the form
  \[ \D_m = \{X \in D^b(\proj\,R)\; \mbox{such that}\; \Lambda(X) \equiv 0 \!\!\! \mod{m} \} \]
for some non-negative integer $m$, where 
$\Lambda(X) = \sum_{-\infty}^\infty (-1)^i \dim_{R/p} H_i(X \otimes R/p)$.
\end{prop}
\begin{proof} Note that any such $R$ is, in particular, a local ring. Therefore, 
\[K_0(D^b(\proj\,R)) (\cong K_0(R)) \cong \mathbb{Z}.\]
It is easily verified that the given Euler characteristic function gives this isomorphism. 
Moreover, since $R$ has a unique prime ideal, it is clear from corollary \ref{cor:thicksubcatrings} that there are no non-trivial
thick subcategories in $D^b(\proj\,R)$. So the dense triangulated
subcategories in $D^b(\proj\,R)$ are all the triangulated subcategories in $D^b(\proj\,R)$.
\end{proof}

\section{Questions}

\subsection{Algebraic $K$-theory for thick subcategories.}

In section \ref{se:algktheory} we have summarised a few results from algebraic $K$-theory. Those
were results about the Grothendieck groups of $D^b(\proj\,R)$. So we now ask if such results also hold for thick subcategories
of $D^b(\proj\,R)$. We ask a very specific question to make this point clear. 
It is well-known that if $J$ is the nilradical of $R$, then
\[ K_0(D^b(\proj\,R)) \cong K_0(D^b(\proj\,R/J)).\]
Since every prime ideal contains the nilradical, the quotient map $R \rar R/J$ induces a homeomorphism on prime spectra:
$\Spec(R) \cong \Spec(R/J)$. This homeomorphism of prime spectra implies that the lattice of specialisation-closed subsets of 
$\Spec(R)$ is isomorphic to that of $\Spec(R/J)$. Now if $R$ is noetherian, we can invoke the Hopkins-Neeman thick subcategory theorem to conclude that 
the same is true for the lattices of thick subcategories of perfect complexes over $R$ and $R/J$. Now the question that arises is whether
the thick subcategories that correspond to each other under this isomorphism have isomorphic Grothendieck groups?

\subsection{The Bousfield lattice for the derived category.}
Another important family of subcategories that are important in stable homotopy theory are the so called Bousfield classes. These were first
introduced in the stable homotopy category by Bousfield \cite{bob} and were later studied in the derived category by Neeman \cite{Ne}.
Working in the derived category, define the \emph{Bousfield class}  of a complex $A$, denoted $\la A \ra$, to be the
collection of all complexes $X$ such that $X \otimes^L_R A = 0$. We say that $\la A \ra \le \la B \ra$ if $\la B \ra \subseteq \la A \ra$.
Amnon Neeman \cite{Ne} classified these Bousfield classes for the derived category  
of a noetherian ring $R$.  He showed that they are in bijection with the subsets of $\Spec(R)$: given a subset $S \subseteq \Spec(R)$, the corresponding 
Bousfield class is the localising subcategory generated by the 
$K(p)$, for $p \in S$. He also gave the following counterexample in \cite{ne2} to show that this theorem breaks for 
non-noetherian rings. If
                     \[R = k[X_2,X_3, \cdots ]/(X_2^2,X_3^3,\cdots),\]
Neeman showed that  $D(R)$ has uncountably many Bousfield classes, although $R$ has only one prime ideal -- a striking contrast from the
noetherian result. This example suggests that the Bousfield lattice of $D(R)$ is very mysterious when $R$ is non-noetherian and it is not 
even clear what a reasonable conjectural classification of this lattice would be. In the non-noetherian case, it is also not known whether
the Bousfield classes form a set. (However, they do from a set when $R$ is a countable ring; see \cite{dwypal} for a short proof.)

From a slightly different but related view point, it is widely acknowledged that the stable homotopy category of spectra closely 
resembles the derived categories of non-noetherian rings. In fact, in \emph{``Brave New Commutative Algebra''}, one views the stable homotopy 
category as a derived category of modules over the brave new ring $S^0$ -- a ``non-noetherian'' ring spectrum. 
Furthermore, the derived category of any ring $R$ is known to be equivalent to the the derived category of 
$HR$-module spectra \cite{ekmm} 
($HR$ denotes the Eilenberg-Mac Lane spectrum). So in view of these results, any new structural information on the derived category may 
have topological implications. To summarise, the hope is that a thorough understanding
of $D(R)$ and its spectral theory, when $R$ is non-noetherian, might shed some new light on the stable homotopy category.
See comments in \cite[Section 4]{Ne}.



\chapter{Krull-Schmidt decompositions for thick subcategories}

Using ideas from modular representation theory of finite groups, Krause \cite{Kr} proved a Krull-Schmidt theorem for thick subcategories of  
the stable module category. More precisely, he showed that if $G$ is any finite group, then the thick ideals of $\stmod(KG)$ decompose uniquely into indecomposable thick 
ideals. In this chapter, we show that such decompositions exist in other stable homotopy categories like the derived categories of commutative
rings and the stable homotopy category of spectra.  To this end, we first generalise his definition to arbitrary triangulated categories. 

\begin{defn} \label{defn:KS} \cite{Kr} Let $\T$ denote a tensor triangulated category. A thick subcategory $\A$ of $\T$ is a 
\emph{thick ideal} if $X \wedge Y$ belongs to $\A$ for all $X \in \A$ and all $Y \in \T$.  If $\A$ is a 
thick (ideal) subcategory of $\T$, a family of thick (ideal) subcategories $(\A_i)_{i \in I}$ is a \emph{decomposition} of $\A$ if 
\begin{enumerate}
\item{the objects of $\A$ are the finite coproducts of objects from the $\A_i$, and}
\item{$\A_i \cap \A_j = 0$ for all $i \ne j$.}
\end{enumerate}
A decomposition $(\A_i)_{i \in I}$ of $\A$ is denoted by $\A= \coprod_{i \in I} \A_i$, and we say that $\A$ is \emph{indecomposable} 
if $\A \ne 0$  and any decomposition $\A = \C \amalg \D$ implies that $\C = 0$ or $\D = 0$. A \emph{Krull-Schmidt decomposition} 
of $\A$ is a unique decomposition $\A = \coprod_{i \in I} \A_i$ where all the $\A_i$ are indecomposable.  Say that the thick 
subcategories of $\T$ 
satisfy the \emph{cancellation property} if $\A \amalg \B \cong \A \amalg \C$ (for thick subcategories $\A$, $\B$, 
and $\C$ ) implies $\B \cong \C$.
\end{defn}

The above definition reminds one 
of the classical Krull-Schmidt theorem which states that any finite length module over a ring admits a unique direct sum decomposition into indecomposable modules. 
Since all thick subcategories studied in this paper consists of finite objects in some triangulated category, it is reasonable to call these
decompositions Krull-Schmidt decompositions.

After developing some necessary $K$-theoretic results for thick subcategories (section \ref{sec:ktheory-thick}), we start proving  Krull-Schmidt theorems 
for  thick subcategories of small objects in some stable homotopy categories. We begin with Krause's Krull-Schmidt theorem
for $\stmod(KG)$ in section \ref{sec:henning}. In section \ref{sec:KS-derivedcategory} we prove a Krull-Schmidt theorem for the derived categories over noetherian
rings. This result is then generalised in section \ref{sec:nsht} to noetherian stable homotopy categories. This paves the way for a Krull-Schmidt theorem for stable
module categories over some finite dimensional co-commutative Hopf algebras in section \ref{sec:KS-hopfalgebras}. Finally in the last section, we show that these
decompositions also hold in the category of finite spectra.

\section{$K$-theory for thick subcategories} \label{sec:ktheory-thick}

In this section we study how $K$-theory of triangulated categories behaves under a Krull-Schmidt decomposition. 
We begin with some lemmas that will be needed in studying $K$-theory.

It is well-known that coproducts of exact triangles are exact in any triangulated category; see \cite[Appendix 2, Prop. 10]{mar}. Under some additional hypothesis,
the next lemma gives a converse to this well-known fact.

\begin{lemma} \label{le:sweet} Let $\T$ denote a triangulated category and let 
\[ A \stk{f} B \rar C \rar \Sigma A  \quad \mbox{and} \quad A' \stk{f'} B' \rar C' \rar \Sigma A' \]
be two sequences of maps in $\T$ such that their sum $ A \amalg A' \stk{f \amalg f'} B\amalg B' \rar C \amalg C' \rar \Sigma (A\oplus A') $
is an exact triangle. If we also have
\begin{eqnarray}
\Hom(\Cone(f),C')  = 0 = \Hom(\Cone (f'),C)  \quad \mbox{and} \\
\Hom(C',\Cone(f))  = 0 = \Hom(C, \Cone(f')), 
\end{eqnarray}
then the given sequences of maps are exact triangles.

\begin{proof} Complete the maps $f$ and $f'$ to exact triangles in $\T$: 
\[ A \stk{f} B \rar \Cone(f) \rar \Sigma A  \quad \mbox{and} \quad A' \stk{f'} B' \rar \Cone(f') \rar \Sigma A'. \]
Since the coproduct of exact triangles is exact, adding these two triangles gives another triangle,  
\[ A \amalg A' \stk{f\oplus f'} B\amalg B' \rar \Cone(f) \amalg \Cone(f') \rar \Sigma (A\oplus A'). \]

We know from the axioms for a triangulated category that there is a fill-in map $H$ in the diagram below.
\[
\xymatrix{ 
A \amalg A' \ar[d]^{=} \ar[r]^{f \amalg f'}&  B \amalg B' \ar[d]^{=} \ar[r] & \Cone(f) \amalg \Cone(f') 
\ar@{..>}[d]^{H} \ar[r] & \Sigma (A \amalg A') \ar[d]^{=} \\ 
A \amalg A' \ar[r]^{f \amalg f'} & B \amalg B' \ar[r] & C \amalg C' \ar[r] &\Sigma(A \amalg A')
 }
\]
Note that three out of the four vertical maps in the above diagram are isomorphisms and therefore so is $H$; see \cite[Appendix 2, Prop. 6]{mar}.
Now the hypothesis 
\[\Hom(\Cone (f), C') = \Hom(\Cone(f'),C)= 0\] 
implies that $H = h \amalg h'$. So we have  $h \amalg h' : \Cone(f) \amalg \Cone(f') \rar C \amalg C'$ is an isomorphism.
The hypothesis 
\[ \Hom(C' , \Cone(f)) = \Hom(C, \Cone(f')) = 0\]
implies that the inverse $G$ of this isomorphism is of the form $ g \amalg g'$.
This forces both $h$ and $h'$ to be isomorphisms and hence $C \cong \Cone(f)$ and $C' \cong \Cone(f')$. 
Since exact triangles in $\T$ form a replete class, the two sequences of maps under consideration are exact triangles.
\end{proof}

\end{lemma}

\begin{lemma} \label{le:short} Let $\C$ and $\D$ be thick ideals in $\T$. If $\C \cap \D = 0$, then $\C \wedge \D = 0$.
($\C \wedge \D$ is the full subcategory of objects of the form $X \wedge Y$ where $X \in \C$ and $Y \in \D$.)
\end{lemma}
\begin{proof} Clear, since $\C \wedge \D \subseteq \C \cap \D = 0$. 
\end{proof}

We are now ready to state and prove our main theorem on $K$-theory for thick subcategories.

\begin{thm} \label{th:main} Let $\A = \underset{i \in I}{\coprod} \A_i$ be a Krull-Schmidt decomposition of a thick (ideal) subcategory 
$\A$ 
in a triangulated category $\T$. If $\A$ is thick, assume that (1) holds; if $\A$ is a thick ideal, assume any one of the following three 
conditions.
\begin{enumerate} 
\item{$\Hom(\A_i,\A_j)=0$ for all $i \ne j$. }
\item{There exist an object $S_i$ in $\A_i$ such that $S_i \wedge X = X$ for all $X \in \A_i$.} 
\item{There exists an object $S_i$ in $\T$ such that $S_i \wedge X = X$ for all $X \in \A_i$, and whenever $j \ne i$, $S_i \wedge X = 0$ 
for all $X \in \A_j$.}
\end{enumerate}
Then the inclusion functors  $\A_i \hookrightarrow \A$ give rise to an isomorphism
\[ \Xi: \underset{i \in I}{\bigoplus} \; K_0(\A_i) \cong K_0(\A). \]
\end{thm}

\begin{proof} The inclusion functors  $\A_i \hookrightarrow \A$ induce maps on the Grothendieck groups which can be assembled 
to obtain the map  $\Xi: \underset{i \in I}{\bigoplus} K_0(\A_i) \rar K_0(\A)$ that we want to show is an isomorphism. Showing that 
$\Xi$ is 
surjective is easy: note that every element in  $K_0(\A)$ is of the form $[X]$ for some $X \in \A$. Since the family of thick ideals 
$(\A_i)_{i \in I}$ is a decomposition for $\A$, we can express $X$ as a finite coproduct; $X = \coprod X_i$ with $X_i \in \A_i$. This 
gives $[X] = [\coprod X_i] = \Sigma_{i \in I} [X_i]$. The last quantity is clearly the image of $\bigoplus [X_i]$ under the map $\Xi$ 
and 
therefore $\Xi$ is surjective. 

To show that $\Xi$ is injective, we use Landsburg's criterion \ref{le:landsburg}.  Suppose $\Xi(\bigoplus [X_i]) = 0$. Then, 
as noted above,  $\Xi(\bigoplus [X_i]) = \Sigma [X_i] = [\coprod X_i]$ and hence 
$[\coprod X_i]=[0]$. This now gives, by  Landsburg's criterion, two exact triangles in $\A$: 
\begin{subequations} \label{triangle} 
\begin{gather}
 A \rar B \amalg \left(\underset{j \in I}{\coprod} X_j \right) \rar C \rar \Sigma A  \label{trianglea} \\
A \rar B \rar C \rar \Sigma A  \label{triangleb} 
\end{gather}
\end{subequations}
We want to show that $\bigoplus [X_i] = 0$ or equivalently $[X_i]=0$ for all $i \in I$. 
First assume that condition (1) holds.  Then consider the decompositions $A = \coprod A_i$, $B = \coprod B_i$, and $C = \coprod C_i$
(which exist because $\A = \underset{i \in I}{\coprod} \A_i$ ). Substituting these decompositions
in the above triangles \eqref{trianglea} \eqref{triangleb} gives
\[\coprod A_i \rar \coprod \left( B_i \amalg X_i \right) \rar \coprod C_i \rar \Sigma \coprod A_i, \]
\[\coprod A_i \rar \coprod B_i \rar \coprod C_i \rar \Sigma \coprod A_i. \]
Now our assumption (1) together with lemma \ref{le:sweet} will enable us to split these two triangles into exact triangles in $\A_i$. 
So for each $i \in I$, we get exact triangles 
\[A_i \rar B_i \amalg X_i \rar C_i \rar \Sigma A_i,\]
\[A_i \rar B_i \rar C_i \rar \Sigma A_i.\]
in $\A_i$. This implies (by Landsburg's criterion) that $[X_i] = 0$ in $K_0(\A_i)$. So $\Xi$ is injective if condition (1) holds.

Now assume that $\A$ is a thick ideal and that for each fixed $i \in I$ either (2) or (3) holds. Smash the above triangles \eqref{trianglea} \eqref{triangleb} 
with $S_i$ 
to get two exact triangles in $\T$:
\[A \wedge S_i \rar (B \wedge S_i) \amalg \left(\underset{j \in I}{\coprod}( X_j \wedge S_i) \right) \rar C \wedge S_i \rar \Sigma( A 
\wedge S_i)\]
\[A \wedge S_i \rar B \wedge S_i \rar C \wedge S_i  \rar \Sigma (A \wedge S_i)\]
It is easily seen that these triangles are in fact triangles in $\A_i$: This is trivial if condition (2) holds (because $\A_i$ are ideals and 
$S_i \in \A_i$). If condition (3) holds, write
$A = \coprod A_i$ with $A_i \in \A_i$, then $A \wedge S_i = \coprod (A_i \wedge S_i) = A_i \in \A_i$. Similarly $B \wedge S_i$ and $C 
\wedge S_i$ also belong to $\A_i$.

Now we claim that $\underset{j \in I}{\coprod} (X_j \wedge S_i) = X_i$. If (2) holds, then by the lemma \ref{le:short}, we get
$X_j \wedge S_i = 0$ whenever $i \ne j$, and since $S_i$ is a unit for  $\A_i$, $X_i \wedge S_i = X_i$. If (3) holds, this is obviously true.
So in both cases (conditions (2) and (3)), the above triangles can be simplified to obtain the following triangles in $\A_i$.
\[A \wedge S_i \rar (B \wedge S_i) \amalg X_i  \rar C \wedge S_i \rar \Sigma( A \wedge S_i),\]
\[A \wedge S_i \rar B \wedge S_i \rar C \wedge S_i  \rar \Sigma (A \wedge S_i).\]
This implies that $[X_i] = 0$ in $K_0(\A_i)$.  This show that $\Xi$ is injective, completing the proof of the theorem.
\end{proof}

Here is another crucial lemma for studying Krull-Schmidt decompositions.

\begin{lemma} \label{le:unclear} Let $\T$ be a triangulated category and let $\A$ and $\B$ be two thick (ideal) subcategories of $\T$. 
If $\;\Hom(\A, \B) = 0 = \Hom(\B, \A)$, then the full subcategory $\A \amalg \B$, consisting of objects of the form $A \amalg B$ with 
$A \in \A$ and $B \in B$, is a thick (ideal) subcategory of $\T$. 
\end{lemma}

\begin{proof} The key observation here is that every map $H: A \amalg B \rar A' \amalg B'$  in $\A \amalg \B$ is forced by the given hypothesis
to be of the form $f \amalg g$.  It is very clear that $\A \amalg \B$ satisfies the ideal condition. To see that $\A \amalg \B$ satisfies
the 2 out of 3 condition, start with a map $H$ as above and complete it to a triangle. Then we have the following diagram where the rows
are triangles in $\T$. (The bottom row is the coproduct of two triangles in $\T$.)
\[
\xymatrix{ 
A \amalg B \ar[d]^{=} \ar[r]^{f \amalg g}&  A' \amalg B' \ar[d]^{=} \ar[r] & \Cone(f \amalg g) 
\ar@{..>}[d]^{H} \ar[r] & \Sigma (A \amalg B) \ar[d]^{=} \\ 
A \amalg B \ar[r]^{f \amalg g} & A' \amalg B' \ar[r] & \Cone(f) \amalg \Cone(g) \ar[r] &\Sigma(A \amalg A')
}
\]
There exits a fill-in map $H$ which turns out to be an isomorphism as before. Therefore $\A \amalg \B$ is a triangulated subcategory.

It remains to show thickness, i.e., $\A \, \amalg \, \B$ is closed under retractions. Consider a retraction map $e: A \amalg B \rar A \amalg B$
(so $e^2=e$). Since $e = a \amalg b$, the equation $e^2 =e$ implies $(a \amalg b)^2 = a^2 \amalg b^2 = a \amalg b$. This shows
both $a$ and $b$ are retractions.
\end{proof}

Having developed all the necessary categorical tools, we now turn our attention to Krull-Schmidt decompositions.

\section{Stable module categories over group algebras} \label{sec:henning}

Consider the stable module category $\StMod(KG)$, where $G$ is a finite group and $K$ is some field. The objects of this
category are the (right) $KG$-modules and morphisms are equivalence classes of $KG$-module homomorphisms where two homomorphisms are equivalent 
if and only if their difference factors through a projective module. This category is well-known to be a triangulated category and has a well-defined 
tensor product (ordinary tensor product of $K$-vector spaces with the diagonal $G$-action) which makes it into a tensor triangulated category. The full subcategory of 
small objects in $\StMod(KG)$  is equivalent to the full subcategory consisting of finitely generated $KG$-modules and is denoted by $\stmod(KG)$. 
The main theorem of \cite{Kr} then states:

\begin{thm} \cite{Kr}
Every thick ideal $\A$ in $\stmod(KG)$ decomposes uniquely into indecomposable thick ideals $(\A_i)_{i \in I}$; $\A = \coprod_{i \in I} 
\A_i.$ 
Conversely given thick  ideals $(\A_i)_{i \in I}$ such that $\A_i \cap \A_j = 0$ for all $i \ne j \in I$, there exists a thick ideal 
$\A$ such that $\A = \coprod_{i \in I} \A_i$.  
\end{thm}

We will briefly outline how Krause arrives at this decomposition. The key idea is to consider the Bousfield localisation with respect to the localising 
subcategory  $\A^{\amalg}$ generated by $\A$. The inclusion $\A^{\amalg} \hookrightarrow \StMod(KG)$ has a right adjoint 
$e: \StMod(KG) \rar {\A}^{\amalg}$. In fact, $e(M)$ is just the fibre of the Bousfield localisation map $M \rar M_{A^{\amalg}} $. Thus for each $KG$-module 
$M$, there is a natural triangle in $\StMod(KG)$ given by $e(M) \stk{\epsilon} M \rar M_{\A^{\amalg}} \rar \Sigma e(M).$
Applying this right adjoint to the trivial module $K$ gives a module $e(K)$, which is denoted by $E_\A$.
The $KG$-module $E_\A$ associated to the thick ideal $\A$ in this way is an idempotent module (i.e., $E_\A \otimes E_\A = E_\A$). (These modules were introduced by
Rickard and they proved to be very useful objects in modular representation theory.) Krause then shows that this 
idempotent module is an endofinite object (see \cite[Definition 1.1]{Kr}) in
$\A^{\amalg}$ and hence admits a splitting into indecomposable modules: $E_\A = \bigoplus_{i \in I} E_i$ with $E_i \otimes E_j =0$ for 
$i \ne j$. 
Let $(\epsilon_i): \coprod_i E_i \rar K$ be the decomposition of $\epsilon: E_\A \rar K$ and define $\A_i$ to be the collection of all 
modules $X$ in 
$\stmod(KG)$ such that $\epsilon_i \otimes X$ is an isomorphism. It then follows that $\A = \coprod_i \A_i$ is the Krull-Schmidt 
decomposition for $\A$.
See \cite{Kr} for more details.

We draw the following corollary by applying theorem \ref{th:main} to the above decomposition.

\begin{prop} \label{prop:slick} Let $\A$ be a thick ideal of $\stmod(KG)$ and let $\amalg \A_i$ be the Krull-Schmidt decomposition of $\A$. Then,
\[ K_0(A) \cong \underset{i \in I}{\bigoplus} \;K_0(\A_i).\]
\end{prop}
\begin{proof} We show that condition (1) of our main theorem is satisfied. This is done in \cite[Lemma 2.4]{Kr} but we include it here 
for the reader's convenience. For any two thick ideals  $\A_1$ and $\A_2$ that occur in the decomposition for $\A$, we want to show that 
$\Hom(\A_1, \A_2)=0$. Towards this, consider objects $X \in \A_1$ and $Y \in \A_2$ and note that every map $X \rar Y$ in $\stmod(KG)$ 
factors through $X \otimes \Hom_K(X,Y)$ in the obvious way. So we will be done if we can show that $\Hom_K(X,Y)$ is a projective $KG$-module 
or equivalently a trivial module in  $\stmod(KG)$. To see this, first note that 
$\Hom_K(X,Y) \cong X^* \otimes Y$, where $X^*$ denotes the $K$-dual of $X$. But $X^*$ is a retract of $X^* \otimes X \otimes X^*$
(see \cite[Lemma A.2.6]{mps}) and therefore belongs to $\A_1$. This implies that $X^* \otimes Y \in \A_1 \otimes \A_2$. 
The latter is zero by lemma \ref{le:short} and therefore $\Hom_K(X,Y) = X^* \otimes Y$ is zero in $\stmod(KG)$. 
\end{proof}

\section{Derived categories of rings} \label{sec:KS-derivedcategory}

In this section we will prove a Krull-Schmidt theorem for the thick subcategories of perfect complexes over a noetherian ring.
We work in some level of formality here so that we can generalise our results easily to noetherian stable homotopy categories in the sense
of Hovey-Palmieri-Strickland \cite{mps}. We begin with some lemmas.

\begin{lemma} \label{le:support} For $X$ in $D^b(\proj\,R)$, let $DX$ denote its Spanier-Whitehead dual, i.e., 
$DX = \RHom(X,R)$. Then $\Supp(X) = \Supp(DX).$ 
\end{lemma}
\begin{proof} $X$ is a retract of $X \otimes DX \otimes X$ \cite[Lemma A.2.6]{mps}, and therefore $DX$ is a retract of $DX \otimes X \otimes DX$. 
Now it is clear that both $X$ and $DX$ have the same support.
\end{proof}

Recall that a subset $S$ of $\Spec(R)$ is a \emph{thick support} if it is a union of (Zariski) closed sets $S_{\alpha}$ such that $\Spec(R) - S_{\alpha}$ 
is a quasi-compact set. Also note that if $R$ is noetherian, a thick support is just a specialisation-closed subset.

Following \cite{mps}, we will denote $\Hom_{D^b(\proj\,R)}(\Sigma^* X ,Y)$ by $[X,Y]_*$ and the internal function spectrum 
$\RHom(X,Y)$ by $F(X,Y)$. With these notations, we have the following natural isomorphism \cite{mps}
\[[X \otimes A , B]_* \cong [X , F(A,B)]_* .\] 
This isomorphism gives the following useful lemma.

\begin{lemma} \label{le:cond1} If $X$ and $Y$ are perfect complexes such  that $\Supp(X) \cap \Supp(Y) = \emptyset$, 
then $[X,Y]_* = 0.$ In particular if $A$ and $B$ are any two disjoint thick supports of $\Spec(R)$, then  $[\T_A, \T_B]_* = 0$.
\end{lemma}

\begin{proof} By Spanier-Whitehead duality, $[X,Y]_* = [R \otimes X, Y_*] = [R, F(X,Y)]_*$. Therefore $[X,Y]_* = 0 $ if and only if
$F(X,Y) = 0$. But since $X$ is a small object in $D(R)$, $F(X,Y) = DX \otimes Y$ \cite[Appendix A.2]{mps}. So we have to show that 
$DX \otimes Y = 0$. Since $\Supp(DX) = \Supp(X)$  is disjoint with $\Supp(Y)$, given any prime $p$, either 
$p$ is not in $\Supp(DX)$ or it is not in $\Supp(Y)$. In the former case, $R_p \otimes DX = 0$, and in the latter, $R_p \otimes Y =0$.
In either case, we get $DX \otimes Y \otimes R_{p} = 0$. Since $p$ is an arbitrary prime, we get $DX \otimes Y = 0$. 
\end{proof}

By a \emph{thick decomposition} of a thick support $S$, we will mean a decomposition $S = \bigcup S_i$ into non-empty thick supports, 
where $S_i \cap S_j = \emptyset$ if $i \ne j$. A thick decomposition is \emph{Krull-Schmidt} if the $S_i$ 
do not admit nontrivial thick decompositions. Here are a few examples.

\begin{example} All rings considered are commutative.
\begin{enumerate}
\item{If $R$ is a PID, then every non-zero prime ideal is a maximal ideal and so the thick supports are:
$\Spec(R)$, and all subsets of the set of maximal ideals. In particular, $\Spec(R)$ is an indecomposable thick support.}
\item{If $R$ is an Artinian ring, it has only finitely many prime ideals and every prime ideal is also a maximal ideal. So every 
non-empty subset $S$ of $\Spec(R)$ is a thick support. It is then clear that $\bigcup_{p \in S} \{ p \}$ is the Krull-Schmidt
decomposition for $S$.}
\end{enumerate}
\end{example}

\n
\emph{Note:} The components of the Krull-Schmidt decomposition of a thick support $S$ are not (in general) the connected 
components of $S$. For example,  $\Spec(\ints) - \{ 0 \}$ is a connected subset of $\Spec(\ints)$; however, it is not
indecomposable as a thick support. In fact, 
\[ \Spec(\ints) - \{ 0 \} = \underset{p \ \text{prime}}{\bigcup}\; \{ (p) \}\]
is its Krull-Schmidt decomposition.

Now we establish the strong connection between the decompositions of thick supports and the decompositions of thick subcategories.

\begin{prop} \label{prop:1} Let $C = A \cup B$ be a thick decomposition of a thick support. Then this induces a decomposition of the associated thick
subcategories:
\[\T_C = \T_A \amalg \T_B.\]
\end{prop}
\begin{proof} Clearly $\T_A \cap \T_B = 0$ because $A$ and $B$ are disjoint by definition. We have to show that the objects of $\T_C$
are coproducts of objects in $\T_A$ and $\T_B$. We give an indirect proof of this statement using the thick subcategory theorem.
Lemma \ref{le:cond1} and lemma \ref{le:unclear} tell us that $\T_A \amalg \T_B$ is a
thick subcategory. So we will be done if we can show that the thick support corresponding to  $\T_A \amalg \T_B$ is $C$.
The thick support corresponding to $\T_A \amalg \T_B$  is  given by $\underset{x \in \T_A \amalg \T_B}{\bigcup} \Supp(x)$.
This clearly contains both $A$ and $B$ and hence their union ($ = C) $.  To see the other inclusion,
just note that $\Supp(a \oplus b) = \Supp(a) \cup \Supp(b) \subseteq A \cup B = C.$
\end{proof}

Now we show that decompositions of thick subcategories gives rise to thick decompositions of the corresponding thick supports.

\begin{prop} \label{prop:2} Let $\C = \A \, \amalg \,\B$ be a decomposition of a thick subcategory $\C$ of $D^b(\proj\,R)$ and let $A$, $B$, and $C$ denote the
corresponding thick supports of these thick subcategories. Then $C = A \cup B$ is a thick decomposition of $C$.
\end{prop} 
\begin{proof} We first observe that the intersection of thick supports is again a thick support: Let $S$ and $T$ be two thick supports of $\Spec(R)$.
Then $S = \bigcup_{\alpha} V(I_{\alpha})$ and $T = \bigcup_{\beta} V(J_{\beta})$ for some finitely  generated ideals $I_{\alpha}$ and $J_{\beta}$
in $R$. Now,
 \[ S \cap T = \bigcup_{\alpha , \beta} \left(V(I_{\alpha}) \cap V(J_{\beta})\right) = \bigcup_{\alpha , \beta} V(I_{\alpha} + J_{\beta}).\] 
Since the sum of two finitely generated ideals is  finitely generated, we conclude that  $S \cap T$ is a thick support.

We now argue that $A$ and $B$ are disjoint. If $A$ and $B$ are not disjoint, then their intersection being a thick support corresponds to a non-zero thick subcategory 
that is contained in both $\A$ and $\B$. This contradicts the fact that $\A \cap \B = 0$, so $A$ and $B$ have to be disjoint. 
It remains to show that $A \cup B$ is $C$.  For this, note that every complex $c \in \C$ splits as $c = a \amalg b$ and recall that 
$\Supp(c) = \Supp(a) \cup \Supp(b)$. It follows that $C = A \cup B.$ 
\end{proof}

\begin{rem} It can be easily seen that the last proposition, together with lemma \ref{le:cond1}, implies that the cancellation property holds for thick 
subcategories of $D^b(\proj\,R)$.
\end{rem}

Combining proposition \ref{prop:1} and \ref{prop:2} we get the following decomposition result.

\begin{thm}\label{thm:aux} Let $R$ be a commutative ring and let $\T_S$ be the thick subcategory of $D^b(\proj \,R)$ corresponding to a thick support $S$.
Then $\T_S$ admits a Krull-Schmidt decomposition if and only if $S$ admits one.
\end{thm}

Applying our main theorem \ref{th:main} to the above decomposition gives,

\begin{cor} \label{cor:Ktheory} Let  $\T_S  = \coprod_{i \in I} \T_{S_i}$ be a decomposition corresponding to a thick decomposition
$S = \bigcup_{i \in I} S_i$. Then,
\[ K_0(\T_S) \cong \underset{i \in I}{\bigoplus} \; K_0(\T_{S_i}).\]
\end{cor}
\begin{proof} Condition (1) of our main theorem holds here by lemma \ref{le:cond1}.
\end{proof}

The question that remains to be addressed is the following. When do thick supports of $\Spec(R)$ admit 
Krull-Schmidt decompositions? We show that this is always possible if $R$ is noetherian. 

\begin{prop} \label{prop:graph} Let $R$ be a noetherian ring and let $S$ be a thick support of $\Spec(R)$. 
Then there exists a unique Krull-Schmidt decomposition $\bigcup S_i$ for $S$.
\end{prop}

\begin{proof} 
It is well-known that the set of prime ideals in a noetherian ring satisfies the descending
chain condition \cite[Corollary 11.12]{am}. (This condition is equivalent to saying that any non-empty collection of primes ideals has a minimal
element.) To start, let $S$ be a thick support in $\Spec(R)$ and let $(p_i)_{i \in I}$ be the collection of all minimal elements in $S$ -- 
i.e., primes $p \in S$ which do not contain any other prime in $S$.  It is now clear (using the above fact about noetherian rings)
that every prime $p \in S$ contains a minimal element $p_i \in S$, therefore $S = \bigcup_{i \in I} V(p_i)$. (Also note that each 
$V(p_i)$ is a closed subset of $\Spec(R)$ and hence a thick support.)  Now define the \emph{thick graph} $G_S$ of $S$ as follows:
The  vertices are the minimal primes $(p_i)_{i \in I}$ in $S$, and two vertices $p_i$ and $p_j$ are adjacent if and only if 
$V(p_i) \cap V(p_j) \ne \emptyset$. 
Let $(C_k)_{k \in K}$ be the connected components of this graph and for each $C_k$ define a thick support
\[S_k:= \underset{p_i \in C_k}{\bigcup}{ V(p_i)}.\]
By construction it is clear that $\bigcup S_k$ is a thick decomposition of $S$. It is not hard to see that each $S_k$ is indecomposable.
This is done by showing that any thick decomposition of $S_k$ disconnects the connected component $C_k$ of the thick graph of $S$.
Finally the uniqueness part: let $\bigcup T_k$ be another Krull-Schmidt decomposition of $S$. 
It can be easily verified that the minimal primes in $T_k$ are precisely the minimal primes of $S$ that are contained in $T_k$. Thus the 
Krull-Schmidt decomposition $\bigcup T_k$ gives a partition of the set of minimal primes in $S$. This partition induces a decomposition of the
thick graph of $S$ into its connected (since each $T_k$ is indecomposable) components. Since the decomposition of a graph into its connected 
components is unique, the uniqueness of Krull-Schmidt decomposition follows.
\end{proof}

\begin{rem} The careful reader will perhaps note that the proof of proposition \ref{prop:graph} makes use of \emph{only} the following two 
conditions on $R$.
\begin{enumerate}
\item{Every open subset of $\Spec(R)$ is compact.}
\item{$\Spec(R)$ satisfies the descending chain condition.}
\end{enumerate}
Therefore the proof generalises to any $R$ which satisfies these two properties. As mentioned above, it is well-known  
that noetherian rings satisfy these two properties. Moreover, any ring which has finitely
many prime ideals automatically satisfies these properties. Note that there are non-noetherian rings which have finitely many primes. For example,
the ring \[ R = \mathbb{C}[x_2,x_3,x_4, \cdots]/(x_2^2,x_3^3,x_4^4,\cdots)\] 
is a non-noetherian ring with only one prime ideal. So everything we are going state for the remainder of this section will work for rings
which have these two properties, but we state our results only for noetherian rings for obvious cognitive reasons.
\end{rem}

After all that work, the following theorem is now obvious.

\begin{thm} \label{th:mainalg} If $R$ is a noetherian ring, then every thick subcategory of $D^b(\proj\,R)$ admits a unique Krull-Schmidt decomposition.
Conversely, given any collection of thick subcategories $(\C_i)_{i \in I}$ in $D^b(\proj\,R)$
such that $\C_i \cap \C_j = 0$ for all $i \ne j$, there exists a thick subcategory $\C$ such that $\C= \coprod_{i \in I} \C_i$.
\end{thm}
\begin{proof} The first part follows from theorem \ref{thm:aux} and proposition \ref{prop:graph}. For the second part, define $\C$ be the
full subcategory of all finite coproducts of objects from the $\C_i$. Thickness of $\C$ follows by combining  theorem \ref{thm:aux}, lemma \ref{le:cond1}, 
and lemma \ref{le:unclear}.
\end{proof}

\begin{cor} A noetherian ring $R$ is local if and only if every thick subcategory of $D^b(\proj\,R)$ is indecomposable.
\end{cor}
\begin{proof} If $R$ is local, then it is clear that every thick support of $\Spec(R)$ contains the unique maximal ideal. In particular,
we cannot have two non-empty disjoint thick supports and therefore every thick subcategory of $D^b(\proj\,R)$ is indecomposable.
Conversely, if $R$ is not local, it is clear that the thick subcategory supported on
the set of maximal ideals is decomposable.
\end{proof}

Theorem \ref{th:mainalg} also gives the following interesting algebraic results.

\begin{cor} \label{cor:purealg}  Let $X$ be a perfect complex over a noetherian ring $R$. Then $X$ admits a unique splitting into perfect complexes,
\[ X \cong \underset{i \in I}{\bigoplus} X_i\]
such that the supports of the $X_i$ are pairwise disjoint and indecomposable.
\end{cor}

\begin{proof} Let $\bigcup S_i$ be the Krull-Schmidt decomposition of $\Supp(X)$ (which exists by proposition \ref{prop:graph}). 
It is clear from  theorem \ref{th:mainalg} that $X$ admits a splitting ($X \cong \oplus X_i$) where $\Supp(X_i) \subseteq S_i$. To see that we have equality  
in this inclusion, observe that 
\[ \Supp (X) = \Supp \left( \underset{i \in I}{\bigoplus} X_i \right) = \underset{i \in I}{\bigcup} \Supp(X_i) 
\subseteq \underset{i \in I}{\bigcup} S_i = \Supp(X) .\]
Uniqueness: If $X$ admits another decomposition $\oplus Y_i$ as above, then by proposition \ref{prop:graph} we know that both $\bigcup\, \Supp(X_i)$ 
and $\bigcup \, \Supp(Y_i)$ are the same decompositions of $\Supp(X)$.  Lemma  \ref{le:cond1} then implies that $X_i \cong Y_i$ for all $i$. So the decomposition is unique.
\end{proof}

Srikanth Iyengar pointed out  that this splitting holds in a much more generality. In fact it holds for all complexes in the derived category which have bounded and finite 
homology (i.e. complexes whose homology groups $H_i(-)$ are finitely generated and are zero for all but finitely many $i$.) The full subcategory of such complexes will be 
denoted by $D^f(R)$. 

\begin{cor} Let $R$ be a noetherian ring. Every complex $X$ in $D^f(R)$ admits a splitting 
\[ X \cong \underset{i \in I}{\bigoplus} X_i\]
such that the supports of the $X_i$ are pairwise disjoint and indecomposable.
\end{cor}

\begin{proof} Since $X$ has bounded and finite homology, $\Supp(X)$ is a closed set and therefore proposition \ref{prop:graph} applies and 
we get a decomposition 
\[ \Supp(X) =  \underset{i \in I}{\bigcup} S_i.\]
We now use a result of Neeman \cite{Ne} which states that the lattice of localising subcategories (thick subcategories that are closed under
arbitrary coproducts) is isomorphic to the lattice of all subsets of $\Spec(R)$.
If $\mathcal{L}$ denotes the localising subcategory generated by $X$, then we get a Krull-Schmidt decomposition 
$ \mathcal{L} \cong \coprod \mathcal{L}_i$ that corresponds to the above decomposition of the support of $X$. (This follows exactly as 
in 
the case of thick subcategories.)  The splitting of $X$ as stated in the corollary is now clear. 
\end{proof}

\begin{example} We illustrate the last two corollaries with some examples.
\begin{enumerate}
\item{For every integer $n > 1$, let $M(n)$ denote the Moore complex $0 \rar \mathbb{Z} \stk{n} \mathbb{Z} \rar 0$
and let $p_1^{i_1}p_2^{i_2}\cdots p_k^{i_k}$ be the unique prime factorisation of $n$. It is easy to see that 
$\bigcup_{t=1}^k \{ (p_t)\} $ is the Krull-Schmidt decomposition of $\Supp (M(n))$.
Then it follows that
\[M(n) \cong \bigoplus_{t=1}^k\, M(p_t^{i_t})\] 
is the splitting of $M(n)$  corresponding to this Krull-Schmidt decomposition.}
\item{ Let $R$ be a self-injective noetherian ring (hence also an Artinian ring). Then it follows (because every prime ideal is also maximal) 
that $\bigcup_{p \in \Spec(R)} \{ p\}$ 
is the Krull-Schmidt decomposition  for $\Supp(R)$.  It can be shown that 
\[ R \cong \bigoplus_{p \in \Spec(R)} E(R/p)\] is the splitting of $R$ corresponding to this Krull-Schmidt decomposition.
(Here $E(R/p)$ denotes the injective hull of $R/p$.)}
\end{enumerate}
\end{example}

\section{Noetherian stable homotopy categories.} \label{sec:nsht}

Motivated by the work in the previous section, we now state a Krull-Schmidt theorem for noetherian stable homotopy categories.
We use the language and results of \cite{mps} freely. We explain how the proofs of the previous section generalise to give us a
 Krull-Schmidt theorem in this more general setting by invoking the appropriate results from \cite{mps}.
We begin with some definitions and preliminaries from Axiomatic Stable Homotopy Theory \cite{mps}.

\begin{defn} \cite{hp} A \emph{unital algebraic stable homotopy category} is a tensor triangulated category $\mathscr{C}$ with the following 
properties.
\begin{enumerate}
\item{Arbitrary products and coproducts of objects in $\mathscr{C}$ exist.}
\item{$\mathscr{C}$ has a finite set $\mathscr{G}$ of weak generators: i.e., $X \cong 0$ if and only if $[A,X]_*=0$ for all $A \in \mathscr{G}$.}
\item{The unit object $S$ and the objects of $\mathscr{G}$ are small.}
\end{enumerate}
$\mathscr{C}$ is a \emph{noetherian stable homotopy category} if in addition the following conditions are satisfied,
\begin{itemize}
\item{ $\pi_*(S) := [S,S]_{*}$ is commutative and noetherian as a bigraded ring. }
\item {For small objects $Y$ and $Z$ of $\C$, $[Y,Z]_{*}$ is a finitely generated module over $[S,S]_{*}$. }
\end{itemize}

 The stable homotopy category of spectra is unital and algebraic but not noetherian. The derived
category  $D(R)$ of a commutative ring $R$ is a noetherian stable homotopy category if and only if $R$ is noetherian
(since $[R,R]_*=R$). $\Spec(\pi_*S)$ will stand for the collection of all homogeneous prime ideals of $\pi_*S$ with the
Zariski topology. 
\end{defn}

Henceforth $\mathscr{C}$  will denote a  bigraded noetherian stable homotopy category. For each  thick ideal $\mathcal{A}$
of finite objects in $\mathscr{C}$, there is a \emph{finite localisation functor} $L_{\A}$ (also denoted $L_{\A}^f$) on $\mathscr{C}$ whose finite acyclics 
are precisely
the objects of $\A$; see \cite[Theorem 2.3]{hp}. For each bihomogeneous prime ideal $\mathbf{p}$ in $\pi_*(S)$, there are finite localisation
functors $L_{\mathbf{p}}$ and $L_{< \mathbf{p}}$ whose finite acyclics are 
$\{X \mbox{\;finite\;}\;|\; X_{\mathbf{q}} = 0 \;\forall \; \mathbf{q} \subseteq \mathbf{p} \}$ and 
$\{X \mbox{\;finite\;}\;| \; X_{\mathbf{q}} = 0 \;\forall \; \mathbf{q} \subsetneq \mathbf{p}  \}$ respectively. 
These functors define, for each $X \in \mathscr{C}$, a natural exact triangle
\[M_{\mathbf{p}} X \rar L_{\mathbf{p}} X \rar L_{<\mathbf{p}} X \rar \Sigma M_{\mathbf{p}} X.\]
We set $M_{\mathbf{p}} : = M_{\mathbf{p}} S$. See \cite{hp} for more details.

We say that $\mathscr{C}$ has the \emph{tensor product property} if for all bihomogeneous prime ideals $\mathbf{p}$ in $\pi_*(S)$ and all objects 
$X, Y \in \mathscr{C}$,  $M_{\mathbf{p}} \wedge X \wedge Y = 0$ if and only if either $M_{\mathbf{p}} \wedge X = 0$ or  $M_{\mathbf{p}} \wedge Y =0$.

%

We now mimic the setup for the derived category of a ring. For $X$ in $\mathscr{C}$, define its \emph{cohomological support variety} as
\[\Supp(X) = \{\mathbf{p}\;|\; M_{\mathbf{p}} \wedge X \ne 0 \} \subseteq \Spec (\pi_*S),\]
and if $\D$ is any thick subcategory of $\mathscr{C}$, define 
\[\Supp(\D) = \underset{X \in \D}{\bigcup} \Supp(X).\]  
It is a fact \cite[Theorem 6.1.7]{mps} that for $X$ a small object in $\mathscr{C}$,  $\Supp(X)$ is a Zariski closed subset and hence also a thick support of
$\Spec \pi_*(S)$.

We are now ready to state the thick subcategory theorem for $\mathscr{C}$. 

\begin{thm} \label{th:tshps} \cite{hp} Let $\mathscr{C}$ be a noetherian stable homotopy category satisfying the tensor product property.
Then the lattice of thick ideals of small objects in $\mathscr{C}$ is isomorphic to the the lattice of thick
supports of $\Spec(\pi_*S)$. Under this isomorphism, a thick ideal $\A$ corresponds to $\Supp(\A)$, and
a thick support $T$ of $\Spec(\pi_*S)$ corresponds to  the thick ideal $\{ Z \in \mathscr{C}\;| \;\Supp(Z) \subseteq T\}.$
\end{thm}

Now using this thick subcategory theorem, we can run our proofs from the previous section which generalise to prove theorem \ref{th:grand} below. 
We will also mention the minor changes that are necessary for this generalisation. We begin with an algebraic lemma.

\begin{lemma} \label{le:benson} Let $R$ be a graded commutative noetherian $K$-algebra. Then $R$ satisfies the descending chain condition on the set 
$\Spec(R)$ of homogeneous prime ideals of $R$.
\end{lemma}

\begin{proof} First note that if char($K$) = 2, then $R$ is a strictly commutative graded ring and the result for such $R$ is well-known; see \cite[Corollary 11.12]{am}.
So assume that $\text{char} (K) \ne 2$. Since every homogeneous prime ideal contains the nilradical, the ring map
\[R \rar R/{\text{nilrad}\, R} \]
induces an order-preserving bijection between $\Spec(R)$ and $\Spec(R/{\text{nilrad}\, R})$. So in view of this remark, we can assume without any
loss of generality that $\text{nilrad} \,R = 0$. But if  $\text{nilrad}\, R = 0$, then it is concentrated in even degrees: If $x$ is an odd 
degree element,  then $x^2 = -x^2$.  This implies $x^2 = 0$ (since char$(K) \ne 2$), or equivalently $x \in \text{nilrad}\, R = 0$.  
This shows that $R$ is concentrated 
in even degrees. In particular it is a strictly commutative graded ring, so we are done.
\end{proof}

\begin{thm} \label{th:grand} Let $\mathscr{C}$ be a noetherian stable homotopy category satisfying the tensor product property.
Then every thick ideal $\A$ of small objects in $\mathscr{C}$ admits a  Krull-Schmidt decomposition: 
$\A = \coprod_{i \in I} \A_i$. Consequently, 
\[ K_0(\A) \cong \underset{i \in I}{\bigoplus}\;\; K_0(\A_i).\] 
Conversely, if $\A_i$ are thick ideals  of $\mathscr{C}$ such that $\A_i \cap \A_j = 0$ for all $i \ne j$, then there exists a thick
ideal $\A$ such that $\A = \coprod_{i \in I} \A_i$.  
\end{thm}

\begin{proof} We will explain why all the lemmas and propositions leading up to theorem \ref{th:mainalg} hold in this generality. 
We proceed in order starting from lemma \ref{le:support}.

Lemma \ref{le:support}: This holds for any strongly dualisable object in a closed symmetric monoidal category \cite[Lemma A.2.6]{mps}. 
By \cite[theorem 2.1.3 (d)]{mps}, every small object is strongly dualisable in a unital algebraic stable homotopy category and therefore lemma \ref{le:support} 
holds for small objects of $\mathscr{C}$.

Lemma \ref{le:cond1}:  Proceeding in the same way, it boils down to showing that $DX \wedge Y$ is zero. Since  $DX$ and $Y$ have disjoint 
supports, for each prime $\mathbf{p}$ in $\Spec(\pi_* S)$, either $Y \wedge M_{\mathbf{p}} = 0$ or 
$DX \wedge M_{\mathbf{p}} = 0$. This implies that $DX \wedge Y \wedge M_{\mathbf{p}} = 0$ for all $\mathbf{p}$. By theorem  \cite[Theorem 6.1.9]{mps}, 
we have the following equality of Bousfield classes:
\[ \la S \ra = \underset{\mathbf{p} \in \Spec(\pi_*S)}{\coprod} \la M_{\mathbf{p}}\ra . \]
This implies $(DX \wedge Y) \wedge S = DX \wedge Y = 0$.

Proposition \ref{prop:1}: This goes through verbatim without any changes.

Proposition \ref{prop:2}: This is actually much easier. Thick supports of $\Spec(\pi_*S)$ are just unions of Zariski-closed sets
and therefore it is obvious that intersection of thick supports is thick. The rest of the proof of this proposition follows exactly the same
way. 

Using these lemmas and proposition, theorem \ref{thm:aux} and corollary \ref{cor:Ktheory} follow. Proposition \ref{prop:graph} also holds in the graded 
noetherian case; see lemma \ref{le:benson}. Thus the first part of theorem follows by combining all these lemmas and propositions.

The converse follows as before by combining theorem \ref{thm:aux}, lemma \ref{le:cond1}, and lemma \ref{le:unclear}.
\end{proof}

\begin{cor}  Let $\mathscr{C}$ be a noetherian stable homotopy category satisfying the tensor product property. Then every object $X$ in $\mathscr{C}$ admits a
unique  decomposition
\[ X \cong \underset{i \in I}{\coprod} X_i,\]
such that the support varieties of the $X_i$ are pairwise disjoint and indecomposable.
\end{cor}

\begin{rem} The derived category $D(R)$ of a noetherian ring $R$ is a noetherian stable homotopy category which satisfies the tensor product property.
(Note that $\pi_*(S)= [R,R]_* = R$.)
So theorem \ref{th:tshps} recovers our Krull-Schmidt theorem for derived categories of noetherian rings.
\end{rem}

We now give another example of a noetherian stable homotopy category for 
which theorem \ref{th:grand} applies.

\begin{example} \cite{hp, hp2} 
Let $B$ denote a finite dimensional graded co-commutative Hopf algebra over a field $k$ ($\text{char}(k) > 0$).  Then the category $K(\Proj\,B)$, whose objects are 
unbounded chain complexes of projective $B$-modules and whose morphisms are chain homotopy classes of chain maps, is a noetherian stable homotopy category.
In this category the unit object is a projective resolution of $k$ by $B$-modules and therefore $\pi_*(S)$ is isomorphic
to $\Ext_B^{*,*}(k,k)$. It is well-known that this ext algebra is a bigraded noetherian algebra; see \cite{fs} and \cite{wil}.
In particular, if $B$  is a finite dimensional sub-Hopf algebra of the $\mbox{mod}\; 2$ Steenrod algebra, it is a deep result 
\cite[Corollary 8.6]{hp} that the category $K(\Proj\,B)$ satisfies the tensor product property. So theorem \ref{th:grand} applies and we get the following 
result.
\end{example}

\begin{thm} \label{thm:KSforstable(B)} Let $B$ denote a finite dimensional graded co-commutative Hopf algebra satisfying the tensor product property.
 Then every thick subcategory of small objects in $K(\Proj\,B)$ is indecomposable.
\end{thm}
 
\begin{proof} 
 The augmentation ideal of $\Ext_B^{*,*}(k,k)$ which consists of all elements in positive degrees is the unique maximal homogeneous ideal and therefore is contained in 
every thick support of $\Spec \Ext_B^{*,*}(k,k)$. Now by theorem \ref{th:grand} it follows that 
every thick subcategory of small objects in this category is indecomposable.
\end{proof}

\section{Stable module categories over Hopf algebras} \label{sec:KS-hopfalgebras}

In this section we give a generalised version of Krause's result \cite[Theorem A]{Kr}. 
To this end, let $B$ denote a finite dimensional graded co-commutative Hopf algebra over a field
$k$. Following \cite{hp}, we say that $B$ satisfies the tensor product if the category $K(\Proj\,B)$ does. We now prove a Krull-Schmidt theorem 
for $\stmod(B)$, whenever $B$ satisfies the tensor product property.  The category $\StMod(B)$ fails to be a noetherian stable homotopy category
in general. (In fact, $\pi_*(S)= \Hom_{\stmod(B)}(k,k)$ is isomorphic to the Tate cohomology of $B$ which is not noetherian in general.) So we take the
standard route of using a finite localisation functor \cite[Section 9]{mps} to get around this problem.

We now review some preliminaries from \cite[Section 9]{mps}. There is a finite localisation functor $L_B^f: K(\Proj\,B) \rar K(\Proj\,B)$, whose finite 
acyclics are precisely the objects in the thick subcategory generated by $B$. It has been shown \cite{mps} that the category $\StMod(B)$ is equivalent 
to the full subcategory of $L_B^f$-local objects. This finite localisation functor also establishes a bijection between the nonzero thick ideals of 
small objects in $K(\Proj\,B)$ and the thick ideals of $\stmod(B)$. In fact, the thick ideals of $\stmod(B)$ are in poset bijection with the non-empty specialisation 
closed subsets of $\Spec \, \Ext_{B}^{*,*}(k,k)$. (Note that for the category $K(\Proj\,B)$, the empty set of $\Spec \, \Ext_{B}^{*,*}(k,k)$ corresponds to the thick 
subcategory consisting of the zero chain complex.) Finally let $\mathfrak{m}$ denote the unique bihomogeneous maximal ideal of $\Ext_B^{*,*}(k,k)$ which corresponds to the
trivial thick subcategory in $\stmod(B)$. Then here is the main result of this section.

\begin{thm} \label{thm:mymainthm} Let $B$ denote a finite dimensional graded co-commutative Hopf algebra satisfying the tensor product property. Then every thick ideal 
of $\stmod(B)$ admits a unique Krull-Schmidt decomposition. Conversely, given thick ideals $\T_i$ in $\stmod(B)$ satisfying $\T_i \cap \T_j = 0$ for $i \ne j$,
there exists a thick ideal $\T$ such that $\T \cong \coprod \T_i$. Moreover, 
\[K_0(\T) \cong \underset{i \in I}{\bigoplus}\; K_0(\T_i).\]
\end{thm}

\begin{proof} Recall once again that the thick ideals of the stable module category are in poset bijection with the non-empty subsets of  
$\Spec \; \Ext_{B}^{*,*}(k,k)$ that are closed under specialisation. To begin, let $\T$ be a thick ideal of $\stmod(B)$ and let $S$ denote the 
corresponding specialisation closed subset. Since $\Ext_B^{*,*}(k,k)$ is a graded noetherian algebra, it follows that 
$S$ admits a unique decomposition $S = \bigcup_i S_i$ such that $S_i \cap S_j = \mathfrak{m}$ 
 for all $i \ne j$ and $S_i$ indecomposable (the proof of
Proposition \ref{prop:graph} can be adapted to the set $S-\mathfrak{m}$ by making use of lemma \ref{le:benson}). 
We now claim that $\T = \coprod \T_{S_i}$, where $\T_{S_i}$ denotes the thick ideal of $\stmod(B)$ corresponding to $S_i$. 

The first thing to note is that for $i \ne j$, $\T_{S_i} \cap \T_{S_j} = 0$:
Suppose not, then the non-zero thick ideal $\T_{S_i} \cap \T_{S_j}$ corresponds to a non-empty subset of $S_i \cap S_j-\mathfrak{m}$ which is impossible.
Therefore it follows that  $\T_{S_i} \cap \T_{S_j} = 0$. 

The next step is to show that the objects in $\T$ are the finite coproducts of the objects in 
$\T_{S_i}$. We give an indirect proof (as before) by showing that $\coprod \T_{S_i}$ is a thick ideal. Then the theorem follows because 
the
specialisation closed subset corresponding to $\coprod \T_{S_i}$ is clearly $\bigcup S_i = S$. To show that $\coprod \T_{S_i}$
is thick (ideal condition holds trivially), by lemma \ref{le:unclear}, all we need to show is that $\Hom_{\stmod}(A,B)= 0$, whenever $A \in \T_{S_i}$
and $B \in \T_{S_j}$ ($i \ne j$).  This was already done for group algebras in proposition \ref{prop:slick}. The proof given there can be easily 
generalised to finite dimensional co-commutative Hopf algebras. So that completes the proof of the existence of Krull-Schmidt decompositions.

For the converse, define $\T$ be the full subcategory of all finite coproducts of objects in $(\T_i)_{i \in I}$ and use the 
fact that $\Hom_{\stmod}(\T_i, \T_j) = 0$ for $i \ne j$, to complete the argument. 

The statement about $K_0$ groups follows from theorem  \ref{th:main}.
\end{proof}

\begin{rem} It is worth pointing out the difference between the decompositions in theorem \ref{thm:KSforstable(B)}
and those in theorem \ref{thm:mymainthm}: Every thick ideal of small objects in $K(\Proj \, B)$ is indecomposable; on the other hand, thick ideals of small 
objects in $\StMod(B)$ admit Krull-Schmidt decompositions. This striking difference is really the effect of
Bousfield localisation ($\StMod(B)$ is a Bousfield localisation of $\Stable(B)$).
\end{rem}

\subsection{Sub-Hopf algebras of the mod-$2$ Steenrod algebra}

It has been shown in \cite{hp} that every finite dimensional sub-Hopf algebra of $A$ (the mod $2$ Steenrod algebra) satisfies the tensor
product property. So theorem \ref{thm:mymainthm} applies and we get interesting Krull-Schmidt decompositions in the stable module categories over these
sub-Hopf algebras. We now look at some standard sub-Hopf algebras of $A$ and analyse decompositions in their stable module categories.
For every non-negative integer $n$, denote by $A(n)$ the finite dimensional sub-Hopf algebra of $A$ generated by $\{ \Sq^k |\, 1 \le k \le 2^n\}$.

\vskip 3mm \n
$A(0)$ - This is the sub-Hopf algebra of $A$ generated by $\Sq^1$ and therefore it is isomorphic to the exterior algebra $\ftwo[\Sq^1]/(\Sq^1)^2$ on $\Sq^1$. It can be easily 
shown that $\Ext_{A(0)}^{*,*} (\ftwo, \ftwo) \cong  \ftwo [h_0]$ with $|h_0|=(1,1)$. Therefore the bihomogeneous prime spectrum of $\Ext_{A(0)}^{*,*} (\ftwo, \ftwo)$
is $\{(0), (h_0)\}$. It is now clear from the thick 
subcategory theorem that there are only two  thick subcategories of $\stmod(A(0))$; the trivial category which corresponds to maximal prime $(h_0)$ and the entire 
category $\stmod(A(0))$ which corresponds to $\{(0), (h_0)\}$.  These subsets of prime ideals are clearly indecomposable closed subsets and therefore we conclude by 
theorem \ref{thm:mymainthm} that all thick subcategories of $\stmod(A(0))$ are indecomposable.
 
\vskip 3mm \n
$A(1)$ - This is the sub-Hopf algebra of $A$ generated by $\Sq^1$ and $\Sq^2$. This algebra can be shown (using the Adem relations) to be isomorphic to the free 
non-commutative  algebra (over $\ftwo$) generated by symbols $\Sq^1$ and $\Sq^2$ modulo the two sided ideal 
\[(\Sq^1 \Sq^1,\, \Sq^2 \Sq^2 + \Sq^1 \Sq^2 \Sq^1).\] 
The cohomology algebra of $A(1)$ has been computed in \cite{wil}:
\[\Ext_{A(1)}^{*,*}(\mathbb{F}_2,\mathbb{F}_2) \cong \mathbb{F}_2[h_0,h_1,w,v]/(h_0h_1, h_1^3,h_1w,w^2-h_0^2v),\]
where the bi-degrees are: $|h_0|=(1,1), |h_1|=(1,2), |w|=(3,7)$, and $|v|=(4,12)$. The bihomogeneous prime lattice of the above cohomology ring consists of the 
minimal prime $(h_1)$; primes $p_0 = (h_1, w, h_0)$ and $p_2 = (h_1, w, v)$ of height one; and the maximal prime $\mathfrak{m} = (h_1, w, h_0, v)$. 

\begin{figure}
\[
\xymatrix{
                           &  \mathfrak{m} \ar@{-}[dl] \ar@{-}[dr] &  \\
p_0  \ar@{-}[dr]  &                                                &  p_1 \ar@{-}[dl] \\
                           &  (h_1)                                &                                  
}
\]
\label{A(1)}
\caption{Bihomogeneous prime lattice of $\Ext_{A(1)}^{*,*} (\ftwo, \, \ftwo)$}
\end{figure}
The thick subcategory theorem now tells us that
there are five thick subcategories in $\stmod(A(1))$ which correspond to the thick supports $\{\mathfrak{m}\}$, $\{p_0, \mathfrak{m}\}$, 
$\{p_1, \mathfrak{m}\}$, $\{p_0, p_1, \mathfrak{m}\}$, and $\{(h_1), p_0, p_1, \mathfrak{m}\}$. It is clear that the only decomposable thick support in this 
list is $\{p_0, p_1, \mathfrak{m}\}$; in fact,
 \[\{p_0, p_1, \mathfrak{m} \} = \{ p_0, \mathfrak{m}\} \, \cup \,\{ p_1, \mathfrak{m}\}\] 
is the desired decomposition. This induces a Krull-Schmidt decomposition of thick subcategories supported on these thick supports:
 \[ \T_{\{p_0, p_1, \mathfrak{m}\}} \cong  \T_{\{p_0, \mathfrak{m}\}} \amalg \T_{\{p_1, \mathfrak{m}\}}.\] 
Note that the remaining four thick subcategories in $\stmod(A(1))$ are indecomposable.

\vskip 3mm 
Decompositions in stable module categories over $A(2)$ and higher are much more complicated.

\section{Finite spectra} \label{sec:KS-torsionspectra}

In this section we work in the triangulated category $\F$ of finite spectra. Recall that a  spectrum $X$ is 
torsion  if $\pi_*(X)$  is a torsion group. In particular, a finite spectrum $X$ is torsion ($p$-torsion) if $\pi_*(X)$ consists of finite abelian groups 
($p$-groups) in every degree. We will now show that the thick subcategories of  finite spectra admit Krull-Schmidt decompositions. 

$\FT$ will denote the category of finite torsion spectra, and 
$\C_{n,p}$ will denote the thick subcategory of finite $p$-torsion spectra consisting of all $K(n-1)_*$ acyclics

\begin{lemma} If $\A$ is a thick subcategory of $\F$, then either $\A = \F$ or $\A \subseteq \FT$.
\end{lemma}

\begin{proof} It suffices to show that if $X$ is any non-torsion finite spectrum, then the integral sphere spectrum $S$ belongs to the thick 
subcategory generated by $X$. Consider the cofibre sequence
\[ X \stk{p} X \rar X \wedge M(p) \rar \Sigma X. \]
It is easy to see that $X \wedge M(p)$ is a type-1 spectrum, therefore the thick subcategory generated by $X$ contains all finite torsion spectra.
The next step is to show that there is a cofibre sequence
\[ \coprod S^r \rar X \rar T \rar \coprod S^{r+1},\]
where $T$ is a finite torsion spectrum; see \cite[Proposition 8.9]{mar}. It is clear from these exact triangles that $S$ belongs to the thick subcategory generated by 
$X$.
\end{proof}

We now start analysing decompositions for thick subcategories of finite spectra.

\begin{prop} The category $\F$ is indecomposable.
\end{prop}
\begin{proof} If $\F = \A \amalg \B$ is any decomposition of $\F$, then that would giving a splitting $A \vee B$ of the sphere spectrum 
($S$).
Since $S$ is an indecomposable spectrum, either $A$ or $B$ must be the sphere spectrum, which implies that either $\A$ or $\B$ is equal to $\F$. 
\end{proof}

We now show that the thick subcategories of finite \emph{torsion} spectra admit non-trivial Krull-Schmidt decompositions.

\begin{thm}  Every thick subcategory $\D$ of $\FT$ admits a unique Krull-Schmidt decomposition.
Conversely, given any collection of thick subcategories $\D_i$ in $\FT$ such that $\D_i \cap \D_j = 0$ for all $i \ne j$, 
there exits a  thick subcategory $\D$ such that $\underset{i}{\coprod}\D_i$ is a Krull-Schmidt decomposition for $\D$.
Moreover $K_0(\D) \cong \bigoplus_i K_0(\D_i).$
\end{thm}

\begin{proof}  
For every prime number $p$, define $\D_p := \D \cap \C_{1,p}$. By \cite[8.2, Theorem 20]{mar}  every torsion spectrum $X$ can be written as   
$X = \bigvee_p X_{(p)}$ where $X_{(p)}$ is the $p$-localisation of $X$. If $X$ is a finite torsion spectrum, it is easy to see that $X_{(p)}$ 
is a finite $p$-torsion spectrum and therefore belongs to $\C_{1,p} (\subseteq \FT)$. Since $\D$ is a thick subcategory of $\FT$, all the 
summands $X_{(p)}$ also belong to $\D$. Therefore $X_{(p)}$ belongs to $ \C_{1,p} \cap \D =\D_p$.
Also since $X$ is a finite spectrum, the above wedge runs over only a finite set of primes. This implies that the spectra in 
$\D$ are the finite coproducts of spectra in $(\D_p)$.

$\D_p$ is a thick subcategory in $\FT$:  Showing the $\C_{1,p}$ is a thick subcategory of $\FT$ is a routine verification. Since intersection of
thick subcategories is thick,  $\D \cap \C_{1,p}$ is thick in $\FT$.

$\D_p \cap \D_q = 0:$   It is clear that if a finite spectrum is both $p$-torsion and $q$-torsion, for primes $p \ne q$, then  
it  is trivial -- i.e., $\C_{1,p} \cap \C_{1,q} = 0$. But $\D_p \cap \D_q \subseteq \C_{1,p} \cap \C_{1,q}$, therefore  $\D_p \cap \D_q = 0.$

$\D_p$ is indecomposable:  $\D_p$ is by definition a thick subcategory of $\C_{1,p}$.  The thick subcategories of $\C_{1,p}$ are well-known to be
nested \cite{hs}. In particular they are all indecomposable.

This completes the proof of the existence of a Krull-Schmidt decomposition for any thick subcategory $\D$ in $\FT$. 

Uniqueness: Let $\coprod \E_k$ be another Krull-Schmidt decomposition for $\D$. Then since each $\E_k$ is an indecomposable thick
subcategory of $\FT$, we conclude that $\E_k = \C_{n_k, p_k}$ for some integer $n_k$ and some prime $p_k$. Thus the decomposition is unique.

For the converse, define $\D$ be the full subcategory of all finite coproducts of objects from $\D_p$. We just have to show that $\D$
is a thick subcategory and then the rest follows.  Note that $[X,Y] = 0$ whenever $X$ is a $p$-torsion spectrum and $Y$ is a $q$-torsion spectrum \cite[8.2, Lemma 21]{mar}. 
Therefore thickness of $\D$ follows from lemma \ref{le:unclear}.
\end{proof}

The following corollary is clear from these decompositions.

\begin{cor} Let $\A$ be be a thick subcategory of $\F$. Then either $\A = \F$, or $\A = \coprod \C_{i,p}$, where $i$ ranges over  
some subset of positive integers and $p$ over some subset of primes.
\end{cor}

\section{Questions}

\subsection{Krull-Schmidt decompositions in $D^b(\proj\,R)$}

We have shown that the thick subcategories of perfect complexes over a \emph{noetherian} ring admit Krull-Schmidt
decompositions. Note that proposition \ref{prop:graph} was crucial for  this decomposition. Now the question that 
remains to be answered  is whether  this proposition holds for arbitrary commutative rings. In other words,
if $R$ is any commutative ring, is it true that every thick support of $\Spec(R)$ admits a unique Krull-Schmidt 
decomposition? An affirmative answer to this question would imply, by theorem \ref{thm:aux}, that thick
subcategories of perfect complexes over commutative rings (not necessarily noetherian) admit Krull-Schmidt decompositions.

\subsection{The tensor product property}

Recall that our decompositions for the thick subcategories of small objects in $K(\Proj\,B)$ and $\StMod(B)$ assumed that 
the finite dimensional co-commutative Hopf algebra $B$ satisfies the tensor product property.  Now the question
that arises is which finite dimensional Hopf algebras satisfy the tensor product property? This seems to be a question of independent interest; 
see \cite[Corollary 3.7]{hp}.
Friedlander and Pevtsova \cite{FriPev} have shown recently that the finite dimensional co-commutative 
\emph{ungraded} Hopf algebras satisfy the tensor product property. Using their ideas, or otherwise, we would like to know if 
this property also holds for finite dimensional graded co-commutative Hopf algebras. Also we do not know if the tensor product product 
property in \cite{FriPev} is equivalent to the one given in \cite{hp}.



\chapter{Wide subcategories of modules}

In this chapter we study analogues of thick subcategories for the category of
$R$-modules.  Following Hovey, we call them wide subcategories. 
After giving some definitions and examples we state Mark Hovey's classification
theorem for wide subcategories. Then we use this theorem to study  $K$-theory
and Krull-Schmidt decompositions for these categories in sections \ref{sec:ktheorywide} and
\ref{sec:KSwide} respectively. All rings will be commutative in this chapter.

\begin{defn} \cite{wide}  A full subcategory $\C$ of an abelian category
is said to be \emph{wide} if it is an abelian subcategory that is closed under extensions. 
\end{defn}

\begin{example} Easy exercises in algebra tell us that the following are wide subcategories.

 1. Rational vector spaces form a wide subcategory of abelian groups.

 2. For every prime $p$, $p$-local abelian groups form a wide subcategory of abelian groups.

 3. Finitely generated modules over a noetherian ring form a wide subcategory of $R$-modules.

 4. Finitely presented modules over a coherent ring form a wide subcategory of $R$-modules.


\end{example}

This definition is motivated by the study of thick subcategories. See \cite{wide} for more details.

In this chapter we want to  understand such wide subcategories of modules. However, in order to get a handle on such categories, we 
restrict ourselves to a smaller category called $\Wide(R)$:
The wide subcategory generated by the ring $R$ in the category of $R$-modules. In other words, $\Wide(R)$ is
the intersection of all the wide subcategories of $R$-modules that contain $R$.  So we now focus our attention on the
wide subcategories of $\Wide(R)$.  The lattice of these
subcategories will be denoted by $\mathcal{L}_{\Wide} (R)$. Analogously, the thick subcategory generated by $R$ in the derived category
will be denoted by $\Thick(R)$ (this is also equivalent to the chain homotopy category of perfect complexes) and the lattice of all the thick 
subcategories of $\Thick(R)$  will be denoted by $\mathcal{L}_{\Thick} (R)$.

It is hard to characterise $\Wide(R)$ for any given ring $R$. However, when 
$R$ is \emph{coherent} (every finitely generated ideal is finitely presented), 
it turns out the $\Wide(R)$ is precisely the category of finitely presented $R$-modules; see \cite[Lemma 1.6]{wide}. 
Similarly, it is an easy exercise to show that if $R$ is noetherian, then $\Wide(R)$ 
consists of the finitely generated $R$-modules.

\section{Classification of wide subcategories}

Recall that a ring is \emph{regular} if every finitely generated ideal has finite
projective dimension.  It is a standard fact \cite[Theorem 6.2.1]{sarah} that a coherent ring $R$ is regular if and only if every 
finitely presented $R$-module has finite projective dimension. Polynomial rings in infinitely many variables over a PID are examples of regular coherent rings that 
are not noetherian. We will denote the collection of thick
supports of $\Spec(R)$ by $\mathbb{S}$. The main result of \cite{wide} then states:

\begin{thm} \cite{wide}  Let $R$ be regular coherent ring. Then the lattice $\mathcal{L}_{\Wide}(R)$
is isomorphic to the lattice $\mathbb{S}$ of thick supports of $\Spec(R)$. Under this isomorphism, a wide subcategory
$\W$ corresponds to the thick support $\bigcup_{M\in \W} \Supp(M)$.
\end{thm}

In fact, there is a commutative diagram of poset bijections:
\begin{figure}[!h]
\[
\xymatrix{
\mathcal{L}_{\Wide}(R)\ar[rrrr]^{f}_{\cong}  &  & & &     \mathcal{L}_{\Thick}(R)  \\
& & \mathbb{S} \ar[ull]^{\xi}_{\cong} \ar[urr]_{\zeta}^{\cong} & & 
}
\]
\caption{Bijection between wide subcategories and thick subcategories}
\label{fig:wideandthick}
\end{figure}

\n
where $\xi(A) = \{ M \in \Wide(R): \Supp(M) \subseteq A \}$, $\zeta (A) = \{ X \in \Thick(R): \Supp(X)\subseteq A\}$, 
and $f(\W) = \{ X \in \Thick(R): H_n(X) \in \W \ \forall \, n \}$.

As indicated in the above diagram, Hovey arrives at this theorem by passing to the derived category and using a result of Thomason \cite{Th} which classifies
thick subcategories of perfect complexes over a commutative ring; see theorem \ref{thm:thomason-thicksubcategory}. 
In fact, it is in the process of this transition (from the module category to the derived category) that the regularity assumption enters into the picture.

\section{$K$-theory for wide subcategories} \label{sec:ktheorywide}

Since a wide subcategory is also an abelian category, it makes sense to talk
of its Grothendieck group. So we now try to relate the Grothendieck groups of wide subcategories and thick subcategories that
correspond to each other under the above bijection $f$. We begin with a simple motivating example.

\begin{example} Let $R = \mathbb{Z}$, or more generally a PID. Since $R$ is clearly regular and coherent, by the above theorem we know that 
the wide subcategories of $\Wide(R)$ and thick subcategories of $\Thick(R)$ are in bijection with the thick supports of 
$Spec(\mathbb{Z})$.  So if $S$ is any proper thick support of $\Spec(\ints)$, then the corresponding wide subcategory is the collection of 
all finite abelian groups whose torsion is contained in  $S$. It can be easily verified that the Grothendieck group of this subcategory 
is a  free abelian group of rank $|S|$. (The $p$-length of a finite abelian group, for each $p$ in $S$, gives the universal Euler 
characteristic function.) If $S = Spec(\mathbb{Z})$, then the corresponding wide subcategory is clearly the collection of all finitely 
generated abelian groups. The Grothendieck group of the latter is well-known to be infinite cyclic.  (The free rank of the finitely generated 
abelian group gives the universal Euler characteristic function.) Observe that these answers agree with our earlier 
computations (theorem \ref{thm:K_0PID}) of the Grothendieck groups for the  thick subcategories of perfect complexes over a PID.
\end{example}

So motivated by this example, we now prove that the above bijection between the lattices $\mathcal{L}_{\Wide} (R)$ 
and $\mathcal{L}_{\Thick}(R)$ preserves $K$-theory.

\begin{thm}\label{thm:K_0wide} Let $R$ be a regular coherent ring. If $\C$ is any wide subcategory of $\Wide(R)$, then
\[ K_0(\C) \cong K_0(f(\C)).\]
\end{thm}

\begin{proof} Since these Grothendieck groups are quotients of free abelian groups on the appropriate isomorphism classes
of objects, we first define maps on the free abelian groups. To this end let $FW$ denote the the free abelian group on isomorphism 
classes of objects in $\C$ and let $FT$ denote the  free abelian group on isomorphism 
classes of objects in $f(\C)$. Now we define maps 
\[  FW \stk{\phi} FT \stk{\psi} FW.  \]
Since $FW$ and $FT$ are free abelian groups, it suffices to define the maps on the basis elements. So we define $\phi(M) := M[0]$ and $\psi(x) := \sum (-1)^i H_i(x)$. 
In order to show that these maps give the isomorphism in $K$-theory, we have to check the following 
things.

\begin{enumerate}
\item Both $\phi$ and $\psi$ are well-defined group homomorphisms.
\item $\phi$ descends to $K_0(\C)$ and $\psi$ descends to $K_0(f(\C))$.
\item $\phi$ and $\psi$ are inverses to each other. 
\end{enumerate}

\noindent
1. $M$ being a finitely presented module over a regular coherent ring, it  has finite projective dimension \cite[Theorem 6.2.1]{sarah}.
This implies \cite[Prop. 3.4]{wide} that $M[0]$ is small. Further $H_*(M[0]) = M$ belongs to $\C$, so by definition
$M[0]$ is in $f(\C)$. This shows that $\phi$ is well defined. As for $\psi$, if $x \in f(\C)$, then by definition 
$H_n(x) \in \C$ for all $n$. So $\sum (-1)^i H_i(x)$ is a well defined element of $FW$.

\noindent
2. To show that $\phi$ descends, all we need to show is that whenever $0 \rightarrow A \rightarrow B \rightarrow C \rightarrow 0$
is a short exact sequence in $\C$, then $\phi(B)=\phi(A) + \phi(C)$. We can replace each of these three modules with their
projective resolutions and these resolutions will be perfect complexes because the ring is coherent (see the proof of \cite[Prop. 3.4]{wide}). 
So let $P_A$ and $P_C$ be projective resolutions of the modules $A$ and $C$. Then by the Horseshoe lemma \cite[Lemma 2.2.8]{wei} we know that 
$P_A \oplus P_C$ is a projective resolution of $B$ and these resolutions fit into a short exact sequences of complexes
$0\rar P_A \rar P_A \oplus P_C \rar P_C \rar 0$. 
This gives $[P_B] = [P_A] + [P_C]$ or equivalently $[B[0]] = [A[0]] + [C[0]]$, which translates to $f(B) = f(A) + f(C)$. 
So we have proved that $\phi$ descends to $K_0 (\C)$.

Now we show that the map $\psi$ descends to $K_0(f(\C))$. Start with an exact triangle 
$x \rar y \rar z \rar \Sigma\, x$ in $f(\C)$. We have to show that 
\[ \sum (-1)^i H_i(y) = \sum (-1)^i H_i(x) + \sum (-1)^i H_i(z)\]
in $FW$. In view of the induced long exact sequence in homology groups (in $\C$), it suffices to show that whenever
\[  0 \rar A_1  \rar  A_2\rar A_3 \rar \cdots \rar A_n \rar 0 \]
is an exact sequence in $\C$, then $\sum (-1)^i A_i = 0$ in $K_0(\C)$. When $n=3$, this holds by definition.
When $n=4$, we break the exact sequence 
\[  0 \rar A_1  \rar  A_2 \rar A_3 \stk{\alpha}   A_4 \rar 0\]
into two exact sequences in $\C$:
$0 \rar A_1 \rar A_2 \rar \Ker(\alpha) \rar 0$ and $ 0 \rar \Ker(\alpha)\rar  A_3 \rar A_4 \rar0$
($\C$ is closed under kernels, therefore $\Ker(\alpha)$ belongs to $\C$). This gives,
in $K_0(\C)$, two equations: $A_1 - A_2 + \Ker(\alpha) = 0$ and $\Ker(\alpha) - A_3 + A_4 = 0$. 
Combining these two 
equations tells us that the alternating sum of the $A_i$ is equal to zero. The general statement
now follows by a straight forward induction.

\noindent
3. Clearly $\psi \phi (M) = \psi(M[0]) = M.$ To show that the other composition is also identity, first note that
$\phi \psi (x) =\phi (\sum(-1)^i H_i(x)) = \sum (-1)^i H_i(x)[0]$. It remains to show that this
last expression is equal to $x$  in $K_0(f(\C))$.

In order to show that, we make use of the natural $t$-structure on the derived category. Let $l$ denote the
smallest integer such that $H_i(x)=0$ for all $i > l $ (homologically graded complexes). Then we get the following 
exact triangles in $D(R)$ using the truncation functors given by the $t$-structure.

\begin{center}
$x^{\le l-1} \rar x \rar x^{= l} \rar \Sigma\, x^{\le l-1}$\\
$x^{\le l-2} \rar x^{\le l-1} \rar x^{= l-1} \rar \Sigma \, x^{\le l-2} $\\
$x^{\le l-3} \rar x^{\le l-2} \rar x^{= l-2} \rar \Sigma \, x^{\le l-3}$  \hspace{3 mm} \text{etc.}
\end{center}
(Here $x^{\le i}$ is a complex whose homology agrees with that of $x$ in dimensions $\le i$ and has no homology in dimensions
$> i$.)
\noindent
Now we claim that these triangles are in fact exact triangles in $f(\C)$. First of all, the homology 
groups of these complexes belong to $\C$, and by properties of the truncation functors,
it is clear that the homology groups of each complex in the above exact triangles are concentrated only in a finite range. Now the smallness
of these complexes follows from the following characterisation of small objects over a regular coherent ring.

\begin{lemma} \cite[Note 9.8]{ch} Let $R$ be a regular coherent ring. Then a complex $x$ in $D(R)$ is small if and only if each
$H_n(x)$ is finitely presented and only finitely many are non-zero.
\end{lemma}

The above exact triangles give the following equations in $K_0 (f(\C))$. 
\begin{eqnarray*}
[x] &=&  [x^{\le l-1}] + [x^{= l}]     \\
    &=&  [x^{\le l-2}] + [x^{= l-1}] + [x^{= l}] \\
    &=&  [x^{\le l-3}] + [x^{= l-2}] + [x^{= l-1}] + [x^{= l}] \\
    &=&  \sum_{i \le l} [x^{= i}]
\end{eqnarray*}
It is easy to see that $x^{=i}$ is quasi-isomorphic to $H_i(x)[i]$.
By replacing this quasi-isomorphism in the last equation we get 
$[x] = \sum [H_i(x)[i]] =  \sum (-1)^i [H_i(x)[0]].$
That proves that $\phi \psi (x) = x$, and we are done.
\end{proof}

\begin{cor} If $R$ is a regular coherent ring, then 
\[ K_0 ( D^b (\Wide(R) )) \cong K_0 (D^b( \proj \,R)).\]
\end{cor}

\begin{proof} It is well-known that if $\A$ is any essentially small abelian category, then
$K_0(\A) \cong K_0(D^b(\A))$; see \cite[Chapter VIII, Page 357]{sga5}. In particular if $\A = \Wide(R)$ (abelian category of finitely
presented $R$-modules), then we get $K_0(\Wide(R)) \cong K_0(D^b(\Wide(R)))$.  By theorem \ref{thm:K_0wide}
we know that $K_0(\Wide(R)) \cong K_0(D^b(\proj \, R))$. So the corollary follows by combining these
two isomorphisms.
\end{proof}

\section{Krull-Schmidt decompositions for wide subcategories} \label{sec:KSwide}

We now study Krull-Schmidt decompositions for wide subcategories.
A \emph{Krull-Schmidt} decomposition of a wide subcategory $\W$  of an abelian category $\A$ is 
a collection of wide subcategories 
$(\W_i)_{i \in I}$ in $\A$ such that, 
\begin{enumerate}
\item $\W_i \bigcap \W_j = 0$ for all $i \ne j$.
\item The objects of $\W$ are the finite coproducts of objects in $\W_i$. 
\item Each $\W_i$ is indecomposable: 
$\W_i \ne 0$, and $\W_i = \A \coprod \B \Rightarrow \A = 0 \ \text{or\;}  \B = 0.$
\item Uniqueness: If $\coprod_{i \in I} \mathcal{V}_i$ is any other collection of wide subcategories that satisfy the above three properties, 
then $\W_i \cong \mathcal{V}_i$ up to a permutation of indices.
\end{enumerate}
A Krull-Schmidt decomposition is denoted by $\W = \amalg_{i \in I} \W_i$. The main theorem of this section then states:

\begin{thm} Let $R$ be a noetherian regular ring. Then every wide subcategory of $\Wide(R)$ admits a Krull-Schmidt
decomposition. Conversely, given wide subcategories  $\W_i$  of $\Wide(R)$ such that $\W_1 \bigcap \W_j = 0$ for $i \ne j$, there 
exists a unique wide subcategory $\W \subseteq \Wide(R)$ such that $\W = \coprod_{i \in I} \W_i$ is a Krull-Schmidt decomposition for $\W$.
\end{thm}

\begin{proof} 
Since $R$ is a noetherian ring, $\Wide(R)$ consists of the finitely generated $R$-modules.
Therefore let $\W$ be a wide subcategory of finitely generated $R$-modules and let $S$ denote the corresponding thick support under the bijection $\xi$ (see figure \ref{fig:wideandthick}).
Since $R$ is a noetherian ring, we know that $S$ admits a unique Krull-Schmidt decomposition $S = \bigcup_{i \in I} S_i$; see proposition \ref{prop:graph}.
Now let $\W_i$ be the wide 
subcategory corresponding to $S_i$. We claim that $\amalg_{i \in I} \W_i$ is a Krull-Schmidt decomposition of $\W$.

The first thing to note here is that $\W_i \bigcap \W_j = 0$ for $i \ne j$: Suppose not, then since the intersection of wide subcategories is a wide 
subcategory, the non-trivial wide subcategory $\W_i \bigcap \W_j$ would then corresponds to a non-empty thick support that is contained in 
$S_i \bigcap S_j$, which is impossible since $S_i \bigcap S_j = \emptyset$. Therefore  $\W_i \bigcap \W_j = 0$ for $i \ne j$. 

We now argue that the objects in $\W$ are the finite coproducts of objects in $\W_i$. We give an indirect method by showing that 
$\amalg\, \W_i$ is a wide subcategory. Once we show this, then since the thick support corresponding to this wide subcategory 
is clearly $\bigcup S_i\;( = S)$, the conclusion follows.
As a further reduction, it suffices to show that $\W_1 \amalg \W_2$ is a wide subcategory if $\W_1$ and $\W_2$ are wide subcategories corresponding
to disjoint thick supports $S_1$ and $S_2$ respectively. First observe that $\Hom(\W_1, \W_2) = 0$: To see this, consider modules $M_i \in \W_i$. 
Then for all primes $p$ in $\Spec(R)$, we have
\[{\Hom_R (M_1, M_2)}_p = \Hom_{R_p}({M_1}_p, {M_2}_{p}) = 0.\]
(The last equality follows from the fact that $S_1 \bigcap S_2 = \emptyset$.)
So $\Hom(\W_1, \W_2)=0$ and consequently every  map $M_1 \oplus M_2 \rar M_1' \oplus M_2'$ in $\W_1 \coprod \W_2$ is of the form $f \oplus f'$. 
It is now clear that  $\W_1 \amalg \W_2$ is an abelian category. The more interesting part is to show that $\W_1 \amalg \W_2$ is closed under extensions. 
Towards this, consider a short exact sequence of finitely generated modules 
\[ 0 \rar M_1 \oplus M_2 \rar M \rar M_1'\oplus M_2'\rar 0 ,\]
where $M_1\oplus M_2$ and $M_1' \oplus M_2'$ belong to $\W_1 \coprod \W_2$. This short exact sequence corresponds to a class in 
$\Ext^1_R (M_1' \oplus M_2', M_1\oplus M_2)$. Since $\Ext$ behaves well with respect to finite coproducts in both components, we get 
\[\Ext^1_R (M_1' \oplus M_2', M_1\oplus M_2) \cong \Ext^1_R(M_1', M_1) \oplus \Ext^1_R(M_1', M_2) \oplus \Ext^1_R(M_2', M_1)\oplus \Ext^1_R(M_2' ,M_2).  \]
We now claim that for $i \ne j$, $\Ext^1_R(M_i', M_j) = 0$. To see this, first observe that the supports of $M_i'$ and $M_j$ are disjoint. Therefore for each $p$
we have
\[ \Ext^1_R {(M_i', M_j)}_p = \Ext^1_{R_p} ({M_i'}_p, {M_j}_p)  = 0. \]
(Note that we used the fact that localisation commutes with the $\Ext$ functor for finitely generated modules over noetherian rings.)
Since all the localisations are zero, the module $\Ext^1_R(M_i', M_j)$ is itself zero.
The last equation of $\Ext$ groups now simplifies to 
\[\Ext^1_R (M_1' \oplus M_2', M_1\oplus M_2) = \Ext^1_R(M_1', M_1) + \Ext^1_R(M_2' ,M_2) .  \]
So the upshot is that every extension of $M_1 \oplus M_2$ and $M_1' \oplus M_2'$ is a sum of two extensions: one obtained from an extension 
of $M_1$ and $M_1'$ and the other from $M_2$ and $M_2'$. This shows that  $\W_1 \coprod \W_2$ is also closed under extensions. 
That completes the proof of the existence of a  Krull-Schmidt decomposition.

We leave the uniqueness part and the converse as easy exercises to the reader.

\end{proof}

\begin{rem} The above can be compared with the analogous Krull-Schmidt theorem \ref{th:mainalg} for the derived category.
Again, it is worth noting that our proofs give the following  stronger statements.

\begin{enumerate}
\item Let $R$ be a regular coherent ring. Then a wide subcategory of finitely presented modules
admits a Krull-Schmidt decomposition if and only if the corresponding thick support admits a Krull-Schmidt decomposition.
\item Let $R$ be a regular coherent ring with the following properties.

(a) Every open subset of $\Spec(R)$ is compact.

(b) $\Spec(R)$ satisfies the descending chain condition.

Then every wide subcategory of finitely presented modules admits a Krull-Schmidt decomposition.
\end{enumerate}
\end{rem}

 The following splitting result  which is a consequence of these decompositions is well-known.

\begin{cor} \label{cor:splittingformodules} Let $R$ be a noetherian ring and $M$ a finitely generated $R$-module. Then $M$
admits a unique splitting
\[M \cong \bigoplus_{i = 1}^n \; M_i, \]
for some $n$ such that the supports of the $M_i$ are pairwise disjoint and indecomposable.
\end{cor}

\begin{proof} Existence of such a splitting:  Since $M$ is a finitely generated module, $\Supp(M)$
is a closed set and therefore admits a Krull-Schmidt decomposition. Let $\mathcal{L}$ denote the smallest wide subcategory of $R$-modules
closed under direct sums that contains $M$. Hovey \cite[Theorem 5.2]{wide} showed that such subcategories are in bijection with subsets of $\Spec(R)$.
It can be shown (exactly as shown above for wide subcategories) that  $\mathcal{L}$ admits a Krull-Schmidt decomposition $\coprod_{i \in I} \mathcal{L}_i$
which corresponds to a Krull-Schmidt decomposition of $\Supp(M)$. The existence of the desired splitting of $M$ is now obtained from this decomposition 
of $\mathcal{L}$. Uniqueness follows from the fact that there are no non-trivial maps from $M$ to $N$ whenever $M$ and $N$ are finitely generated modules 
with disjoint supports.
\end{proof}

\begin{cor} Let $R$ be a noetherian regular ring and let $\W$ be a wide subcategory of finitely generated $R$ modules. Then
\[K_0(\W) \cong \underset{i \in I}{\bigoplus} \; K_0(\W_i),\]
where $\coprod_{i \in I} \W_i$ is the Krull-Schmidt decomposition of $\W$.
\end{cor}

\begin{proof} One can give a  direct proof of this corollary but we prefer to use our earlier results on Krull-Schmidt decompositions for 
thick subcategories. Let $S$ be the thick support corresponding to $\W$ and let $\bigcup_{i \in I} S_i$ be the Krull-Schmidt 
decomposition for $S$ (proposition \ref{prop:graph}). Now $f(\W)$ is a thick subcategory of perfect complexes and therefore it admits a unique Krull-Schmidt 
decomposition (theorem \ref{th:mainalg}) $\coprod \T_i$ corresponding to the decomposition $\bigcup_{i \in I} S_i$ (in fact $\T_i = f(\W_i)$). 
We have also seen (corollary \ref{cor:Ktheory}) that $K_0(f(\W)) \cong \bigoplus_{i \in I} K_0(f(W_i))$.  The main theorem of the previous section tells that
$K_0(\W_i) = K_0(f(\W_i)).$ By combining all these facts, we get the desired isomorphism. 
\end{proof}

\section{Questions}

Note that our commutative ring $R$ was assumed to be regular for most of the results in this chapter. That seems to be somewhat 
unnatural restriction on the ring for these results. At least some of these results ought to be true without this regularity condition.  
In order to address this issue, we first have to find out whether Hovey's classification of 
wide subcategories holds for all  coherent rings. The classification of wide subcategories is completely a module theoretic question. 
However, Hovey's proof of this classification detours into the derived category and makes use of Thomason's classification of thick 
subcategories. It is precisely in this passage, from the module category to the derived category, the regularity assumption came in.
(Basically the point is that we are replacing a finitely presented module by its projective resolution. Unless the ring is regular,
there is no guarantee that this resolution will land in the category of perfect complexes.) So maybe one should try to give a direct proof 
of Hovey's classification. Once the regularity assumption can be relaxed from his theorem,
we can hope to relax it from the other results in this chapter. 


\printendnotes

%
%
\bibliographystyle{alpha}
\bibliography{uwthesis}

\begin{thebibliography}{ATJLSS04}

\bibitem[Ada69]{adams}
J.~F. Adams.
\newblock Lectures on generalised cohomology.
\newblock In {\em Category Theory, Homology Theory and their Applications, III
  (Battelle Institute Conference, Seattle, Wash., 1968, Vol. Three)}, pages
  1--138. Springer, Berlin, 1969.

\bibitem[Ada92]{MayTho}
J.~F. Adams.
\newblock {\em The selected works of {J}. {F}rank {A}dams. {V}ol. {II}}.
\newblock Cambridge University Press, Cambridge, 1992.

\bibitem[Ada95]{bluebook}
J.~F. Adams.
\newblock {\em Stable homotopy and generalised homology}.
\newblock Chicago Lectures in Mathematics. University of Chicago Press,
  Chicago, IL, 1995.

\bibitem[AK89]{adamskuhn}
J.~F. Adams and N.~J. Kuhn.
\newblock Atomic spaces and spectra.
\newblock {\em Proc. Edinburgh Math. Soc. (2)}, 32(3):473--481, 1989.

\bibitem[AM69]{am}
M.~F. Atiyah and I.~G. Macdonald.
\newblock {\em Introduction to commutative algebra}.
\newblock Addison-Wesley Publishing Co., Reading, Mass.-London-Don Mills, Ont.,
  1969.

\bibitem[ATJLSS04]{AlLeJe}
Leovigildo Alonso~Tarr{\'{\i}}o, Ana Jerem{\'{\i}}as~L{\'o}pez, and
  Mar{\'{\i}}a~Jos{\'e} Souto~Salorio.
\newblock Bousfield localization on formal schemes.
\newblock {\em J. Algebra}, 278(2):585--610, 2004.

\bibitem[BCR97]{bcr}
D.~J. Benson, Jon~F. Carlson, and Jeremy Rickard.
\newblock Thick subcategories of the stable module category.
\newblock {\em Fund. Math.}, 153(1):59--80, 1997.

\bibitem[Ben98]{be}
D.~J. Benson.
\newblock {\em Representations and cohomology. {I}}, volume~30 of {\em
  Cambridge Studies in Advanced Mathematics}.
\newblock Cambridge University Press, Cambridge, 1998.

\bibitem[BH93]{cm}
Winfried Bruns and J{\"u}rgen Herzog.
\newblock {\em Cohen-{M}acaulay rings}, volume~39 of {\em Cambridge Studies in
  Advanced Mathematics}.
\newblock Cambridge University Press, Cambridge, 1993.

\bibitem[Boa99]{boardman}
J.~Michael Boardman.
\newblock Conditionally convergent spectral sequences.
\newblock In {\em Homotopy invariant algebraic structures (Baltimore, MD,
  1998)}, volume 239 of {\em Contemp. Math.}, pages 49--84. Amer. Math. Soc.,
  Providence, RI, 1999.

\bibitem[Bou79]{bob}
A.~K. Bousfield.
\newblock The {B}oolean algebra of spectra.
\newblock {\em Comment. Math. Helv.}, 54(3):368--377, 1979.

\bibitem[Car81]{car}
Jon~F. Carlson.
\newblock Complexity and {K}rull dimension.
\newblock In {\em Representations of algebras (Puebla, 1980)}, volume 903 of
  {\em Lecture Notes in Math.}, pages 62--67. Springer, Berlin, 1981.

\bibitem[Chr98]{ch}
J.~Daniel Christensen.
\newblock Ideals in triangulated categories: phantoms, ghosts and skeleta.
\newblock {\em Adv. Math.}, 136(2):284--339, 1998.

\bibitem[CKN01]{ChKeNe}
J.~Daniel Christensen, Bernhard Keller, and Amnon Neeman.
\newblock Failure of {B}rown representability in derived categories.
\newblock {\em Topology}, 40(6):1339--1361, 2001.

\bibitem[DHS88]{dhs}
Ethan~S. Devinatz, Michael~J. Hopkins, and Jeffrey~H. Smith.
\newblock Nilpotence and stable homotopy theory. {I}.
\newblock {\em Ann. of Math. (2)}, 128(2):207--241, 1988.

\bibitem[DP01]{dwypal}
William~G. Dwyer and John~H. Palmieri.
\newblock Ohkawa's theorem: there is a set of {B}ousfield classes.
\newblock {\em Proc. Amer. Math. Soc.}, 129(3):881--886, 2001.

\bibitem[EKMM97]{ekmm}
A.~D. Elmendorf, I.~Kriz, M.~A. Mandell, and J.~P. May.
\newblock {\em Rings, modules, and algebras in stable homotopy theory},
  volume~47 of {\em Mathematical Surveys and Monographs}.
\newblock American Mathematical Society, Providence, RI, 1997.
\newblock With an appendix by M. Cole.

\bibitem[Fau03]{fau}
H.~Fausk.
\newblock Picard groups of derived categories.
\newblock {\em J. Pure Appl. Algebra}, 180(3):251--261, 2003.

\bibitem[FP04]{FriPev}
Eric~M. Friedlander and Julia Pevtsova.
\newblock $\pi$-supports for modules for finite group schemes over a field.
\newblock {\em Preprint}, 2004.

\bibitem[FS97]{fs}
Eric~M. Friedlander and Andrei Suslin.
\newblock Cohomology of finite group schemes over a field.
\newblock {\em Invent. Math.}, 127(2):209--270, 1997.

\bibitem[Gla92]{sarah}
Sarah Glaz.
\newblock Commutative coherent rings: historical perspective and current
  developments.
\newblock {\em Nieuw Arch. Wisk. (4)}, 10(1-2):37--56, 1992.

\bibitem[Gro77]{sga5}
A~Grothendieck.
\newblock {\em Cohomologie L-adique et Fonctions L}.
\newblock Springer-Verlag, 1977.

\bibitem[Har77]{Ha}
Robin Hartshorne.
\newblock {\em Algebraic geometry}.
\newblock Springer-Verlag, New York, 1977.

\bibitem[Hop87]{Ho}
Michael~J. Hopkins.
\newblock Global methods in homotopy theory.
\newblock In {\em Homotopy theory (Durham, 1985)}, volume 117 of {\em London
  Math. Soc. Lecture Note Ser.}, pages 73--96. Cambridge Univ. Press,
  Cambridge, 1987.

\bibitem[Hov99]{hoveymodel}
Mark Hovey.
\newblock {\em Model categories}, volume~63 of {\em Mathematical Surveys and
  Monographs}.
\newblock American Mathematical Society, Providence, RI, 1999.

\bibitem[Hov01]{wide}
Mark Hovey.
\newblock Classifying subcategories of modules.
\newblock {\em Trans. Amer. Math. Soc.}, 353(8):3181--3191 (electronic), 2001.

\bibitem[HP00]{hp2}
Mark Hovey and John~H. Palmieri.
\newblock Galois theory of thick subcategories in modular representation
  theory.
\newblock {\em J. Algebra}, 230(2):713--729, 2000.

\bibitem[HP01]{hp}
Mark Hovey and John~H. Palmieri.
\newblock Stably thick subcategories of modules over {H}opf algebras.
\newblock {\em Math. Proc. Cambridge Philos. Soc.}, 130(3):441--474, 2001.

\bibitem[HPS97]{mps}
Mark Hovey, John~H. Palmieri, and Neil~P. Strickland.
\newblock Axiomatic stable homotopy theory.
\newblock {\em Mem. Amer. Math. Soc.}, 128(610):x+114, 1997.

\bibitem[HS98]{hs}
Michael~J. Hopkins and Jeffrey~H. Smith.
\newblock Nilpotence and stable homotopy theory. {II}.
\newblock {\em Ann. of Math. (2)}, 148(1):1--49, 1998.

\bibitem[Kra99]{Kr}
Henning Krause.
\newblock Decomposing thick subcategories of the stable module category.
\newblock {\em Math. Ann.}, 313(1):95--108, 1999.

\bibitem[Kra02]{kr1}
Henning Krause.
\newblock A {B}rown representability theorem via coherent functors.
\newblock {\em Topology}, 41(4):853--861, 2002.

\bibitem[Lan91]{la}
Steven~E. Landsburg.
\newblock {$K$}-theory and patching for categories of complexes.
\newblock {\em Duke Math. J.}, 62(2):359--384, 1991.

\bibitem[Mar83]{mar}
H.~R. Margolis.
\newblock {\em Spectra and the {S}teenrod algebra}, volume~29 of {\em
  North-Holland Mathematical Library}.
\newblock North-Holland Publishing Co., Amsterdam, 1983.

\bibitem[Mit85]{mit}
Stephen~A. Mitchell.
\newblock Finite complexes with {$A(n)$}-free cohomology.
\newblock {\em Topology}, 24(2):227--246, 1985.

\bibitem[Nee92]{Ne}
Amnon Neeman.
\newblock The chromatic tower for {$D(R)$}.
\newblock {\em Topology}, 31(3):519--532, 1992.

\bibitem[Nee96]{ne5}
Amnon Neeman.
\newblock The {G}rothendieck duality theorem via {B}ousfield's techniques and
  {B}rown representability.
\newblock {\em J. Amer. Math. Soc.}, 9(1):205--236, 1996.

\bibitem[Nee00]{ne2}
Amnon Neeman.
\newblock Oddball {B}ousfield classes.
\newblock {\em Topology}, 39(5):931--935, 2000.

\bibitem[Nee01]{ne1}
Amnon Neeman.
\newblock {\em Triangulated categories}.
\newblock Princeton University Press, Princeton, NJ, 2001.

\bibitem[Pup62]{puppe}
D~Puppe.
\newblock On the formal structure of stable homotopy theory.
\newblock {\em Colloq. Alg. Topology}, pages 65--71, 1962.

\bibitem[Rav84]{rav}
Douglas~C. Ravenel.
\newblock Localization with respect to certain periodic homology theories.
\newblock {\em Amer. J. Math.}, 106(2):351--414, 1984.

\bibitem[Ric00]{Ri1}
Jeremy Rickard.
\newblock Bousfield localization for representation theorists.
\newblock In {\em Infinite length modules (Bielefeld, 1998)}, Trends Math.,
  pages 273--283. Birkh\"auser, Basel, 2000.

\bibitem[Tho97]{Th}
R.~W. Thomason.
\newblock The classification of triangulated subcategories.
\newblock {\em Compositio Math.}, 105(1):1--27, 1997.

\bibitem[Wei94]{wei}
Charles~A. Weibel.
\newblock An introduction to homological algebra.
\newblock 38:xiv+450, 1994.

\bibitem[Wei03]{wei1}
Charles~A. Weibel.
\newblock {\em Algebraic K-theory}.
\newblock http://math.rutgers.edu/\~{}weibel, 2003.

\bibitem[Wil81]{wil}
Clarence Wilkerson.
\newblock {\em The cohomology algebras of finite-dimensional {H}opf algebras},
  volume 264.
\newblock 1981.

\bibitem[Xu95]{xu}
Kai Xu.
\newblock A note on atomic spectra.
\newblock In {\em The \v Cech centennial (Boston, MA, 1993)}, volume 181 of
  {\em Contemp. Math.}, pages 419--422. Amer. Math. Soc., Providence, RI, 1995.

\end{thebibliography}

%
%
\appendix
\raggedbottom\sloppy
 

%
%

\vita{
Sunil Chebolu was born in the Visakhapatnam district of south India on May 11, 1978. As an undergraduate he attended the Indian Statistical Institute 
where he received a Bachelor of Statistics degree (B.Stat) in 1998 and a Master of Statistics degree (M.Stat) in 2000. He moved to the U.S. in the Fall of
2000 to pursue graduate studies in Mathematics. He received a Master of Science degree in 2002 
and a Doctor of Philosophy in 2005, both from the University of Washington.}


\end{document}